\documentclass[A4j,11pt]{article}
\usepackage{amsfonts,amsmath,amscd}
\usepackage{graphics}
\usepackage{amssymb}
\begin{document}

\title{{\bf Collapse of the mean curvature flow\\
for isoparametric submanifolds\\
in non-compact symmetric spaces}}
\author{{\bf Naoyuki Koike}}
\date{}
\maketitle

\begin{abstract}
It is known that principal orbits of Hermann actions on a symmetric space 
of non-compact type are curvature-adapted isoparametric submanifolds 
having no focal point of non-Euclidean type on the ideal boundary of the ambient 
symmetric space.  
In this paper, we investigate the mean curvature flows for such a 
curvature-adapted isoparametric submanifold and its focal submanifold.  
Concretely the investigation is performed by investigating the mean curvature 
flows for the lift of the submanifold to an infinite dimensional 
pseudo-Hilbert space through a pseudo-Riemannian submersion.  
\end{abstract}

\section{Introduction}
Let $f_t$'s ($t\in[0,T)$) be a one-parameter $C^{\infty}$-family of immersions 
of a manifold $M$ into a Riemannian manifold $N$, where $T$ is a positive 
constant or $T=\infty$.  Define a map 
$\widetilde f:M\times[0,T)\to N$ by $\widetilde f(x,t)=f_t(x)$ 
($(x,t)\in M\times[0,T)$).  If, for each $t\in[0,T)$,  
$\widetilde f_{\ast}((\frac{\partial}{\partial t})_{(x,t)})$ is the mean 
curvature vector of $f_t:M\hookrightarrow N$, then $f_t$'s ($t\in[0,T)$) is 
called a mean curvature flow.  
In particular, if $f_t$'s are embeddings, then we call $M_t:=f_t(M)$'s 
$(0\in[0,T))$ rather than $f_t$'s $(0\in[0,T))$ a mean curvature flow.  
Liu-Terng [LT] investigated the mean curvature flows for an isoparametric 
submanifold in a Euclidean space and its focal submanifold 
and obtained the following facts.  

\vspace{0.5truecm}


\noindent
$\overline{\qquad\qquad\qquad\qquad\qquad\qquad}$

\vspace{0.1truecm}

\noindent
{\footnotesize 2010 {\sl Mathematics Subject classification.} Primary 53C44; Secondly 53C35.

\noindent
{\sl Key words and phrases.} mean curvature flow, isoparametric submanifold, symmetric space.}

\newpage

\noindent
{\bf Fact 1([LT]).} {\sl Let $M$ be a compact isoparametric submanifold in a 
Euclidean space.  Then the following statements ${\rm (i)}$ and ${\rm (ii)}$ hold:

${\rm (i)}$ The mean curvature flow $M_t$ for $M$ collapses 
to a focal submanifold $F$ of $M$ in finite time.  
If the natural fibration of $M$ onto $F$ 
is spherical, then the mean curvature flow $M_t$ has type I singularity, 
that is, $\lim\limits_{t\to T-0}
{\rm max}_{v\in S^{\perp}M_t}\vert\vert A^t_v\vert\vert^2_{\infty}(T-t)$ 
$<\infty$, where $A^t_v$ is the shape operator of $M_t$ for $v$, 
$\vert\vert A^t_v\vert\vert_{\infty}$ is the sup norm of $A^t_v$ and 
$S^{\perp}M_t$ is the unit normal bundle of $M_t$.

{\rm(ii)} For any focal submanifold $F$ of $M$, the set of all parallel 
submanifolds of $M$ collapsing to $F$ along the mean curvature flow 
is a one-parameter $C^{\infty}$-family.}

\vspace{0.5truecm}

\noindent
{\bf Fact 2([LT]).} {\sl Let $M$ be as in Fact 1, $C$ the Weyl domain of $M$ at $x_0\,(\in M)$ and $\sigma$ a 
simplex of dimension greater than zero of $\partial C$.  
Then the following statements ${\rm (i)}$ and ${\rm (ii)}$ hold:

{\rm(i)} For any focal submanifold $F$ (of $M$) through $\sigma$, the maen 
curvature flow $F_t$ for $F$ collapses to a focal submanifold $F'$ (of $M$) 
through $\partial{\sigma}$ in finite time.  
If the natural fibration of $F$ onto $F'$ is spherical, then the mean 
curvature flow $F_t$ has type I singularity.

(ii) For any focal submanifold $F$ (of $M$) through $\partial\sigma$, the set 
of all focal submanifolds of $M$ through $\sigma$ collapsing to $F$ along the 
mean curvature flow is a one-parameter $C^{\infty}$-family.}

\vspace{0.5truecm}

Since the focal submanifold of $M$ through the only $0$-dimensional simplex of 
$\partial C$ is a one-point set, it follows from the statement (i) of Facts 1 
and 2 that $M$ collapses to a one-point set after finitely many times of 
collapsings along the mean curvature flows.  

As a generalized notion of compact isoparametric 
hypersurfaces in a sphere and a hyperbolic space, and a compact isoparametric 
submanifolds in a Euclidean space, Terng-Thorbergsson [TT] introduced 
the notion of an equifocal submanifold in a symmetric space $G/K$.  
This notion is defined as a compact submanifold (which we denote by $M$) in 
$G/K$ with flat section, trivial normal holonomy group and parallel focal 
structure.  Here the parallel focal structure means that the 
tangential focal structures of $M$ move to one another under the parallel 
translations with respect to the normal connection of $M$.  
For a compact submanifold $M$ with flat section and trivial normal 
holonomy group, it is equifocal if and only if, for any parallel normal vector 
field $\widetilde v$ of $M$, the set of all the focal radii of $M$ along 
the normal geodesic $\gamma_{\widetilde v_x}$ with 
$\gamma'_{\widetilde v_x}(0)=\widetilde v_x$ is independent of the choice of 
$x\in M$.  
On the other hand, Heintze-Liu-Olmos [HLO] introduced the notion of 
isoparametric submanifold with flat section in a (general) Riemannian 
manifold as a (properly embedded) submanifold with flat section and 
trivial normal holonomy group whose sufficiently close parallel submanifolds 
are of constant mean curvature with respect to the radial direction.  
In the sequel, we assume that all isoparametric submanifolds with flat section are 
complete.  

\vspace{0.4truecm}

\noindent
{\bf Terminology.} In this paper, we shall call an isoparametric submanifold 
with flat section an {\it isoparametric submanifold} simply.  

\vspace{0.4truecm}

\noindent
For a compact submanifold in a symmetric space $G/K$ of compact type, 
the isoparametricness is equivalent to the equifocality (see [HLO]).  
The author has recently investigated the mean curvature flows for an equifocal 
submanifold in a symmetric space of compact type and its focal submanifold, 
and obtained the following facts.  

\vspace{0.5truecm}

\noindent
{\bf Fact 3([Koi10]).} 
{\sl Let $M$ be an equifocal submanifold in a symmetric space 
$G/K$ of compact type.  
Then the following statements ${\rm (i)}$ and ${\rm (ii)}$ hold:

${\rm (i)}$ If $M$ is not minimal, then the mean curvature flow $M_t$ for $M$ 
collapses to a focal submanifold $F$ of $M$ in finite 
time.  Furthermore, if $M$ is irreducible, if the codimension of $M$ is geater 
than one and if the natural fibration of $M$ onto $F$ is spherical, then 
$M_t$ has type I singularity.  

${\rm (ii)}$ For any focal submanifold $F$ of $M$, 
the set of all parallel submanifolds of $M$ collapsing to $F$ along the mean 
curvature flow is a one-parameter $C^{\infty}$-family.}

\vspace{0.5truecm}

\noindent
{\bf Fact 4([Koi10]).} 
{\sl Let $M$ be as in Fact 3, $C$ the image of the fundamental domain of 
the Coxeter group of $M$ at $x_0\,(\in M)$ by the normal exponential map 
and $\sigma$ a stratum of dimension greater than zero of $\partial C$ 
(which is a stratified space).  
Then the following statements ${\rm (i)}$ and ${\rm (ii)}$ hold:

${\rm (i)}$ For any non-minimal focal submanifold $F$ of $M$ through $\sigma$, 
the mean curvature flow $F_t$ for $F$ collapses to a focal submanifold $F'$ 
of $M$ through $\partial\sigma$ in finite time.  
If $M$ is irreducible, if the codimension of $M$ is greater than one and if 
the natural fibration of $F$ onto $F'$ is spherical, then the mean curvature 
flow $F_t$ has type I singularity.

${\rm (ii)}$ For any focal submanifold $F$ of $M$ through $\partial\sigma$, 
the set of all focal submanifolds of $M$ through $\sigma$ collapsing to $F$ 
along the mean curvature flow is a one-parameter $C^{\infty}$-family.}

\vspace{0.5truecm}

Since focal submanifolds of $M$ through the lowest dimensional stratum of 
$\partial C$ are minimal, it follows from the statement 
(i) of Facts 3 and 4 that $M$ collapses to a minimal focal submanifold of $M$ 
after finitely many times of collapsings along the mean curvature flows.  

\vspace{0.5truecm}

\noindent
{\bf Assumption.} In the sequel, we assume that all submanifolds are 
real anlaytic.  

\vspace{0.5truecm}

\noindent
We [Koi1,2] introduced the notion of a complex equifocal submanifold 
as a (properly embedded) complete submanifold with flat section, 
trivial normal holonomy group and parallel complex focal structure, 
where the parallel complex focal structure means that 
the tangential focal structures of the complexification 
$M^{\bf c}(\subset G^{\bf c}/K^{\bf c})$ of $M(\subset G/K)$ 
move to one another under the parallel translations with respect to 
the normal connection of $M^{\bf c}$.  
For a submanifold $M$ with flat section and trivial normal 
holonomy group, it is complex equifocal if and only if, for any parallel 
normal vector field $\widetilde v$ of $M$, the set of all the complex focal 
radii of $M$ along the normal geodesic $\gamma_{\widetilde v_x}$ with 
$\gamma'_{\widetilde v_x}(0)=\widetilde v_x$ is independent of the choice of 
$x\in M$.  

We ([Koi3]) introduced the notion of a proper complex equifocal submanifold 
as a complex equifocal submanifold having a good complex focal structure, 
where "good complex focal structure" means that the focal structure of 
the complexification of the submanifold at any point $x_0$ consists of 
infinitely many complex hyperplanes in the normal space at $x_0$ and that 
the group generated by the complex reflections of order two with respect to 
the complex hyperplanes is discrete.  

Next we recall the notion of a focal point of non-Euclidean type on the ideal 
boundary $N(\infty)$ of a submanifold $M$ in a Hadamard manifold $N$ 
which was introduced in [Koi6].  
Let $v$ be a unit normal vector of $M$ and 
$\gamma_v:[0,\infty)\to N$ the normal geodesic 
of $M$ of direction $v$.  If there exists a $M$-Jacobi field 
$Y$ along $\gamma_v$ satisfying 
$\lim\limits_{t\to\infty}\frac{\vert\vert Y(t)\vert\vert}{t}=0$, then we call 
$\gamma_v(\infty)\,(\in N(\infty))$ 
a {\it focal point} of $M$ {\it on the ideal boundary} $N(\infty)$ 
{\it along} $\gamma_v$, 
where $\gamma_v(\infty)$ is the asymptotic class of $\gamma_v$.  
Also, if there exists a $M$-Jacobi field $Y$ along $\gamma_v$ satisfying 
$\lim\limits_{t\to\infty}\frac{\vert\vert Y(t)\vert\vert}{t}=0$ and 
${\rm Sec}(v,Y(0))\not=0$, then we call $\gamma_v(\infty)$ a 
{\it focal point of non-Euclidean type of} $M$ {\it on} $N(\infty)$ 
{\it along} $\gamma_v$, where ${\rm Sec}(v,Y(0))$ is the 
sectional curvature for the $2$-plane spanned by $v$ and $Y(0)$.  
If, for any unit normal vector $v$ of $M$, 
$\gamma_v(\infty)$ is not a focal point of 
non-Euclidean type of $M$ on $N(\infty)$, then we say 
that $M$ {\it has no focal point of non-Euclidean type on the ideal boundary} 
$N(\infty)$.  
It is known that principal orbits of Hermann actions on symmetric spaces of 
non-compact type are curvature-adapted isoparametric submanifolds and 
they have no focal point of non-Euclidean type on the ideal boundary 
(see Theorem B in [Koi3] and its proof and so on).  
According to Theorem 15 in [Koi2] and Theorem A in [Koi6], we have the 
following fact.  

\vspace{0.5truecm}

\noindent
{\bf Fact 5.} {\sl For a curvature-adapted isoparametric submanifold $M$ in 
a symmetric space $N$ of non-compact type, it has no focal point of 
non-Euclidean type on the ideal boundary $N(\infty)$ if and only if it is proper complex 
equifocal.}

\vspace{0.5truecm}


\noindent
Let $M$ be a curvature-adapted isoparametric submanifold in a symmetric space $N=G/K$ 
of non-compact type having no focal point of non-Euclidean type on the ideal boundary 
$N(\infty)$ of $N$.  Assume that a focal submanifold of $M$ exists.  
Note that a focal submanifold of $M$ exists if $G/K$ is other than a hyperbolic space.  
Let $F_{\it l}$ be one of the lowest dimensional focal submanifolds of $M$.  
Without loss of generality, 
we may assume that $eK\in F_{\it l}$.  
Note that $F_{\it l}$ passes through 
$\exp^{\perp}(\widetilde{\sigma})$ for one $\widetilde{\sigma}$ of the lowest 
dimensional simplex of the boundary $\partial\widetilde C$ of the fundamental 
domain $\widetilde C$ of the real Coxeter group associated with $M$.  
See the next paragraph about the definition of the real Coxeter group 
associated with $M$.  
Set $\mathfrak p:=T_{eK}(G/K)$ and $\mathfrak p':=T_{eK}^{\perp}F_{\it l}$.  
Take a maximal abelian subspace $\mathfrak b$ of $\mathfrak p'$ and 
a maximal abelian subspace $\mathfrak a$ of $\mathfrak p$ containing 
$\mathfrak b$.  Let $\triangle$ be the root system of $G/K$ with respect to 
$\mathfrak a$ and $\triangle'$ be that of $F_{\it l}^{\perp}$ with respect to 
$\mathfrak b$.  Also, let $\mathfrak p_{\alpha}$ be the root space for 
$\alpha\in\triangle$.  Note that, if ${\rm rank}\,F_{\it l}^{\perp}={\rm rank}(G/K)$, 
then we have $\mathfrak a=\mathfrak b$ and $\triangle'\subset\triangle$.  
Since $M$ is curvature-adapted, so is also $F_{\it l}$.  
Hence we have 
$\mathfrak p'=\sum_{\alpha\in\triangle_+}(\mathfrak p_{\alpha}\cap\mathfrak p')$, 
where $\triangle_+$ is the positive root system of $\triangle$ with respect to a 
lexicographic ordering of $\mathfrak b^{\ast}$.  
In this paper, we prove the following fact for the mean curvature flows 
for a curvature-adapted isoparametric submanifold in a symmetric space $N=G/K$ of 
non-compact type having no focal point of non-Euclidean type on the ideal boundary 
$N(\infty)$ and its focal submanifold.  

\vspace{0.5truecm}

\noindent
{\bf Theorem A.} 
{\sl Let $M$ be a curvature-adapted isoparametric submanifold in 
a symmetric space $N=G/K$ of non-compact type 
having no focal point of non-Euclidean type on the ideal boundary 
$N(\infty)$, $M_t$ ($0\leq t<T$) the mean curvature flow for $M$, 
$\triangle,\,\mathfrak p_{\alpha}$ and $\mathfrak p'$ be as above.  
Assume that ${\rm codim}\,M={\rm rank}\,N$ and that 
${\rm dim}(\mathfrak p_{\alpha}\cap\mathfrak p')\geq\frac12{\rm dim}\,
\mathfrak p_{\alpha}\,(\alpha\in\triangle)$.  
Then the following statements {\rm(i)}, {\rm(ii)} and {\rm (iii)} hold.  

{\rm (i)} $M$ is not minimal and $M_t$ collapses to a focal submanifold of $M$ 
in finite time.  

{\rm (ii)} If $M_t$ collapses to a focal submanifold $F$ of $M$ in finite time 
and if the natural fibration of $M$ onto $F$ is spherical, then 
$M_t$ has type I singularity.  

{\rm (iii)} For any focal submanifold $F$ of $M$, 
the set of all parallel submanifolds of $M$ collapsing 
to $F$ along the mean curvature flow is a one-parameter $C^{\infty}$-family.}

\vspace{0.5truecm}

\noindent
{\it Remark 1.1.} 
The principal orbits of the isotropy action (of a symmetric space of non-compact type) 
and Hermann actions in Table 1 (see Section 5) satisfy all the conditions in Theorem A.  

\vspace{0.5truecm}

The focal set of a curvature-adapted proper complex equifocal submanifold $M$ 
at any point $x(\in M)$ consists of the images of finitely many (real) 
hyperplanes in the normal space $T^{\perp}_xM$ by the normal exponential map 
$\exp^{\perp}$ of $M$ and the group generated by the reflections with 
respect to the hyperplanes is a (finite) Coxeter group.  In [Koi6], we called 
this group the {\it real Coxeter group associated with} $M$.  

\vspace{0.5truecm}

\noindent
{\bf Theorem B.} 
{\sl Let $M$ be a curvature-adapted isoparametric submanifold in 
a symmetric space $N=G/K$ of non-compact type 
having no focal point of non-Euclidean type on the ideal boundary $N(\infty)$ 
and $M_t$ ($0\leq t<T$) the mean curvature flow for $M$.  
Assume that ${\rm codim}\,M={\rm rank}\,N$ and that the lowest dimensional 
focal submanifold of $M$ is an one-point set.  
Let $\widetilde{\sigma}$ be a stratum of dimension greater 
than zero of the fundamental domain $\widetilde C$ (which is a stratified space) 
of the real Coxeter group of $M$.  
Then, the following statements {\rm(i)} and {\rm (ii)} hold.  

${\rm (i)}$ Any focal submanifold $F$ of $M$ through 
$\exp^{\perp}(\widetilde{\sigma})$ is not minimal and the mean curvature flow 
$F_t$ for $F$ collapses to a focal submanifold $F'$ of $M$ 
through $\exp^{\perp}(\partial\widetilde{\sigma})$ in finite time.  
If the natural fibration of $F$ onto $F'$ is spherical, then 
$F_t$ has type I singularity.  

${\rm (ii)}$ For any focal submanifold $F$ of $M$ through 
$\exp^{\perp}(\partial\widetilde{\sigma})$, 
the set of all focal submanifolds of $M$ through 
$\exp^{\perp}(\widetilde{\sigma})$ collapsing to $F$ 
along the mean curvature flow is a one-parameter $C^{\infty}$-family.}

\vspace{0.5truecm}

According to the statement (i) of Theorems A and B, if $M$ is a 
curvature-adapted isoparametric submanifold 
having no focal point of non-Euclidean type on the ideal boundary, 
if ${\rm codim}\,M={\rm rank}(G/K)$ and if $F_{\it l}$ is 
one-point set, then $M$ collapses to one-point set after finitely many times 
of collapsings along the mean curvature flows.  
$$\begin{array}{c}
\displaystyle{
\begin{array}{llll}
\displaystyle{M_t\mathop{\longrightarrow}_{(t\to T_1)}}&
\displaystyle{F^1}&&\\
&\displaystyle{F^1_t\mathop{\longrightarrow}_{(t\to T_2)}F^2}&&\\
&&\displaystyle{\ddots}&\\
&&&\displaystyle{F^{k-1}_t\mathop{\longrightarrow}_{(t\to T_k)}
{\rm one}\,\,{\rm point}\,\,{\rm set}}
\end{array}
}\\
\displaystyle{
\left(
\begin{array}{l}
\displaystyle{F^1\,:\,{\rm a}\,\,{\rm focal}\,\,{\rm submanifold}\,\,
{\rm of}\,\,M}\\
\displaystyle{F^i\,:\,{\rm a}\,\,{\rm focal}\,\,{\rm submanifold}\,\,
{\rm of}\,\,F^{i-1}\,\,(i=2,\cdots,k-1)}
\end{array}
\right)}
\end{array}$$










\section{Basic notions and facts}
In this section, we briefly recall the notions of a proper complex 
equifocal submanifold in a symmetric space $G/K$ of non-compact type and 
a proper complex isoparametric submanifold in an 
(infinite dimensional) pseudo-Hilbert space.  
First we recall the notion of a complex equifocal submanifold in $G/K$.  
Let $M$ be a submanifold with flat section in $G/K$, where 
"$M$ has flat section" means that, for each $x=gK\in M$, 
$\exp^{\perp}(T^{\perp}_{gK}M)$ is a flat totally geodesic submanifold in 
$G/K$.  Denote by $A$ the shape tensor of $M$ and $R$ the curvature tensor of 
$G/K$.  Let $v\in T^{\perp}_xM$ and $X\in T_xM$ ($x=gK$).  
Set $R(v):=R(\cdot,v)v$.  Denote by $\gamma_v$ the geodesic in $G/K$ with 
$\gamma'_v(0)=v$.  The strongly $M$-Jacobi field $Y$ along $\gamma_v$ with 
$Y(0)=X$ (hence $Y'(0)=-A_vX$) is given by 
$$Y(s)=(P_{\gamma_v\vert_{[0,s]}}\circ(\cos(s\sqrt{R(v)})
-\frac{\sin(s\sqrt{R(v)})}{\sqrt{R(v)}}\circ A_v))(X),$$
where $Y'(0)=\widetilde{\nabla}_vY$, $P_{\gamma_v\vert_{[0,s]}}$ is 
the parallel translation along $\gamma_v\vert_{[0,s]}$, 
and 
$\cos(s\sqrt{R(v)})$ and $\frac{\sin(s\sqrt{R(v)})}{s\sqrt{R(v)}}$ are defined by 
$$
\cos(s\sqrt{R(v)}):=\sum_{k=0}^{\infty}
\frac{(-1)^ks^{2k}R(v)^k}{(2k)!}\quad\,\,{\rm and}\quad\,\,
\frac{\sin(s\sqrt{R(v)})}{s\sqrt{R(v)}}
:=\sum_{k=0}^{\infty}\frac{(-1)^ks^{2k}R(v)^k}{(2k+1)!},
$$
respectively.  
Since $M$ has flat section, all focal radii of $M$ along $\gamma_v$ are given 
as zero points of strongly $M$-Jacobi fields along $\gamma_v$.  
Hence all focal radii of $M$ along $\gamma_v$ coincide with the zero points of 
the real-valued function $F_v$ over ${\Bbb R}$ defined by 
$$F_v(s):={\rm det}\left(\cos(s\sqrt{R(v)})
-\frac{\sin(s\sqrt{R(v)})}{\sqrt{R(v)}}\circ A_v\right).$$
So we defined the notion of a complex focal radius of $M$ along $\gamma_v$ as 
the zero points of the complex-valued function $F_v^{\bf c}$ over ${\Bbb C}$ 
defined by 
$$F_v^{\bf c}(z):={\rm det}\left(\cos(z\sqrt{R(v)^{\bf c}})
-\frac{\sin(z\sqrt{R(v)^{\bf c}})}{\sqrt{R(v)^{\bf c}}}\circ A_v^{\bf c}
\right),$$
where $R(v)^{\bf c}$ (resp. $A_v^{\bf c}$) is the complexification of $R(v)$ 
(resp. $A_v$).  
Also, for a complex focal radius $z$ of $M$ along $\gamma_v$, we call 
$$
{\rm dim}_{\bf c}{\rm Ker}\left(
\cos(z\sqrt{R(v)^{\bf c}})
-\frac{\sin(z\sqrt{R(v)^{\bf c}})}{\sqrt{R(v)^{\bf c}}}\circ A_v^{\bf c}
\right)
$$
the {\it multiplicity} of the complex focal radius $z$.  
Here we note that complex focal radii along $\gamma_v$ 
indicate the positions of focal points of the extrinsic 
complexification $M^{\bf c}(\hookrightarrow G^{\bf c}/K^{\bf c})$ of $M$ 
along the complexified geodesic $\gamma_{\iota_{\ast}v}^{\bf c}$, where 
$G^{\bf c}/K^{\bf c}$ is the anti-Kaehlerian symmetric space associated with 
$G/K$ and $\iota$ is the natural embedding of $G/K$ into 
$G^{\bf c}/K^{\bf c}$.  
See [Koi2] about the definitions of $G^{\bf c}/K^{\bf c},\,
M^{\bf c}(\hookrightarrow G^{\bf c}/K^{\bf c})$ and 
$\gamma_{\iota_{\ast}v}^{\bf c}$.  
Furthermore, assume that the normal holonomy group of $M$ is trivial.  
Let $\tilde v$ be a parallel unit normal vector field of $M$.  
Assume that the number (which may be $\infty$) of distinct complex 
focal radii along $\gamma_{\tilde v_x}$ is independent of the choice of 
$x\in M$.  
Let $\{r_{i,x}\,\vert\,i=1,2,\cdots\}$ 
be the set of all complex focal radii along $\gamma_{\tilde v_x}$, where 
$\vert r_{i,x}\vert\,<\,\vert r_{i+1,x}\vert$ or 
"$\vert r_{i,x}\vert=\vert r_{i+1,x}\vert\,\,\&\,\,{\rm Re}\,r_{i,x}
>{\rm Re}\,r_{i+1,x}$" or 
"$\vert r_{i,x}\vert=\vert r_{i+1,x}\vert\,\,\&\,\,
{\rm Re}\,r_{i,x}={\rm Re}\,r_{i+1,x}\,\,\&\,\,
{\rm Im}\,r_{i,x}=-{\rm Im}\,r_{i+1,x}<0$".  
Let $r_i$ ($i=1,2,\cdots$) be complex 
valued functions on $M$ defined by assigning $r_{i,x}$ to each $x\in M$.  
We call these functions $r_i$ ($i=1,2,\cdots$) {\it complex 
focal radius functions for} $\tilde v$.  
If, for each parallel 
unit normal vector field $\tilde v$ of $M$, the set of all complex 
focal radii along $\gamma_{\tilde v_x}$ is independent of the choice of 
$x\in M$, if each complex focal radius function for $\tilde v$ 
is constant on $M$ and if it has constant multiplicity, then 
we call $M$ a {\it complex equifocal submanifold}.  

Next we recall the notion of a proper complex isoparametric submanifold in 
an (infinite dimensional) pseudo-Hilbert space.  
Let $M$ be a pseudo-Riemannian submanifold of finite codimension in 
a pseudo-Hilbert space $(V,\langle\,\,,\,\,\rangle)$.  See [Koi1] about this 
definition.  
We call $M$ a {\it Fredholm pseudo-Riemannian submanifold} 
(or simply {\it Fredholm submanifold}) if 
there exists an orthogonal time-space decomposition 
$V=V_-\oplus V_+$ such that $(V,\langle\,\,,\,\,\rangle_{V_{\pm}})$ is 
a Hilbert space 
and that, for each $v\in T^{\perp}M$, $A_v$ is a compact operator 
with respect to $f^{\ast}\langle\,\,,\,\,\rangle_{V_{\pm}}$, where 
an orthogonal time-space decomposition $V=V_-\oplus V_+$ means that 
$\langle\,\,,\,\,\rangle\vert_{V_-\times V_-}$ is negative definite, 
$\langle\,\,,\,\,\rangle\vert_{V_+\times V_+}$ is positive definite and that 
$\langle\,\,,\,\,\rangle\vert_{V_-\times V_+}=0$, and 
$\langle\,\,,\,\,\rangle_{V_{\pm}}:=-\pi_{V_-}^{\ast}\langle\,\,,\,\,\rangle
+\pi_{V_+}^{\ast}\langle\,\,,\,\,\rangle$ ($\pi_{V_-}$ (resp. $\pi_{V_+}$)
$\,:\,$the orthogonal projection of $V$ onto $V_-$ (resp. $V_+$)).  
Since $A_v$ is a compact operator with respect to 
$f^{\ast}\langle\,\,,\,\,\rangle_{V_{\pm}}$, for each $v\in T^{\perp}M$, 
the operator ${\rm id}-A_v$ is a Fredholm operator with respect to 
$f^{\ast}\langle\,\,,\,\,\rangle_{V_{\pm}}$ and hence the normal exponential 
map $\exp^{\perp}\,:\,T^{\perp}M\to V$ of $M$ is a Fredholm map with respect 
to the metric of $T^{\perp}M$ naturally defined from 
$f^{\ast}\langle\,\,,\,\,\rangle_{V_{\pm}}$ and 
$\langle\,\,,\,\,\rangle_{V_{\pm}}$, 
where ${\rm id}$ is the identity transformation of $TM$.  
In the sequel, set 
$\langle\,\,,\,\,\rangle:=f^{\ast}\langle\,\,,\,\,\rangle$.  
The set of all eigenvalues of the complexification 
$A_v^{\bf c}$ of $A_v$ is described as 
$\{0\}\cup\{\mu_i\,\vert\,i=1,2,\cdots\}$, where "$\vert\mu_i\vert\,>\,
\vert\mu_{i+1}\vert$" or 
"$\vert\mu_i\vert=\vert\mu_{i+1}\vert\,\,\&\,\,
{\rm Re}\,\mu_i>{\rm Re}\,\mu_{i+1}$" or 
"$\vert\mu_i\vert=\vert\mu_{i+1}\vert\,\,\&\,\,
{\rm Re}\,\mu_i={\rm Re}\,\mu_{i+1}\,\,\&\,\,
{\rm Im}\,\mu_i=-{\rm Im}\,\mu_{i+1}>0$".  
We call $\mu_i$ the $i$-th complex principal curvature for $v$.  
Assume that the normal holonomy group $M$ is trivial.  
Fix a parallel normal vector field $\tilde v$ on $M$.  
Assume that the number (which may be $\infty$) 
of distinct complex principal curvatures of $\tilde v_x$ is 
independent of $x\in M$.  Then we define functions $\tilde{\mu}_i$ 
($i=1,2,\cdots$) on $M$ by assigning the $i$-th complex principal curvature 
for $\tilde v_x$ to each $x\in M$.  We call this function 
$\tilde{\mu}_i$ the $i$-th {\it complex principal curvature function for} 
$\tilde v$.  
A Fredholm submanifold $M$ is called a 
{\it complex isoparametric submanifold} if the normal holonomy group of $M$ is 
trivial and if, for each parallel normal vector field $\tilde v$, 
the number of distinct complex principal curvatures of direction 
$\tilde v_x$ is independent of the choice of $x\in M$ and if 
each complex principal curvature function of direction $\tilde v$ 
is constant on $M$.  
Assume that $M$ is a complex isoparametric submanifold.  
If, for each $v\in T^{\perp}M$, the complexified shape operator $A_v^{\bf c}$ 
is diagonalizable with respect to a a pseudo-orthonormal base of 
$(T_xM)^{\bf c}$ ($x\,:\,$ the base point of $v$), that is, 
there exists a pseudo-orthonormal base of $(T_xM)^{\bf c}$ consisting 
of the eigenvectors of $A_v^{\bf c}$, then we call $M$ a 
{\it proper complex isoparametric submanifold}, where a pseudo-orthonormal 
base means a linearly independent system $\{e_i\}_{i=1}^{\infty}$ of 
a pseudo-Hilbert space $(T_xM,\langle\,\,,\,\,\rangle)$ such that, for each 
$i\in {\Bbb N}$, there exists $\hat i\in{\Bbb N}$ 
satisfying $\vert\langle v_i,v_j\rangle\vert=\delta_{\hat ij}$ 
($j\in{\Bbb N}$) ($\delta_{\cdot\cdot}\,:\,$ the Kronecker's symbol) and that 
$\displaystyle{\overline{\mathop{\oplus}_{i=1}^{\infty}{\rm Span}\{v_i\}}
=T_xM}$ ($\overline{\,\,\cdot\,\,}\,:\,$ the closure of $\cdot$ with respect 
to the original topology of $T_xM$).  
Then, for each $x\in M$, $A_v^{\bf c}$'s ($v\in T^{\perp}_xM$) are 
simultaneously diagonalizable with respect to a pseudo-orthonormal base 
of $(T_xM)^{\bf c}$ because $A_v^{\bf c}$'s commute.  
There exists a family $\{E_i\,\vert\,i\in I\}$ 
($I\subset {\Bbb N}$) of parallel subbundles of $(TM)^{\bf c}$ such that, 
for each $x\in M$, 
$\displaystyle{(T_xM)^{\bf c}=\overline{\mathop{\oplus}_{i\in I}(E_i)_x}}$ 
holds and that this decomposition is a common-eigenspace decomposition of 
$A_v^{\bf c}$'s ($v\in T_x^{\perp}M$).  
Also, there exist smooth sections $\lambda_i$ ($i\in I$) of 
$((T^{\perp}M)^{\bf c})^{\ast}$ such that $A_v^{\bf c}
=\lambda_i(v){\rm id}$ on $(E_i)_x$ for each 
$v\in (T^{\perp}M)^{\bf c}$, where 
$x$ is the base point of $v$.  The subbundles $E_i$ ($i\in I$) are called 
{\it complex curvature distributions of} $M$ and $\lambda_i$ ($i\in I$) 
are called {\it complex principal curvatures of} $M$.  
Define a complex normal vector field ${\bf n}_i$ ($i\in I$) by 
$\lambda_i(\cdot)=\langle {\bf n}_i,\cdot\rangle^{\bf c}$, where 
$\langle\,\,,\,\,\rangle^{\bf c}$ is the complexification of 
$\langle\,\,,\,\,\rangle$.  
Note that each ${\bf n}_i$ is parallel with respect to the complexification 
$\nabla^{\perp{\bf c}}$ of $\nabla^{\perp}$.  The normal vector fields 
${\bf n}_i$ ($i\in I$) are called {\it complex curvature normals of} $M$.  

Let $G/K$ be a symmetric space of non-compact type and $\pi:G\to G/K$ be 
the natural projection.  The parallel transport map $\phi$ 
for the semi-simple Lie group $G$ is defined by 
$\phi(u):=g_u(1)$ ($u\in H^0([0,1],\mathfrak g))$, where $g_u$ is the element 
of $H^1([0,1],G)$ with $g_u(0)=e$ ($e,:\,$the identity element of $G$) and 
$g_{u\ast}^{-1}g'_u=u$.  Here we note that $H^0([0,1],\mathfrak g)$ is 
a pseudo-Hilbert space.  
See [Koi1] the detail of the definition of the pseudo-Hilbert space 
$H^0([0,1],\mathfrak g)$ and $\phi$.  
Let $M$ be a complex equifocal submanifold in $G/K$.  
Since $M$ is complex equifocal, $\widetilde M:=(\pi\circ\phi)^{-1}(M)$ is 
complex isoparametric.  
In particular, if $\widetilde M$ is proper complex isoparametric, 
then $M$ is called a {\it proper complex equifocal submanifold}.  
Let $M$ be a proper complex equifocal submanifold in a symmetric space $G/K$ 
of non-compact type.  Denote by $A$ (resp. $\widetilde A$) the shape tensor of 
$M$ (resp. $\widetilde M)$.  
Since $\widetilde M$ is proper complex isoparametric, 
the complexified shape operators of $\widetilde M$ is simultaneously 
diagonalizeble with respect to a pseudo-orthonormal base.  
Hence the complex focal set of $\widetilde M$ at any point 
$u(\in\widetilde M)$ consists of infinitely many complex hyperplanes in the 
complexified normal space $(T^{\perp}_u\widetilde M)^{\bf c}$ and the group 
generated by the complex reflections of order two with respect to the complex 
hyperplanes is discrete.  
From this fact, it follows that, for the complex 
focal set of the proper complex equifocal submanifold $M$, 
the following fact holds:

\vspace{0.3truecm}

$(\ast)\,\,$
{\sl The complex focal set of $M$ at any point $x(\in M)$ consists of 
infinitely many}

{\sl complex hyperplanes in the complexified normal space 
$(T^{\perp}_x M)^{\bf c}$ and the group}

{\sl generated by the complex reflections of order two with respect to 
the complex}

{\sl hyperplanes is discrete.}

\vspace{0.3truecm}


Let $H$ be a symmetric subgroup of $G$ 
(i.e., there exists an involution of $G$ with 
$({\rm Fix}\,\tau)_0\subset H\subset{\rm Fix}\,\tau$), where ${\rm Fix}\,\tau$ 
is the fixed point group of $\tau$ and $({\rm Fix}\,\tau)_0$ is the identity 
component of ${\rm Fix}\,\tau$.  
The natural action $H$ on $G/K$ is called a {\it Hermann type action}.  
It is shown that a principal orbit of a Hermann type ation is a proper 
complex equifocal and curvature-adapted ([Koi3]), where the 
curvature-adaptedness means that, for each normal vector $v$ of $M$, 
$R(\cdot,v)v$ preserves $T_xM$ ($x\,:\,$the base point of $v$) invariantly 
and that $[R(\cdot,v)v,\,A_v]=0$ ($R\,:\,$the curvature tensor of $G/K$).  
Let $P(G,H\times K):=\{g\in H^1([0,1],G)\,\vert\,(g(0),g(1))\in H\times K\}$, 
where $H^1([0,1],G)$ is a pseudo-Hilbert Lie group of all $H^1$-paths in 
$G$ having $[0,1]$ as the domain.  
See [Koi1] about the detail of the definition of $H^1([0,1],G)$.  
This group $P(G,H\times K)$ acts on 
$H^0([0,1],\mathfrak g)$ as gauge action.  It is shown that orbits of the 
$P(G,H\times K)$-action are the inverse images of the $H$-orbits by 
$\pi\circ\phi$ (see [Koi2]).  

\section{The regularized mean curvature vector of Fredholm submanifold with 
proper shape operators} 
In this section, we shall define the regularized mean curvature vector of 
a certain kind of Fredholm submanifold in a pseudo-Hlbert space.  Let $M$ be 
a Fredholm submanifold in a pseudo-Hilbert space 
$(V,\langle\,\,,\,\,\rangle)$.  
Denote by $A$ the shape tensor of $M$.  Fix $v\in T^{\perp}M$.  
If the complexified shape operator $A_v^{\bf c}$ is diagoalizable with respect 
to a pseudo-orthonormal base, then $A_v^{\bf c}$ is said to be {\it proper}.  
If $A_v^{\bf c}$ is proper for any $v\in T^{\perp}M$, then we say that $M$ 
{\it has proper shape operators}.  Assume that $M$ has proper shape 
operators.  Fix $v\in T_x^{\perp}M$.  
Let $\{\mu_i\,\vert\,i=1,2,\cdots\}$ 
("$\vert\mu_i\vert\,>\,\vert\mu_{i+1}\vert$" or 
"$\vert\mu_i\vert=\vert\mu_{i+1}\vert\,\,\&\,\,
{\rm Re}\,\mu_i>{\rm Re}\,\mu_{i+1}$" or 
"$\vert\mu_i\vert=\vert\mu_{i+1}\vert\,\,\&\,\,
{\rm Re}\,\mu_i={\rm Re}\,\mu_{i+1}\,\,\&\,\,
{\rm Im}\,\mu_i=-{\rm Im}\,\mu_{i+1}>0$") be the set of all eigenvalues of 
$A_v^{\bf c}$ other than zero and $m_i$ the multiplicity of $\mu_i$.  
Then we define the regularized trace 
${\rm Tr}_rA_v^{\bf c}$ of $A_v^{\bf c}$ by 
${\rm Tr}_rA_v^{\bf c}:=\sum\limits_{i}m_i\mu_i$.  
Also, we define the trace ${\rm Tr}_{\rm abs}(A_v^{\bf c})^2$ by 
${\rm Tr}_{\rm abs}(A_v^{\bf c})^2:=\sum\limits_{i}m_i\vert\mu_i\vert^2$.  
If there exist ${\rm Tr}_rA_v^{\bf c}$ and ${\rm Tr}_{\rm abs}(A_v^{\bf c})^2$ 
for each $v\in T^{\perp}M$, then we say that $M$ is {\it regularizable}.  
It is shown that, if $\mu$ is an eigenvalue of $A^{\bf c}_v$ with multiplicity 
$m$, then so is also the conjugate $\bar{\mu}$ of $\mu$.  Hence we have 
${\rm Tr}_rA^{\bf c}_v\in{\Bbb R}$.  
Define $H_x\in T^{\perp}_xM$ by 
$\langle H_x,v\rangle={\rm Tr}_rA_v^{\bf c}$ ($\forall\,v\in T^{\perp}_xM$).  
We call the normal vector field $H\,(:x\mapsto H_x)$ of $M$ the 
{\it regularized mean curvature vector} of $M$.  
Let $f_t:M\hookrightarrow V$ ($0\leq t<T$) be a $C^{\infty}$-family of 
regularizable Fredholm submanifolds with proper shape operators and 
$H_t$ be the regularized mean curvature vector of $f_t$.  Define by 
$\widetilde f:M\times[0,T)\to V$ by $\widetilde f(x,t):=f_t(x)$ ($(x,t)\in 
M\times[0,T)$).  If $\widetilde f_{\ast}(\frac{\partial}{\partial t})=H_t$, 
then we call $f_t$ ($0\leq t<T$) the {\it regularized mean curvature flow}.  
In the sequel, we call this flow the mean curvature flow for simplicity.  

Let $G/K$ be a symmetric space of non-compact type, 
$\pi:G\to G/K$ the natural projection and 
$\phi:H^0([0,1],\mathfrak g)\to G$ the parallel tansport map for $G$.  Let 
$M$ be a curvature-adapted isoprametric submanifold in $G/K$ having no focal point of non-Euclidean type 
on the ideal boundary of $G/K$ and set $\widetilde M:=(\pi\circ\phi)^{-1}(M)$, 
which is proper complex isoparametric (hence has proper shape operators).  
Denote by $H$ the mean curvature vector of $M$.  Then we have the following 
fact.  

\vspace{0.5truecm}

\noindent
{\bf Lemma 3.1.} {\sl The submanifold $\widetilde M$ is regularizable and 
the regularized mean curvature vector $\widetilde H$ is equal to 
the horizontal lift $H^L$ of $H$.}

\vspace{0.5truecm}

\noindent
{\it Proof.} Without loss of generality, we may assume that $eK\in M$.  
For simplicity, set $\mathfrak m:=T_{eK}M$ and 
$\mathfrak b:=T_{eK}^{\perp}M$.  
Since $M$ is flat section (hence $\mathfrak b$ is abelian), the normal 
connection of $M$ is flat and since $M$ is curvature-adapted, the operators 
$R(\cdot,v)v$'s ($v\in\mathfrak b$) and $A_v$'s ($v\in\mathfrak b$) commute 
to one another.  Also they are diagonalizable with respect to an orthonormal 
base, respectively.  Therefore they are simultaneously diagonalizable with 
respect to an orthonormal base.  Let $\mathfrak m=\mathfrak m^R_0
+\sum\limits_{i\in I^R}\mathfrak m^R_i$ be the common eigenspace decomposition 
of $R(\cdot,v)v$'s ($v\in\mathfrak b$) and $\mathfrak m=\mathfrak m^A_0
+\sum\limits_{i\in I^A}\mathfrak m^A_i$ be the common eigenspace 
decomposition of $A_v$'s ($v\in\mathfrak b$), where 
$\displaystyle{\mathfrak m^R_0:=\mathop{\cap}_{v\in\mathfrak b}{\rm Ker}
R(\cdot,v)v}$ and 
$\displaystyle{\mathfrak m^A_0
:=\mathop{\cap}_{v\in\mathfrak b}{\rm Ker}\,A_v}$.  
Set $m_i^R:={\rm dim}\,\mathfrak m^R_i$ and 
$m_i^A:={\rm dim}\,\mathfrak m^A_i$.  
Also, set $I^A_i:=\{j\in I^A\cup\{0\}\,\vert\,\mathfrak m^A_j\cap
\mathfrak m^R_i\not=\{0\}\}$ ($i\in I^R\cup\{0\}$).  Since $R(\cdot,v)v$ 
($v\in\mathfrak b$) and $A_v$'s ($v\in\mathfrak b$) are simultaneously 
diagonalizable, we have 
$\mathfrak m=\sum\limits_{i\in I^R\cup\{0\}}\sum\limits_{j\in I^A_i}
(\mathfrak m^A_j\cap\mathfrak m^R_i)$.  
Let $\beta_i(\geq0)$ ($i\in I^R$) and $\lambda_i$ ($i\in I^A$) be 
linear functions over $\mathfrak b$ defined by 
$R(\cdot,v)v\vert_{\mathfrak m^R_i}=-\beta_i(v)^2{\rm id}$ ($v\in\mathfrak b$) 
and $A_v\vert_{\mathfrak m^A_i}=\lambda_i(v){\rm id}$ ($v\in\mathfrak b$).  
Denote by $\mathfrak b_r$ the set of all $v\in\mathfrak b$ such that 
$\beta_i(v)\not=0,\,\,\lambda_i(v)\not=0,\,\,$ $\beta_i(v)$'s ($i\in I^R$) 
are mutually distinct and that so are also $\lambda_i(v)$'s ($i\in I^A$).  
Note that $\mathfrak b_r$ is open and dense in $\mathfrak b$.  
Fix $v\in\mathfrak b_r$.  Denote by $\widetilde A$ the shape tensor of 
$\widetilde M$ and ${\rm Spec}\,\widetilde A^{\bf c}_{v^L}$ the spectrum of 
$\widetilde A^{\bf c}_{v^L}$, where $v^L$ is the horizontal lift of $v$ to 
the constant path $\hat 0$ at the zero vector $0$ of $\mathfrak g$.  
Set $I^A_{i,v,+}:=\{j\in I^A_i\,\vert\,\vert\lambda_j(v)\vert>\vert\beta_i(v)
\vert\},\,\,$ $I^A_{i,v,-}:=\{j\in I^A_i\,\vert\,\vert\lambda_j(v)\vert<
\vert\beta_i(v)\vert\}$ and 
$I^A_{i,v,0}:=\{j\in I^A_i\,\vert\,\vert\lambda_j(v)\vert
=\vert\beta_i(v)\vert\}$.  Since $M$ is a curvature-adapted isoparametric 
submanifold satisfying the condition $(\ast_{\bf c})$, we can show that 
$I^A_{i,v,0}=\emptyset$ (see Theorem A of [Koi1]) and that 
$I^A_{i,v,+}$ and $I^A_{i,v,-}$ are at most one point sets, respectively 
(see the proof of Theorems B and C of [Koi6]).  
When $I^A_{i,v,+}\not=\emptyset$ (resp. $I^A_{i,v,-}\not=\emptyset$), 
denote by $j^+_{i,v}$ (resp. $j^-_{i,v}$) the only element.  
Set $I^R_{v,+}:=\{i\in I^R\,\vert\,I^A_{i,v,+}\not=\emptyset\}$ and 
$I^R_{v,-}:=\{i\in I^R\,\vert\,I^A_{i,v,-}\not=\emptyset\}$.  
For simplicity, set $\mathfrak m_{i,v}^{\pm}
:=\mathfrak m^A_{j^{\pm}_{i,v}}\cap\mathfrak m^R_i$ and 
$m_{i,v}^{\pm}:={\rm dim}\,\mathfrak m_{i,v}^{\pm}$.  
According to Theorem 5.9 in [Koi1], 
${\rm Spec}\widetilde A^{\bf c}_{v^L}\setminus
\{0\}$ is given by 
$$\begin{array}{l}
\displaystyle{{\rm Spec}\,\widetilde A^{\bf c}_{v^L}\setminus\{0\}=
\left\{\left.\frac{\beta_i(v)}{{\rm arctanh}(\beta_i(v)/\lambda_{j^+_{i,v}}(v))
+k\pi\sqrt{-1}}\,\right\vert\,i\in I^R_{v,+},\,k\in{\Bbb Z}\right\}}\\
\hspace{1.5truecm}\displaystyle{\cup
\left\{\left.\frac{\beta_i(v)}{{\rm arctanh}(\lambda_{j^-_{i,v}}(v)/\beta_i(v))
+(k+\frac12)\pi\sqrt{-1}}\,\right\vert\,i\in I^R_{v,-},\,k\in{\Bbb Z}\right\}}
\end{array}
\leqno{(3.1)}$$
Hence we have 
$$\begin{array}{l}
\displaystyle{{\rm Tr}_r\widetilde A^{\bf c}_{v^L}=\sum_{i\in I^R_{v,+}}
\sum_{k\in{\Bbb Z}}\frac{\beta_i(v)}{{\rm arctanh}(\beta_i(v)/
\lambda_{j^+_{i,v}}(v))+k\pi\sqrt{-1}}\times m_{i,v}^+}\\
\hspace{2truecm}\displaystyle{
+\sum_{i\in I^R_{v,-}}
\sum_{k\in{\Bbb Z}}\frac{\beta_i(v)}{{\rm arctanh}
(\lambda_{j^-_{i,v}}(v)/\beta_i(v))+(k+\frac12)\pi\sqrt{-1}}
\times m_{i,v}^-}\\
\hspace{1.5truecm}\displaystyle{
=\sum_{i\in I^R_{v,+}}m_{i,v}^+\lambda_{j^+_{i,v}}(v)
+\sum_{i\in I^R_{v,-}}m_{i,v}^-\lambda_{j^-_{i,v}}(v)}\\
\hspace{1.5truecm}\displaystyle{
=\sum_{j\in I^A}m_j\lambda_j(v)={\rm Tr}\,A_v\,\,\,\,(\in{\Bbb R})}
\end{array}$$
in terms of ${\rm coth}\theta=\sum_{j\in{\Bbb Z}}
\frac{1}{\theta+j\pi\sqrt{-1}}$ 
and ${\rm coth}(\theta+\frac{\pi\sqrt{-1}}{2})=\tanh\theta$.  
Also, we have 
$${\rm Tr}_{\rm abs}(\widetilde A^{\bf c}_{v^L})^2
\leq\sum_{i\in I^R_{v,+}}\sum_{k\in{\Bbb Z}}
\frac{\vert\beta_i(v)\vert^2}{k^2}
+\sum_{i\in I^R_{v,-}}\sum_{k\in{\Bbb Z}}\frac{\vert\beta_i(v)\vert^2}{k^2}
<\infty.$$
Hence $\widetilde M$ is regularizable and 
${\rm Tr}_r\widetilde A^{\bf c}_{v^L}={\rm Tr}\,A_v$ holds.  
This implies $\langle\widetilde H_{\hat 0},v^L\rangle=\langle H_{eK},v\rangle
(=\langle(H^L)_{\hat 0},v^L\rangle)$.  Since this relation holds for any 
$v\in\mathfrak b_r$ and $\mathfrak b_r$ is dense in $\mathfrak b$, we obtain 
$\widetilde H_{\hat0}=(H^L)_{\hat 0}$.  
Similarly we can show 
$\widetilde H_u=(H^L)_u$ for any $u\in\widetilde M$.  
Thus we obtain $\widetilde H=H^L$.
\hspace{9.35truecm}q.e.d.

\vspace{0.5truecm}

By using Lemma 3.1 and imitating the proof of Lemma 3.1 of [Koi10], we 
can show the following fact.  

\vspace{0.5truecm}

\noindent
{\bf Lemma 3.2.} {\sl The mean curvature flow $\widetilde M_t$ 
(resp. $M_t$) for $\widetilde M$ (resp. $M$) exists in short time and 
$\widetilde M_t=(\pi\circ\phi)^{-1}(M_t)$ holds.}

\section{Proofs of Theorems A and B} 
Let $M$ be a curvature-adapted isoparametric submanifold 
in a symmetric space $G/K$ of non-compact type as in Theorem A.  
Without loss of generality, we may assume that $eK\in M$.  
We use the notations in the previous section.  
Denote by $\mathfrak F$ the focal set of $M$ at $eK$.  
Since $\pi\circ\phi$ is a pseudo-Riemannian submersion, the focal set 
$\widetilde{\mathfrak F}$ of $\widetilde M$ at $\hat 0$ is equal to 
$\{v^L_{\hat 0}\,\vert\,\exp^{\perp}(v)\in\mathfrak F\}$, where $\exp^{\perp}$ 
is the normal exponential map of $M$ and $v_{\hat 0}^L$ is the horizontal llft 
of $v$ to $\hat 0$.  Here we regard the normal space 
$T^{\perp}_{\hat 0}\widetilde M$ of $\widetilde M$ at $\hat 0$ as 
a subspace of $H^0([0,1],\mathfrak g)$.  
In the sequel, we identify $v^L_{\hat 0}$ with $v$ 
through $(\pi\circ\phi)_{\ast\hat 0}$.  The focal set 
$\widetilde{\mathfrak F}$ is equal to 
$\{v\,\vert\,{\rm Ker}(\widetilde A_v-{\rm id})\not=\{0\}\}\,
(\subset\mathfrak b)$.  
The complex focal structure $\widetilde{\mathfrak F}^{\bf c}$ of 
$\widetilde M$ at $\hat 0$ is defined by 
$\widetilde{\mathfrak F}^{\bf c}:=
\{v\,\vert\,{\rm Ker}(\widetilde A_v^{\bf c}-{\rm id})\not=\{0\}\}\,
(\subset\mathfrak b^{\bf c})$.  According to the proof of Theorems B and C 
of [Koi6], by using $(3.1)$ and discussing delicately, 
we can show that $\beta_i(v)/\lambda_{j_{i,v}^+}(v)$ and 
$\lambda_{j_{i,v}^-}(v)/\beta_i(v)$ are independent of the choice of $v$ 
(in the sequel, we denote these constants by $c_i^+$ and $c_i^-$, 
respectively), $I^R_{v,+}$ and $I^R_{v,-}$ are independent of the choice of 
$v$ (in the sequel, we denote these sets by $I^R_+$ and $I^R_-$) and that 
$\widetilde{\mathfrak F}^{\bf c}$ is described as follows:
$$\begin{array}{l}
\displaystyle{\widetilde{\mathfrak F}^{\bf c}
=\left(\mathop{\cup}_{i\in I^R_+}\mathop{\cup}_{j\in{\Bbb Z}}
\beta_i^{{\bf c}-1}({\rm arctanh}\,c_i^++j\pi\sqrt{-1})\right)}\\
\hspace{1truecm}\displaystyle{\cup\left(\mathop{\cup}_{i\in I^R_-}
\mathop{\cup}_{j\in{\Bbb Z}}
\beta_i^{{\bf c}-1}({\rm arctanh}\,c_i^-+(j+\frac12)\pi\sqrt{-1})\right).}
\end{array}
\leqno{(4.1)}$$
Since $\lambda_{j_{i,v}^+}=\frac{1}{c_i^+}\beta_i$ and 
$\lambda_{j_{i,v}^-}=c_i^-\beta_i$, they are independent of the choice of 
$v\in\mathfrak b$.  Hence we denote $\lambda_{j_{i,v}^{\pm}}$ by 
$\lambda_i^{\pm}$.  
Therefore, $\widetilde{\mathfrak F}$ is given by 
$$\widetilde{\mathfrak F}
=\mathop{\cup}_{i\in I^R_+}\beta_i^{-1}({\rm arctanh}\,c_i^+).\leqno{(4.2)}$$
Denote by $\Lambda$ the set of all complex principal curvatures of 
$\widetilde M$.  According to $(3.1)$, $\Lambda$ is given by 
$$\begin{array}{l}
\displaystyle{\Lambda=\left\{
\left.\frac{\widetilde{\beta_i^{\bf c}}}{{\rm arctanh}\,c_i^++j\pi\sqrt{-1}}
\,\right\vert\,i\in I^R_+,\,\,j\in{\Bbb Z}\right\}}\\
\displaystyle{\cup\left\{
\left.\frac{\widetilde{\beta_i^{\bf c}}}
{{\rm arctanh}\,c_i^-+(j+\frac12)\pi\sqrt{-1}}\,\right\vert\,i\in I^R_-,\,\,
j\in{\Bbb Z}\right\},}
\end{array}$$
where $\widetilde{\beta_i^{\bf c}}$ is the parallel section of 
$((T^{\perp}\widetilde M)^{\bf c})^{\ast}$ with 
$(\widetilde{\beta_i^{\bf c}})_{\hat 0}=\beta_i^{\bf c}$.  
From the assumption, $M$ admits a focal submanifold.  
Hence we have $I^R_+\not=\emptyset$ and 
$\displaystyle{\mathop{\cap}_{i\in I^R_+}\beta_i^{-1}
({\rm arctanh}\,c_i^+)\not=\emptyset}$ (see the proof of Theorems B and C 
in [Koi6]).  Fix $\displaystyle{O\in\mathop{\cap}_{i\in I^R_+}\beta_i^{-1}
({\rm arctanh}\,c_i^+)}$.  
Set 
$$\begin{array}{l}
\displaystyle{\triangle':=\{\beta_i\,\vert\,i\in I^R\}\cup
\{-\beta_i\,\vert\,i\in I^R\},}\\
\displaystyle{{\triangle'}^V:=\{\beta_i\,\vert\,i\in I^R_+\}\cup
\{-\beta_i\,\vert\,i\in I^R_+\},}\\
\displaystyle{{\triangle'}^H:=\{\beta_i\,\vert\,i\in I^R_-\}\cup
\{-\beta_i\,\vert\,i\in I^R_-\}.}
\end{array}$$
Let $\mathfrak a$ be a maximal abelian subspace of $\mathfrak p$ containing 
$\mathfrak b(=T^{\perp}_{eK}M)$ and $\triangle$ be the root system of $G/K$ 
with respect to $\mathfrak a$.  
Then we have $\triangle'=\{\alpha\vert_{\mathfrak b}\,\vert\,
\alpha\in\triangle\,\,{\rm s.t.}\,\,\alpha\vert_{\mathfrak b}\not=0\}$.  
Let $F_{\it l}$ be the focal submanifold of $M$ through 
$x_0:=\exp^{\perp}(O)$, 
which is one of the lowest dimensional focal submanifolds of $M$.  
For simplicity, we set 
$$\begin{array}{ll}
\displaystyle{\widetilde{\lambda}_i^+:=\frac{\widetilde{\beta_i^{\bf c}}}
{{\rm arctanh}\,c_i^+}} & \displaystyle{(i\in I^R_+)}\\
\displaystyle{\widetilde{\lambda}_i^-:=\frac{\widetilde{\beta_i^{\bf c}}}
{{\rm arctanh}\,c_i^-+\frac12\pi\sqrt{-1}}} & \displaystyle{(i\in I^R_-)}
\end{array}$$
and 
$$\begin{array}{ll}
\displaystyle{b_i^+:=\frac{\pi}{{\rm arctanh}\,c_i^+}} & 
\displaystyle{(i\in I^R_+)}\\
\displaystyle{b_i^-:=\frac{\pi}{{\rm arctanh}\,c_i^-+\frac12\pi\sqrt{-1}}} & 
\displaystyle{(i\in I^R_-).}
\end{array}$$
Then we have 
$$\begin{array}{l}
\displaystyle{\Lambda=\left\{\left.\frac{\widetilde{\lambda}_i^+}
{1+b_i^+j\sqrt{-1}}\,\right\vert\,i\in I^R_+,\,\,j\in{\Bbb Z}\right\}}\\
\hspace{1truecm}\displaystyle{
\cup\left\{\left.\frac{\widetilde{\lambda}_i^-}{1+b_i^-j\sqrt{-1}}\,
\right\vert\,i\in I^R_-,\,\,j\in{\Bbb Z}\right\}.}
\end{array}$$
For simplicity, we set $\widetilde{\lambda}_{ij}^{\pm}:=
\frac{\widetilde{\lambda}^{\pm}_i}{1+b_i^{\pm}j\sqrt{-1}}$ 
($i\in I^R_{\pm},\,j\in{\Bbb Z}$).  Take $v\in\mathfrak b_r\setminus
\widetilde{\mathfrak F}$, where we note that 
$\displaystyle{\widetilde{\mathfrak F}=
\mathop{\cup}_{i\in I^R_+}(\widetilde{\lambda}_i^+)_{\hat0}^{-1}(1)}$.  
We have ${\rm dim}\,{\rm Ker}\,({\widetilde A}_v^{\bf c}
-(\widetilde{\lambda}_{ij}^{\pm})_{\hat 0}(v){\rm id})=m_i^{\pm}$ 
($i\in I^R_{\pm},\,j\in{\Bbb Z}$).  
Set $E_{ij}^{\pm}:={\rm Ker}
(\widetilde A^{\bf c}_v-(\widetilde{\lambda}_{ij}^{\pm})_{\hat 0}(v){\rm id})$ 
($i\in I^R_{\pm},\,j\in{\Bbb Z}$), which are independent of the choice of 
$v\in\mathfrak b_r\setminus\widetilde{\mathfrak F}$.  
Take another $w\in\mathfrak b_r\setminus\widetilde{\mathfrak F}$.  
Let $\widetilde w$ be the parallel 
normal vector field of $\widetilde M$ with $\widetilde w_{\hat0}=w$.  
Denote by $\eta_{\widetilde w}$ the end-point map for $\widetilde w$ and 
$\widetilde M_w:=\eta_{\widetilde w}(\widetilde M)$, which is a parallel 
submanifold of $\widetilde M$.  
We have 
$$\begin{array}{l}
\displaystyle{T_w\widetilde M_w=\eta_{\widetilde w\ast}
(T_{\hat0}\widetilde M)}\\
\hspace{1.3truecm}\displaystyle{
=\left(\mathop{\oplus}_{i\in I^R_+}\mathop{\oplus}_{j\in{\Bbb Z}}
\eta_{\widetilde w\ast}(E_{ij}^+)\right)\oplus
\left(\mathop{\oplus}_{i\in I^R_-}\mathop{\oplus}_{j\in{\Bbb Z}}
\eta_{\widetilde w\ast}(E_{ij}^-)\right).}
\end{array}$$
Denote by $\widetilde A^w$ the shape tensor 
of $\widetilde M_w$.  We have 
$$(\widetilde A^w)_v^{\bf c}\vert_{\eta_{\widetilde w\ast}(E^{\pm}_{ij})}
=\frac{(\widetilde{\lambda}_{ij}^{\pm})_{\hat 0}(v)}
{1-(\widetilde{\lambda}_{ij}^{\pm})_{\hat 0}(w)}{\rm id}
=\frac{(\widetilde{\lambda}_i^{\pm})_{\hat 0}(v)}
{(1+b^{\pm}_ij\sqrt{-1})-(\widetilde{\lambda}_i^{\pm})_{\hat 0}(w)}{\rm id}.$$
Hence the set $\Lambda^w$ of all complex principal curvatures of 
$\widetilde M_w$ is given by 
$$\begin{array}{l}
\displaystyle{\Lambda^w=\left\{\left.\frac{\widetilde{\lambda}_{ij}^+}
{1-(\widetilde{\lambda}_{ij}^+)_{\hat0}(w)}\,\right\vert\,i\in I^R_+,\,
j\in{\Bbb Z}\right\}}\\
\hspace{1truecm}\displaystyle{\cup
\left\{\left.\frac{\widetilde{\lambda}_{ij}^-}
{1-(\widetilde{\lambda}_{ij}^-)_{\hat0}(w)}\,\right\vert\,i\in I^R_-,\,
j\in{\Bbb Z}\right\}}\\
\displaystyle{=\left\{\left.\frac{\widetilde{\lambda}_i^+}{1+b_i^+j\sqrt{-1}
-(\widetilde{\lambda}_i^+)_{\hat 0}(w)}\,\right\vert\,i\in I^R_+,\,\,
j\in{\Bbb Z}\right\}}\\
\hspace{1truecm}\displaystyle{\cup
\left\{\left.\frac{\widetilde{\lambda}_i^-}{1+b_i^-j\sqrt{-1}
-(\widetilde{\lambda}_i^-)_{\hat 0}(w)}\,\right\vert\,i\in I^R_-,\,\,
j\in{\Bbb Z}\right\}}\\
\displaystyle{=\left\{\left.\frac{\widetilde{\beta_i^{\bf c}}}
{{\rm arctanh}\,c_i^++j\pi\sqrt{-1}-\beta_i(w)}\,\right\vert\,i\in I^R_+,\,\,
j\in{\Bbb Z}\right\}}\\
\hspace{1truecm}\displaystyle{\cup
\left\{\left.\frac{\widetilde{\beta_i^{\bf c}}}
{{\rm arctanh}\,c_i^-+(j+\frac12)\pi\sqrt{-1}-\beta_i(w)}\,\right\vert\,i\in 
I^R_-,\,\,j\in{\Bbb Z}\right\}.}
\end{array}$$
Hence we have 
$$\begin{array}{l}
\displaystyle{{\rm Tr}(\widetilde A^w_v)^{\bf c}
=\sum_{i\in I^R_+}\sum_{j\in{\Bbb Z}}\frac{\beta_i(v)}{{\rm arctanh}\,c_i^+
+j\pi\sqrt{-1}-\beta_i(w)}\times m_i^+}\\
\hspace{2truecm}\displaystyle{
+\sum_{i\in I^R_-}\sum_{j\in{\Bbb Z}}\frac{\beta_i(v)}
{{\rm arctanh}\,c_i^-+(j+\frac12)\pi\sqrt{-1}-\beta_i(w)}\times m_i^-}\\
\hspace{1.5truecm}\displaystyle{
=\sum_{i\in I^R_+}m_i^+\coth({\rm arctanh}\,c_i^+-\beta_i(w))\beta_i(v)}\\
\hspace{2truecm}\displaystyle{
+\sum_{i\in I^R_-}m_i^-\tanh({\rm arctanh}\,c_i^+-\beta_i(w))\beta_i(v)
\,\,\,\,(v\in\mathfrak b),}
\end{array}\leqno{(4.3)}$$
where we use $\sum\limits_{j\in{\Bbb Z}}\frac{1}{\theta+j\pi\sqrt{-1}}
=\coth\theta$ and $\coth(\theta+\frac{\pi\sqrt{-1}}{2})=\tanh\theta$.  
Hence we have 
$$\begin{array}{l}
\displaystyle{\langle(\widetilde H^w)_w,v\rangle
=\langle\sum_{i\in I^R_+}m_i^+\coth({\rm arctanh}\,c_i^+-\beta_i(w))
\beta_i^{\sharp},\,v\rangle}\\
\hspace{2.8truecm}\displaystyle{
+\langle\sum_{i\in I^R_-}m_i^-\tanh({\rm arctanh}\,c_i^--\beta_i(w))
\beta_i^{\sharp},\,v\rangle,}
\end{array}$$
where $\beta_i^{\sharp}$ is defined by $\beta_i(\cdot)
=\langle\beta_i^{\sharp},\cdot\rangle$ ($i\in I^R$).  
Since this relation holds for any $v\in\mathfrak b_r\setminus
\widetilde{\mathfrak F}$, we have 
$$\begin{array}{l}
\displaystyle{(\widetilde H^w)_w=
\sum_{i\in I^R_+}m_i^+\coth({\rm arctanh}\,c_i^+-\beta_i(w))
\beta_i^{\sharp}}\\
\hspace{1.8truecm}\displaystyle{
+\sum_{i\in I^R_-}m_i^-\tanh({\rm arctanh}\,c_i^--\beta_i(w))
\beta_i^{\sharp}.}
\end{array}\leqno{(4.4)}$$
Set 
$$\begin{array}{l}
\displaystyle{\widetilde C:=\{w\in\mathfrak b\,\vert\,
(\widetilde{\lambda}_i^+)_{\hat 0}(w)<1\,\,(i\in I^R_+)\}}\\
\hspace{0.6truecm}\displaystyle{=\{w\in\mathfrak b\,\vert\,
\beta_i(w)<{\rm arctanh}\,c_i^+\,\,(i\in I^R_+)\},}
\end{array}$$
which is a fundamental domain of the real Coxeter group associated with 
$\widetilde M$.  Each parallel submanifold of $M$ passes through the only 
one point of $\exp^{\perp}(\widetilde C)$ and each focal submanifold of $M$ 
passes through the only one point of $\exp^{\perp}(\partial\widetilde C)$.  
Define a vector field $X$ on $\widetilde C$ by 
$X_w:=(\widetilde H^w)_w$ ($w\in\widetilde C$).  
Let $\{\psi_t\}$ be the local one-parameter transformation group of $X$.  
Now we shall prove the statements (i) and (iii) of Theorem A.  

\vspace{0.3truecm}

\noindent
{\it Proof of {\rm(i)} and {\rm(iii)} of Theorem A.}
First we shall show the statement (i).  
Denote by $\widetilde{\sigma}_i$ ($i\in I^R_+$) the maximal dimensional 
simplex of $\partial\widetilde C$ contained in 
$\beta_i^{-1}({\rm arctanh}\,c_i^+)$.  
Fix $i_0\in I^R_+$.  Take $w_0\in\widetilde{\sigma}_{i_0}$ and 
and $w'_0\in\widetilde C$ near $w_0$ such that $w_0-w_0'$ is 
normal to $\widetilde{\sigma}_{i_0}$.  Set 
$w_0^{\varepsilon}:=\varepsilon w'_0+(1-\varepsilon)w_0$ for 
$\varepsilon\in(0,1)$.  
Then we have $\lim\limits_{\varepsilon\to+0}\beta_{i_0}
(w^{\varepsilon}_0)={\rm arctanh}\,c_{i_0}^+$ and 
$\displaystyle{\mathop{\sup}_{0<\varepsilon<1}
\beta_{i_0}(w^{\varepsilon}_0)<{\rm arctanh}\,c_i^+}$ for each 
$i\in I^R_+\setminus\{i_0\}$.  
Hence we have 
$$
\lim_{\varepsilon\to+0}\coth({\rm arctanh}\,c_{i_0}^+-\beta_{i_0}
(w^{\varepsilon}_0))=\infty$$
and
$$\mathop{\sup}_{0<\varepsilon<1}\coth
({\rm arctanh}\,c_i^+-\beta_i(w^{\varepsilon}_0))\,<\,\infty\,\,
(i\in I^R_+\setminus\{i_0\}).$$
Therefore, we have $\lim\limits_{\varepsilon\to+0}
\frac{X_{w_0^{\varepsilon}}}{\vert\vert X_{w_0^{\varepsilon}}\vert\vert}$ is 
the outward unit normal vector of $\widetilde{\sigma}_{i_0}$.  Also we have 
$\lim\limits_{\varepsilon\to +0}\vert\vert X_{w_0^{\varepsilon}}\vert\vert
=\infty$.  From these facts, $X$ is as in the first figure of Figure 1 on 
a sufficiently small collar neighborhood of $\widetilde{\sigma}_{i_0}$.  
Define a function $\rho$ over $\widetilde C$ by 
$$\begin{array}{l}
\displaystyle{\rho(w):=-\sum_{i\in I^R_+}m_i^+\log\sinh
({\rm arctanh}\,c_i^+-\beta_i(w))}\\
\hspace{1.2truecm}\displaystyle{-\sum_{i\in I^R_-}m_i^-\log\cosh
({\rm arctanh}\,c_i^--\beta_i(w))\qquad(w\in\widetilde C).}
\end{array}$$

\vspace{0.2truecm}

\centerline{
\unitlength 0.1in
\begin{picture}( 52.3400, 24.3000)(  9.7600,-35.7000)
%
\special{pn 8}%
\special{pa 4220 2720}%
\special{pa 6210 2720}%
\special{fp}%
%
\special{pn 8}%
\special{pa 4220 2720}%
\special{pa 5220 1140}%
\special{fp}%
%
\special{pn 8}%
\special{pa 4666 2700}%
\special{pa 4666 2850}%
\special{fp}%
\special{sh 1}%
\special{pa 4666 2850}%
\special{pa 4686 2784}%
\special{pa 4666 2798}%
\special{pa 4646 2784}%
\special{pa 4666 2850}%
\special{fp}%
%
\special{pn 8}%
\special{pa 5116 2700}%
\special{pa 5116 2850}%
\special{fp}%
\special{sh 1}%
\special{pa 5116 2850}%
\special{pa 5136 2784}%
\special{pa 5116 2798}%
\special{pa 5096 2784}%
\special{pa 5116 2850}%
\special{fp}%
%
\special{pn 8}%
\special{pa 5906 2700}%
\special{pa 5906 2850}%
\special{fp}%
\special{sh 1}%
\special{pa 5906 2850}%
\special{pa 5926 2784}%
\special{pa 5906 2798}%
\special{pa 5886 2784}%
\special{pa 5906 2850}%
\special{fp}%
%
\special{pn 8}%
\special{pa 4440 2410}%
\special{pa 4316 2326}%
\special{fp}%
\special{sh 1}%
\special{pa 4316 2326}%
\special{pa 4360 2380}%
\special{pa 4360 2356}%
\special{pa 4382 2346}%
\special{pa 4316 2326}%
\special{fp}%
%
\special{pn 8}%
\special{pa 4710 2010}%
\special{pa 4586 1926}%
\special{fp}%
\special{sh 1}%
\special{pa 4586 1926}%
\special{pa 4630 1980}%
\special{pa 4630 1956}%
\special{pa 4652 1946}%
\special{pa 4586 1926}%
\special{fp}%
%
\special{pn 8}%
\special{pa 5160 1270}%
\special{pa 5036 1186}%
\special{fp}%
\special{sh 1}%
\special{pa 5036 1186}%
\special{pa 5080 1240}%
\special{pa 5080 1216}%
\special{pa 5102 1206}%
\special{pa 5036 1186}%
\special{fp}%
%
\special{pn 8}%
\special{pa 5180 2500}%
\special{pa 5116 2636}%
\special{fp}%
\special{sh 1}%
\special{pa 5116 2636}%
\special{pa 5164 2584}%
\special{pa 5140 2588}%
\special{pa 5126 2566}%
\special{pa 5116 2636}%
\special{fp}%
%
\special{pn 8}%
\special{pa 4920 2040}%
\special{pa 4772 2036}%
\special{fp}%
\special{sh 1}%
\special{pa 4772 2036}%
\special{pa 4838 2058}%
\special{pa 4824 2038}%
\special{pa 4838 2018}%
\special{pa 4772 2036}%
\special{fp}%
%
\special{pn 8}%
\special{pa 5970 2520}%
\special{pa 5906 2656}%
\special{fp}%
\special{sh 1}%
\special{pa 5906 2656}%
\special{pa 5954 2604}%
\special{pa 5930 2608}%
\special{pa 5916 2586}%
\special{pa 5906 2656}%
\special{fp}%
%
\special{pn 8}%
\special{pa 4760 2530}%
\special{pa 4696 2666}%
\special{fp}%
\special{sh 1}%
\special{pa 4696 2666}%
\special{pa 4744 2614}%
\special{pa 4720 2618}%
\special{pa 4706 2596}%
\special{pa 4696 2666}%
\special{fp}%
%
\special{pn 8}%
\special{pa 5370 1290}%
\special{pa 5222 1286}%
\special{fp}%
\special{sh 1}%
\special{pa 5222 1286}%
\special{pa 5288 1308}%
\special{pa 5274 1288}%
\special{pa 5288 1268}%
\special{pa 5222 1286}%
\special{fp}%
%
\special{pn 8}%
\special{pa 4620 2420}%
\special{pa 4472 2416}%
\special{fp}%
\special{sh 1}%
\special{pa 4472 2416}%
\special{pa 4538 2438}%
\special{pa 4524 2418}%
\special{pa 4538 2398}%
\special{pa 4472 2416}%
\special{fp}%
%
\special{pn 8}%
\special{pa 4700 2490}%
\special{pa 4570 2570}%
\special{fp}%
\special{sh 1}%
\special{pa 4570 2570}%
\special{pa 4638 2552}%
\special{pa 4616 2542}%
\special{pa 4616 2518}%
\special{pa 4570 2570}%
\special{fp}%
%
\special{pn 8}%
\special{pa 5120 2230}%
\special{pa 4990 2310}%
\special{fp}%
\special{sh 1}%
\special{pa 4990 2310}%
\special{pa 5058 2292}%
\special{pa 5036 2282}%
\special{pa 5036 2258}%
\special{pa 4990 2310}%
\special{fp}%
%
\special{pn 8}%
\special{pa 5800 1870}%
\special{pa 5670 1950}%
\special{fp}%
\special{sh 1}%
\special{pa 5670 1950}%
\special{pa 5738 1932}%
\special{pa 5716 1922}%
\special{pa 5716 1898}%
\special{pa 5670 1950}%
\special{fp}%
%
\special{pn 8}%
\special{pa 5526 2700}%
\special{pa 5526 2850}%
\special{fp}%
\special{sh 1}%
\special{pa 5526 2850}%
\special{pa 5546 2784}%
\special{pa 5526 2798}%
\special{pa 5506 2784}%
\special{pa 5526 2850}%
\special{fp}%
%
\special{pn 8}%
\special{pa 5590 2510}%
\special{pa 5526 2646}%
\special{fp}%
\special{sh 1}%
\special{pa 5526 2646}%
\special{pa 5574 2594}%
\special{pa 5550 2598}%
\special{pa 5536 2576}%
\special{pa 5526 2646}%
\special{fp}%
%
\special{pn 8}%
\special{pa 5480 2060}%
\special{pa 5350 2140}%
\special{fp}%
\special{sh 1}%
\special{pa 5350 2140}%
\special{pa 5418 2122}%
\special{pa 5396 2112}%
\special{pa 5396 2088}%
\special{pa 5350 2140}%
\special{fp}%
%
\special{pn 8}%
\special{pa 4970 1600}%
\special{pa 4846 1516}%
\special{fp}%
\special{sh 1}%
\special{pa 4846 1516}%
\special{pa 4890 1570}%
\special{pa 4890 1546}%
\special{pa 4912 1536}%
\special{pa 4846 1516}%
\special{fp}%
%
\special{pn 8}%
\special{pa 5190 1620}%
\special{pa 5042 1616}%
\special{fp}%
\special{sh 1}%
\special{pa 5042 1616}%
\special{pa 5108 1638}%
\special{pa 5094 1618}%
\special{pa 5108 1598}%
\special{pa 5042 1616}%
\special{fp}%
%
\special{pn 8}%
\special{pa 5490 1530}%
\special{pa 5370 1630}%
\special{fp}%
\special{sh 1}%
\special{pa 5370 1630}%
\special{pa 5434 1604}%
\special{pa 5412 1596}%
\special{pa 5408 1572}%
\special{pa 5370 1630}%
\special{fp}%
%
\special{pn 8}%
\special{pa 5230 1920}%
\special{pa 5110 2020}%
\special{fp}%
\special{sh 1}%
\special{pa 5110 2020}%
\special{pa 5174 1994}%
\special{pa 5152 1986}%
\special{pa 5148 1962}%
\special{pa 5110 2020}%
\special{fp}%
%
\special{pn 8}%
\special{pa 4800 2290}%
\special{pa 4680 2390}%
\special{fp}%
\special{sh 1}%
\special{pa 4680 2390}%
\special{pa 4744 2364}%
\special{pa 4722 2356}%
\special{pa 4718 2332}%
\special{pa 4680 2390}%
\special{fp}%
%
\special{pn 8}%
\special{pa 4220 2720}%
\special{pa 4080 2630}%
\special{dt 0.045}%
\special{sh 1}%
\special{pa 4080 2630}%
\special{pa 4126 2684}%
\special{pa 4126 2660}%
\special{pa 4148 2650}%
\special{pa 4080 2630}%
\special{fp}%
%
\special{pn 8}%
\special{pa 4220 2720}%
\special{pa 4080 2790}%
\special{dt 0.045}%
\special{sh 1}%
\special{pa 4080 2790}%
\special{pa 4150 2778}%
\special{pa 4128 2766}%
\special{pa 4132 2742}%
\special{pa 4080 2790}%
\special{fp}%
%
\special{pn 8}%
\special{pa 4220 2720}%
\special{pa 4220 2870}%
\special{dt 0.045}%
\special{sh 1}%
\special{pa 4220 2870}%
\special{pa 4240 2804}%
\special{pa 4220 2818}%
\special{pa 4200 2804}%
\special{pa 4220 2870}%
\special{fp}%
\put(39.7000,-27.7000){\makebox(0,0)[rt]{?}}%
\put(52.3000,-31.0000){\makebox(0,0)[rt]{$\frac{X}{\vert\vert X\vert\vert}$}}%
%
\special{pn 8}%
\special{pa 990 2730}%
\special{pa 2980 2730}%
\special{fp}%
%
\special{pn 8}%
\special{pa 990 2730}%
\special{pa 1990 1150}%
\special{fp}%
\put(21.0000,-20.5000){\makebox(0,0)[lt]{$\widetilde C$}}%
%
\special{pn 8}%
\special{ar 3210 2970 2590 480  3.6717604 3.7108483}%
\special{ar 3210 2970 2590 480  3.7343011 3.7733890}%
\special{ar 3210 2970 2590 480  3.7968418 3.8359297}%
\special{ar 3210 2970 2590 480  3.8593825 3.8984705}%
\special{ar 3210 2970 2590 480  3.9219232 3.9610112}%
\special{ar 3210 2970 2590 480  3.9844639 4.0235519}%
\special{ar 3210 2970 2590 480  4.0470047 4.0860926}%
\special{ar 3210 2970 2590 480  4.1095454 4.1486333}%
\special{ar 3210 2970 2590 480  4.1720861 4.2111740}%
\special{ar 3210 2970 2590 480  4.2346268 4.2737148}%
\special{ar 3210 2970 2590 480  4.2971675 4.3362555}%
\special{ar 3210 2970 2590 480  4.3597082 4.3987962}%
\special{ar 3210 2970 2590 480  4.4222490 4.4613369}%
\special{ar 3210 2970 2590 480  4.4847897 4.5238776}%
\special{ar 3210 2970 2590 480  4.5473304 4.5864183}%
\special{ar 3210 2970 2590 480  4.6098711 4.6112912}%
%
\special{pn 8}%
\special{pa 1468 2632}%
\special{pa 1440 3490}%
\special{fp}%
\special{sh 1}%
\special{pa 1440 3490}%
\special{pa 1462 3424}%
\special{pa 1442 3438}%
\special{pa 1422 3424}%
\special{pa 1440 3490}%
\special{fp}%
%
\special{pn 8}%
\special{pa 1800 2580}%
\special{pa 1772 3440}%
\special{fp}%
\special{sh 1}%
\special{pa 1772 3440}%
\special{pa 1794 3374}%
\special{pa 1774 3386}%
\special{pa 1754 3372}%
\special{pa 1772 3440}%
\special{fp}%
%
\special{pn 8}%
\special{pa 2090 2550}%
\special{pa 2062 3410}%
\special{fp}%
\special{sh 1}%
\special{pa 2062 3410}%
\special{pa 2084 3344}%
\special{pa 2064 3358}%
\special{pa 2044 3344}%
\special{pa 2062 3410}%
\special{fp}%
%
\special{pn 8}%
\special{pa 2360 2530}%
\special{pa 2332 3390}%
\special{fp}%
\special{sh 1}%
\special{pa 2332 3390}%
\special{pa 2354 3324}%
\special{pa 2334 3338}%
\special{pa 2314 3324}%
\special{pa 2332 3390}%
\special{fp}%
%
\special{pn 8}%
\special{pa 2660 2510}%
\special{pa 2632 3370}%
\special{fp}%
\special{sh 1}%
\special{pa 2632 3370}%
\special{pa 2654 3304}%
\special{pa 2634 3316}%
\special{pa 2614 3302}%
\special{pa 2632 3370}%
\special{fp}%
%
\special{pn 8}%
\special{pa 2960 2950}%
\special{pa 2830 2730}%
\special{dt 0.045}%
\special{sh 1}%
\special{pa 2830 2730}%
\special{pa 2848 2798}%
\special{pa 2858 2776}%
\special{pa 2882 2778}%
\special{pa 2830 2730}%
\special{fp}%
\put(29.4000,-30.1000){\makebox(0,0)[lt]{$\widetilde{\sigma}_{i_0}$}}%
\put(20.1000,-35.7000){\makebox(0,0)[lt]{$X$}}%
%
\special{pn 8}%
\special{pa 4940 2460}%
\special{pa 4800 2510}%
\special{fp}%
\special{sh 1}%
\special{pa 4800 2510}%
\special{pa 4870 2506}%
\special{pa 4850 2492}%
\special{pa 4856 2470}%
\special{pa 4800 2510}%
\special{fp}%
%
\special{pn 8}%
\special{pa 5490 2300}%
\special{pa 5350 2350}%
\special{fp}%
\special{sh 1}%
\special{pa 5350 2350}%
\special{pa 5420 2346}%
\special{pa 5400 2332}%
\special{pa 5406 2310}%
\special{pa 5350 2350}%
\special{fp}%
%
\special{pn 8}%
\special{pa 5900 2270}%
\special{pa 5760 2320}%
\special{fp}%
\special{sh 1}%
\special{pa 5760 2320}%
\special{pa 5830 2316}%
\special{pa 5810 2302}%
\special{pa 5816 2280}%
\special{pa 5760 2320}%
\special{fp}%
%
\special{pn 8}%
\special{pa 4480 2600}%
\special{pa 4350 2680}%
\special{fp}%
\special{sh 1}%
\special{pa 4350 2680}%
\special{pa 4418 2662}%
\special{pa 4396 2652}%
\special{pa 4396 2628}%
\special{pa 4350 2680}%
\special{fp}%
\end{picture}%
\hspace{1.5truecm}}

\vspace{0.2truecm}

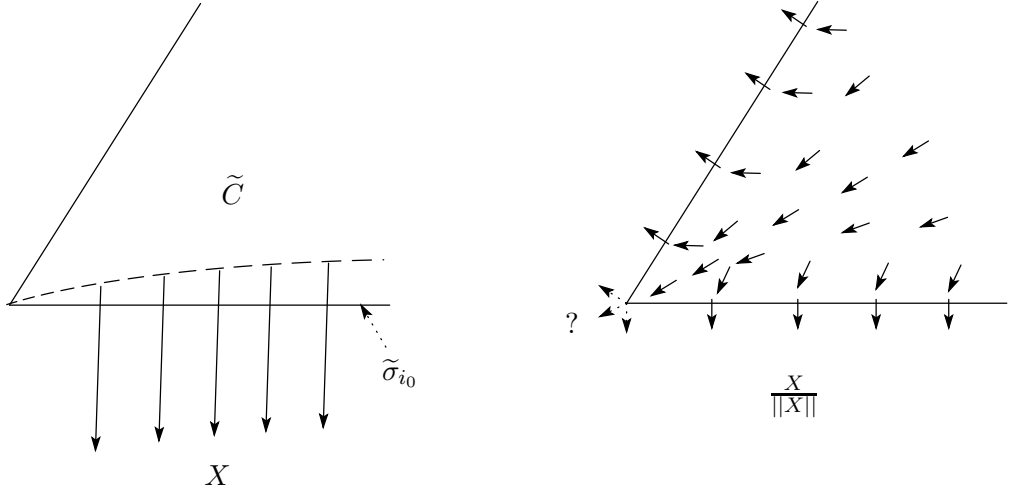
\centerline{{\bf Figure 1.}}

\vspace{0.3truecm}

\noindent
From the definition of $X$ and $(4.4)$, we have ${\rm grad}\,\rho=X$.  
For simplicity, set 
$\partial_i:=\frac{\partial}{\partial x_i}$ ($i=1,\cdots,r$).  
Then we have 
$$\begin{array}{l}
\displaystyle{(\partial_j\partial_k\rho)(w)=
\sum_{i\in I^R_+}\frac{m_i^+}
{\sinh^2({\rm arctanh}\,c_i^+-\beta_i(w))}
\beta_i(\partial_j)\beta_i(\partial_k)}\\
\hspace{2.5truecm}\displaystyle{
-\sum_{i\in I^R_-}\frac{m_i^-}
{\cosh^2({\rm arctanh}\,c_i^--\beta_i(w))}
\beta_i(\partial_j)\beta_i(\partial_k).}
\end{array}
\leqno{(4.5)}$$ 
As $w\to\partial\widetilde C$,
$\frac{1}{\sinh^2({\rm arctanh}\,c_i^+-\beta_i(w))}\to\infty$ for at least 
one $i\in I^R_+$ and 
$\frac{1}{\cosh^2({\rm arctanh}\,c_i^--\beta_i(w))}$\newline
$\leq1$ for all 
$i\in I^R_-$.  Hence we see that $\rho$ is downward convex 
on a sufficiently small collar neighborhood of $\partial\widetilde C$.  
Furthermore, since ${\rm codim}\,M={\rm rank}(G/K)$ and 
${\rm dim}(\mathfrak p_{\alpha}\cap\mathfrak p')\geq\frac12{\rm dim}\,
\mathfrak p_{\alpha}$ ($\alpha\in\triangle$) by the assumption, 
we have $I^R_+=I^R$, $m_i^+\geq m_i^-$ and $c_i^+=c_i^-$ ($i\in I^R_-$).  
From the relation $(4.5)$, we have 
$$\begin{array}{l}
\displaystyle{(\partial_j\partial_k\rho)(w)\geq
\sum_{i\in I^R\setminus I^R_-}\frac{m_i^+}
{\sinh^2({\rm arctanh}\,c_i^+-\beta_i(w))}
\beta_i(\partial_j)\beta_i(\partial_k)}\\
\hspace{2.5truecm}\displaystyle{
+\sum_{i\in I^R_-}\frac{4m_i^+}
{\sinh^22({\rm arctanh}\,c_i^+-\beta_i(w))}
\beta_i(\partial_j)\beta_i(\partial_k).}
\end{array}$$
Hence we see that $\rho$ is downward convex on the whole of $\widetilde C$.  
Also, it is clear that $\rho(w)\to\infty$ as $w\to\partial\widetilde C$ 
and that $\rho(tw)\to-\infty$ as $t\to\infty$ for each $w\in\widetilde C$.  
From these facts, $\rho$ and $X$ are as in Figure 2.  
Hence $t\mapsto\psi_t(\hat0)$ converges to a point $w_2$ of 
$\partial\widetilde C$ in a finite time $T$.   
Therefore $M$ is not minimal and 
the mean curvature flow $M_t$ collapses to the focal submanifold of $M$ 
through $\exp^{\perp}(w_2)$ in finite time.  
Thus the statement (i) is shown.  

Next we shall show the statement (iii) of Theorem A.  
Since $X$ is as in the second figure of Figure 2, 
we obtain the following fact:

\vspace{0.2truecm}

($\ast_1$) For each $w\in\partial\widetilde C$, there exists 
$w'\in\widetilde C$ such that the flow $\psi_t(w')$ converges 

\hspace{0.5truecm}to $w$.  

\vspace{0.3truecm}


\centerline{
\unitlength 0.1in
\begin{picture}( 56.1600, 25.5300)(  1.4400,-33.6300)
%
\special{pn 8}%
\special{pa 1164 2252}%
\special{pa 2866 2252}%
\special{fp}%
%
\special{pn 8}%
\special{pa 1164 3226}%
\special{pa 1164 810}%
\special{fp}%
\special{sh 1}%
\special{pa 1164 810}%
\special{pa 1144 878}%
\special{pa 1164 864}%
\special{pa 1184 878}%
\special{pa 1164 810}%
\special{fp}%
%
\special{pn 8}%
\special{ar 3224 122 2002 3194  1.7974986 2.9185575}%
%
\special{pn 8}%
\special{ar 1554 810 238 1102  0.0683985 1.5384335}%
%
\special{pn 8}%
\special{ar 1622 826 180 1128  1.8846315 1.9030082}%
\special{ar 1622 826 180 1128  1.9581384 1.9765151}%
\special{ar 1622 826 180 1128  2.0316453 2.0500220}%
\special{ar 1622 826 180 1128  2.1051522 2.1235289}%
\special{ar 1622 826 180 1128  2.1786590 2.1970358}%
\special{ar 1622 826 180 1128  2.2521659 2.2705427}%
\special{ar 1622 826 180 1128  2.3256728 2.3440495}%
\special{ar 1622 826 180 1128  2.3991797 2.4175564}%
\special{ar 1622 826 180 1128  2.4726866 2.4910633}%
\special{ar 1622 826 180 1128  2.5461935 2.5645702}%
\special{ar 1622 826 180 1128  2.6197004 2.6380771}%
\special{ar 1622 826 180 1128  2.6932073 2.7115840}%
\special{ar 1622 826 180 1128  2.7667142 2.7850909}%
\special{ar 1622 826 180 1128  2.8402211 2.8585978}%
\special{ar 1622 826 180 1128  2.9137280 2.9321047}%
\special{ar 1622 826 180 1128  2.9872348 3.0056116}%
\special{ar 1622 826 180 1128  3.0607417 3.0791185}%
\special{ar 1622 826 180 1128  3.1342486 3.1526254}%
%
\special{pn 8}%
\special{ar 2286 908 374 2036  6.2824821 6.2831853}%
\special{ar 2286 908 374 2036  0.0000000 1.5497222}%
%
\special{pn 8}%
\special{ar 2448 900 342 2302  2.0590642 2.0681448}%
\special{ar 2448 900 342 2302  2.0953866 2.1044672}%
\special{ar 2448 900 342 2302  2.1317089 2.1407895}%
\special{ar 2448 900 342 2302  2.1680313 2.1771119}%
\special{ar 2448 900 342 2302  2.2043537 2.2134343}%
\special{ar 2448 900 342 2302  2.2406760 2.2497566}%
\special{ar 2448 900 342 2302  2.2769984 2.2860790}%
\special{ar 2448 900 342 2302  2.3133208 2.3224013}%
\special{ar 2448 900 342 2302  2.3496431 2.3587237}%
\special{ar 2448 900 342 2302  2.3859655 2.3950461}%
\special{ar 2448 900 342 2302  2.4222878 2.4313684}%
\special{ar 2448 900 342 2302  2.4586102 2.4676908}%
\special{ar 2448 900 342 2302  2.4949326 2.5040131}%
\special{ar 2448 900 342 2302  2.5312549 2.5403355}%
\special{ar 2448 900 342 2302  2.5675773 2.5766579}%
\special{ar 2448 900 342 2302  2.6038996 2.6129802}%
\special{ar 2448 900 342 2302  2.6402220 2.6493026}%
\special{ar 2448 900 342 2302  2.6765444 2.6856250}%
\special{ar 2448 900 342 2302  2.7128667 2.7219473}%
\special{ar 2448 900 342 2302  2.7491891 2.7582697}%
\special{ar 2448 900 342 2302  2.7855114 2.7945920}%
\special{ar 2448 900 342 2302  2.8218338 2.8309144}%
\special{ar 2448 900 342 2302  2.8581562 2.8672368}%
\special{ar 2448 900 342 2302  2.8944785 2.9035591}%
\special{ar 2448 900 342 2302  2.9308009 2.9398815}%
\special{ar 2448 900 342 2302  2.9671233 2.9762038}%
\special{ar 2448 900 342 2302  3.0034456 3.0125262}%
\special{ar 2448 900 342 2302  3.0397680 3.0488486}%
\special{ar 2448 900 342 2302  3.0760903 3.0851709}%
\special{ar 2448 900 342 2302  3.1124127 3.1214933}%
\special{ar 2448 900 342 2302  3.1487351 3.1578156}%
%
\special{pn 8}%
\special{pa 1164 2252}%
\special{pa 1572 2042}%
\special{fp}%
\special{pa 856 1548}%
\special{pa 856 1548}%
\special{fp}%
%
\special{pn 8}%
\special{pa 1658 2010}%
\special{pa 2594 1548}%
\special{fp}%
%
\special{pn 8}%
\special{pa 1164 2252}%
\special{pa 2754 1864}%
\special{dt 0.045}%
%
\special{pn 8}%
\special{pa 1564 1912}%
\special{pa 1564 2148}%
\special{dt 0.045}%
%
\special{pn 8}%
\special{pa 2296 2942}%
\special{pa 2296 1978}%
\special{dt 0.045}%
%
\special{pn 13}%
\special{sh 1}%
\special{ar 1564 2156 10 10 0  6.28318530717959E+0000}%
\special{sh 1}%
\special{ar 1564 2156 10 10 0  6.28318530717959E+0000}%
%
\special{pn 13}%
\special{sh 1}%
\special{ar 2296 1986 10 10 0  6.28318530717959E+0000}%
\special{sh 1}%
\special{ar 2296 1986 10 10 0  6.28318530717959E+0000}%
%
\special{pn 13}%
\special{sh 1}%
\special{ar 1564 1912 10 10 0  6.28318530717959E+0000}%
\special{sh 1}%
\special{ar 1564 1912 10 10 0  6.28318530717959E+0000}%
%
\special{pn 13}%
\special{sh 1}%
\special{ar 2296 2942 10 10 0  6.28318530717959E+0000}%
\special{sh 1}%
\special{ar 2296 2942 10 10 0  6.28318530717959E+0000}%
%
\special{pn 13}%
\special{sh 1}%
\special{ar 1674 2132 10 10 0  6.28318530717959E+0000}%
\special{sh 1}%
\special{ar 1674 2132 10 10 0  6.28318530717959E+0000}%
%
\special{pn 8}%
\special{pa 4058 2252}%
\special{pa 5760 2252}%
\special{fp}%
%
\special{pn 8}%
\special{pa 4058 2252}%
\special{pa 4058 2252}%
\special{fp}%
\special{pa 4058 2252}%
\special{pa 5326 1572}%
\special{fp}%
%
\special{pn 13}%
\special{pa 4244 2212}%
\special{pa 3350 2424}%
\special{fp}%
\special{sh 1}%
\special{pa 3350 2424}%
\special{pa 3420 2428}%
\special{pa 3402 2412}%
\special{pa 3410 2388}%
\special{pa 3350 2424}%
\special{fp}%
%
\special{pn 8}%
\special{pa 4058 2252}%
\special{pa 5548 1928}%
\special{dt 0.045}%
%
\special{pn 13}%
\special{pa 5394 1962}%
\special{pa 5148 2018}%
\special{fp}%
\special{sh 1}%
\special{pa 5148 2018}%
\special{pa 5216 2022}%
\special{pa 5200 2006}%
\special{pa 5208 1984}%
\special{pa 5148 2018}%
\special{fp}%
%
\special{pn 13}%
\special{pa 5334 1872}%
\special{pa 5004 1864}%
\special{fp}%
\special{sh 1}%
\special{pa 5004 1864}%
\special{pa 5070 1886}%
\special{pa 5056 1866}%
\special{pa 5070 1846}%
\special{pa 5004 1864}%
\special{fp}%
%
\special{pn 13}%
\special{pa 5446 2034}%
\special{pa 5148 2172}%
\special{fp}%
\special{sh 1}%
\special{pa 5148 2172}%
\special{pa 5216 2162}%
\special{pa 5196 2150}%
\special{pa 5200 2126}%
\special{pa 5148 2172}%
\special{fp}%
\put(10.4400,-8.1000){\makebox(0,0)[rt]{${\Bbb R}$}}%
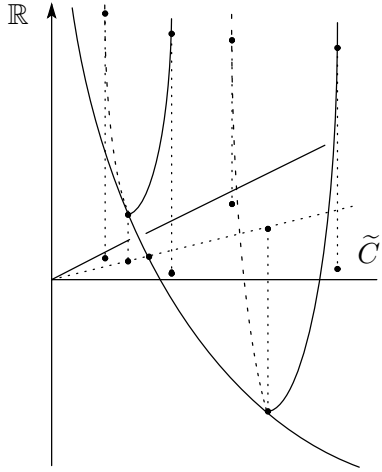
\put(24.9900,-33.6300){\makebox(0,0)[rt]{the graph of $\rho$}}%
\put(51.3000,-29.5800){\makebox(0,0)[rt]{$X$}}%
%
\special{pn 13}%
\special{sh 1}%
\special{ar 1792 964 10 10 0  6.28318530717959E+0000}%
\special{sh 1}%
\special{ar 1792 964 10 10 0  6.28318530717959E+0000}%
%
\special{pn 8}%
\special{pa 1792 972}%
\special{pa 1792 2220}%
\special{dt 0.045}%
%
\special{pn 8}%
\special{sh 1}%
\special{ar 1792 2212 10 10 0  6.28318530717959E+0000}%
\special{sh 1}%
\special{ar 1792 2212 10 10 0  6.28318530717959E+0000}%
%
\special{pn 13}%
\special{sh 1}%
\special{ar 1444 858 10 10 0  6.28318530717959E+0000}%
\special{sh 1}%
\special{ar 1444 858 10 10 0  6.28318530717959E+0000}%
%
\special{pn 13}%
\special{sh 1}%
\special{ar 1792 2220 10 10 0  6.28318530717959E+0000}%
\special{sh 1}%
\special{ar 1792 2220 10 10 0  6.28318530717959E+0000}%
%
\special{pn 8}%
\special{pa 1444 858}%
\special{pa 1444 2132}%
\special{dt 0.045}%
%
\special{pn 13}%
\special{sh 1}%
\special{ar 1444 2140 10 10 0  6.28318530717959E+0000}%
\special{sh 1}%
\special{ar 1444 2140 10 10 0  6.28318530717959E+0000}%
%
\special{pn 13}%
\special{sh 1}%
\special{ar 2660 1038 10 10 0  6.28318530717959E+0000}%
\special{sh 1}%
\special{ar 2660 1038 10 10 0  6.28318530717959E+0000}%
%
\special{pn 13}%
\special{sh 1}%
\special{ar 2108 998 10 10 0  6.28318530717959E+0000}%
\special{sh 1}%
\special{ar 2108 998 10 10 0  6.28318530717959E+0000}%
%
\special{pn 8}%
\special{pa 2660 1038}%
\special{pa 2660 2212}%
\special{dt 0.045}%
%
\special{pn 8}%
\special{pa 2108 1006}%
\special{pa 2108 1848}%
\special{dt 0.045}%
%
\special{pn 13}%
\special{sh 1}%
\special{ar 2660 2196 10 10 0  6.28318530717959E+0000}%
\special{sh 1}%
\special{ar 2660 2196 10 10 0  6.28318530717959E+0000}%
%
\special{pn 13}%
\special{sh 1}%
\special{ar 2108 1856 10 10 0  6.28318530717959E+0000}%
\special{sh 1}%
\special{ar 2108 1856 10 10 0  6.28318530717959E+0000}%
\put(28.8300,-20.1700){\makebox(0,0)[rt]{$\widetilde C$}}%
\put(56.8300,-20.1700){\makebox(0,0)[rt]{$\widetilde C$}}%
%
\special{pn 13}%
\special{pa 5080 2220}%
\special{pa 4812 2816}%
\special{fp}%
\special{sh 1}%
\special{pa 4812 2816}%
\special{pa 4858 2764}%
\special{pa 4834 2768}%
\special{pa 4820 2748}%
\special{pa 4812 2816}%
\special{fp}%
%
\special{pn 13}%
\special{pa 4738 2220}%
\special{pa 4426 2814}%
\special{fp}%
\special{sh 1}%
\special{pa 4426 2814}%
\special{pa 4474 2764}%
\special{pa 4450 2768}%
\special{pa 4438 2746}%
\special{pa 4426 2814}%
\special{fp}%
%
\special{pn 13}%
\special{pa 4884 1872}%
\special{pa 4348 1362}%
\special{fp}%
\special{sh 1}%
\special{pa 4348 1362}%
\special{pa 4382 1422}%
\special{pa 4386 1398}%
\special{pa 4410 1394}%
\special{pa 4348 1362}%
\special{fp}%
%
\special{pn 13}%
\special{pa 4586 2010}%
\special{pa 3922 1516}%
\special{fp}%
\special{sh 1}%
\special{pa 3922 1516}%
\special{pa 3964 1572}%
\special{pa 3964 1548}%
\special{pa 3986 1540}%
\special{pa 3922 1516}%
\special{fp}%
%
\special{pn 13}%
\special{pa 5650 1832}%
\special{pa 5470 1864}%
\special{fp}%
\special{sh 1}%
\special{pa 5470 1864}%
\special{pa 5540 1872}%
\special{pa 5522 1854}%
\special{pa 5532 1832}%
\special{pa 5470 1864}%
\special{fp}%
%
\special{pn 13}%
\special{pa 5684 1954}%
\special{pa 5540 2010}%
\special{fp}%
\special{sh 1}%
\special{pa 5540 2010}%
\special{pa 5608 2004}%
\special{pa 5590 1990}%
\special{pa 5594 1966}%
\special{pa 5540 2010}%
\special{fp}%
%
\special{pn 13}%
\special{pa 4994 2042}%
\special{pa 4722 2108}%
\special{fp}%
\special{sh 1}%
\special{pa 4722 2108}%
\special{pa 4790 2112}%
\special{pa 4774 2096}%
\special{pa 4782 2072}%
\special{pa 4722 2108}%
\special{fp}%
%
\special{pn 13}%
\special{pa 5088 2092}%
\special{pa 4790 2196}%
\special{fp}%
\special{sh 1}%
\special{pa 4790 2196}%
\special{pa 4860 2194}%
\special{pa 4840 2178}%
\special{pa 4846 2156}%
\special{pa 4790 2196}%
\special{fp}%
%
\special{pn 13}%
\special{pa 5028 1994}%
\special{pa 4722 2018}%
\special{fp}%
\special{sh 1}%
\special{pa 4722 2018}%
\special{pa 4790 2032}%
\special{pa 4774 2014}%
\special{pa 4786 1992}%
\special{pa 4722 2018}%
\special{fp}%
%
\special{pn 13}%
\special{pa 5548 2204}%
\special{pa 5280 2800}%
\special{fp}%
\special{sh 1}%
\special{pa 5280 2800}%
\special{pa 5326 2748}%
\special{pa 5302 2752}%
\special{pa 5288 2732}%
\special{pa 5280 2800}%
\special{fp}%
%
\special{pn 13}%
\special{pa 5250 1686}%
\special{pa 4714 1174}%
\special{fp}%
\special{sh 1}%
\special{pa 4714 1174}%
\special{pa 4748 1234}%
\special{pa 4752 1212}%
\special{pa 4776 1206}%
\special{pa 4714 1174}%
\special{fp}%
%
\special{pn 13}%
\special{pa 4356 2220}%
\special{pa 3820 2796}%
\special{fp}%
\special{sh 1}%
\special{pa 3820 2796}%
\special{pa 3880 2760}%
\special{pa 3856 2756}%
\special{pa 3850 2734}%
\special{pa 3820 2796}%
\special{fp}%
%
\special{pn 13}%
\special{pa 4322 2156}%
\special{pa 3496 1880}%
\special{fp}%
\special{sh 1}%
\special{pa 3496 1880}%
\special{pa 3554 1920}%
\special{pa 3548 1898}%
\special{pa 3566 1882}%
\special{pa 3496 1880}%
\special{fp}%
\end{picture}%
\hspace{3truecm}}

\vspace{0.3truecm}

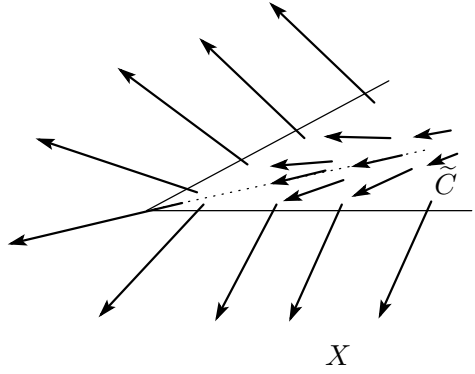
\centerline{{\bf Figure 2.}}

\vspace{0.3truecm}

\noindent
Now we shall show the following statement:

\vspace{0.2truecm}

\noindent
($\ast_2$) For any $w\in\partial\widetilde C$, the set 
$\{w'\in\widetilde C\,\vert\,{\rm the}\,\,{\rm flow}\,\,\psi_t(w')
\,\,{\rm converges}\,\,{\rm to}\,\,w\}$ 
is equal to 

\hspace{0.24truecm}the image of a flow of $X$.  

\vspace{0.2truecm}

\noindent
That is, we shall show that the situation as in Figure 3 cannot happen.  
Let $W$ be the real Coxeter group of $\widetilde M$ at $\hat 0$, that is, 
the group generated by the reflections with respect to 
the (real) hyperplanes ${\it l}_i$'s ($i\in I^R_+$) in $\mathfrak b$ 
containing $\widetilde{\sigma}_i$.  This group $W$ is a finite Coxeter 
group.  Set $V:={\rm Span}\{\beta_i^{\sharp}\,\vert\,i\in I^R_+\}$ and 
$\widetilde C_V:=\widetilde C\cap V$ (see Figure 4).  

\vspace{0.3truecm}

\centerline{
\unitlength 0.1in
\begin{picture}( 23.8000, 14.0000)( 18.3000,-28.0000)
%
\special{pn 8}%
\special{pa 2220 2800}%
\special{pa 4210 2800}%
\special{fp}%
%
\special{pn 8}%
\special{pa 2220 2800}%
\special{pa 3170 1400}%
\special{fp}%
%
\special{pn 8}%
\special{ar 800 940 2740 2170  0.6168532 1.0293852}%
%
\special{pn 8}%
\special{pa 2230 2800}%
\special{pa 2260 2788}%
\special{pa 2288 2774}%
\special{pa 2318 2760}%
\special{pa 2346 2748}%
\special{pa 2376 2734}%
\special{pa 2406 2722}%
\special{pa 2436 2712}%
\special{pa 2466 2700}%
\special{pa 2496 2688}%
\special{pa 2526 2678}%
\special{pa 2556 2668}%
\special{pa 2586 2658}%
\special{pa 2616 2648}%
\special{pa 2648 2638}%
\special{pa 2678 2630}%
\special{pa 2710 2622}%
\special{pa 2740 2612}%
\special{pa 2772 2604}%
\special{pa 2802 2598}%
\special{pa 2834 2590}%
\special{pa 2864 2584}%
\special{pa 2896 2578}%
\special{pa 2928 2572}%
\special{pa 2960 2566}%
\special{pa 2990 2562}%
\special{pa 3022 2556}%
\special{pa 3054 2552}%
\special{pa 3086 2550}%
\special{pa 3118 2546}%
\special{pa 3150 2542}%
\special{pa 3182 2540}%
\special{pa 3214 2538}%
\special{pa 3218 2538}%
\special{sp}%
%
\special{pn 20}%
\special{sh 1}%
\special{ar 2870 2360 10 10 0  6.28318530717959E+0000}%
\special{sh 1}%
\special{ar 2870 2360 10 10 0  6.28318530717959E+0000}%
%
\special{pn 20}%
\special{sh 1}%
\special{ar 2980 2560 10 10 0  6.28318530717959E+0000}%
\special{sh 1}%
\special{ar 2980 2560 10 10 0  6.28318530717959E+0000}%
%
\special{pn 8}%
\special{pa 2690 2500}%
\special{pa 2630 2550}%
\special{fp}%
\special{sh 1}%
\special{pa 2630 2550}%
\special{pa 2694 2524}%
\special{pa 2672 2516}%
\special{pa 2668 2492}%
\special{pa 2630 2550}%
\special{fp}%
%
\special{pn 8}%
\special{pa 2830 2590}%
\special{pa 2750 2600}%
\special{fp}%
\special{sh 1}%
\special{pa 2750 2600}%
\special{pa 2820 2612}%
\special{pa 2804 2594}%
\special{pa 2814 2572}%
\special{pa 2750 2600}%
\special{fp}%
\put(32.6000,-24.8000){\makebox(0,0)[lt]{$\psi_t(w')$}}%
\put(30.5000,-20.6000){\makebox(0,0)[lt]{$\psi_t({\bar w}')$}}%
\put(29.3000,-26.0000){\makebox(0,0)[lt]{$w'$}}%
\put(28.2000,-23.6000){\makebox(0,0)[rb]{${\bar w}'$}}%
\put(37.6000,-18.2000){\makebox(0,0)[lt]{$\widetilde C$}}%
\end{picture}%
\hspace{2truecm}}

\vspace{0.5truecm}

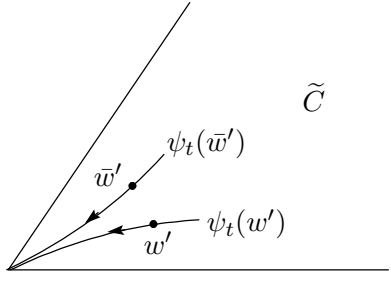
\centerline{{\bf Figure 3.}}

\vspace{1truecm}

\centerline{
\unitlength 0.1in
\begin{picture}( 32.0000, 14.7000)( 14.7000,-22.0000)
%
\special{pn 8}%
\special{pa 1800 2200}%
\special{pa 3610 2200}%
\special{fp}%
\special{pa 3610 2200}%
\special{pa 4200 1790}%
\special{fp}%
%
\special{pn 8}%
\special{pa 1800 2200}%
\special{pa 2600 1800}%
\special{fp}%
%
\special{pn 8}%
\special{pa 2710 2110}%
\special{pa 3730 2110}%
\special{fp}%
%
\special{pn 8}%
\special{pa 2700 2110}%
\special{pa 3470 1800}%
\special{fp}%
%
\special{pn 8}%
\special{pa 2700 2110}%
\special{pa 2700 1020}%
\special{fp}%
%
\special{pn 8}%
\special{pa 2600 1800}%
\special{pa 2660 1800}%
\special{fp}%
%
\special{pn 8}%
\special{pa 2750 1800}%
\special{pa 4200 1800}%
\special{fp}%
%
\special{pn 8}%
\special{pa 2700 1030}%
\special{pa 3730 1030}%
\special{da 0.070}%
%
\special{pn 8}%
\special{pa 3720 1030}%
\special{pa 3720 2110}%
\special{da 0.070}%
%
\special{pn 8}%
\special{pa 2700 1030}%
\special{pa 3450 750}%
\special{ip}%
%
\special{pn 8}%
\special{pa 3470 1800}%
\special{pa 3470 750}%
\special{da 0.070}%
%
\special{pn 8}%
\special{pa 3720 1030}%
\special{pa 4200 740}%
\special{da 0.070}%
%
\special{pn 8}%
\special{pa 3470 740}%
\special{pa 4200 740}%
\special{da 0.070}%
%
\special{pn 8}%
\special{pa 4200 740}%
\special{pa 4200 1800}%
\special{da 0.070}%
%
\special{pn 8}%
\special{pa 2700 2110}%
\special{pa 3720 2110}%
\special{pa 4200 1790}%
\special{pa 3470 1800}%
\special{pa 3470 1790}%
\special{pa 2700 2110}%
\special{fp}%
%
\special{pn 8}%
\special{pa 2700 1030}%
\special{pa 2700 2100}%
\special{pa 3730 2110}%
\special{pa 4200 1790}%
\special{pa 4200 750}%
\special{pa 3470 740}%
\special{pa 3480 730}%
\special{pa 3480 730}%
\special{pa 2700 1030}%
\special{ip}%
%
\special{pn 8}%
\special{pa 1940 1840}%
\special{pa 2260 2120}%
\special{da 0.070}%
\special{sh 1}%
\special{pa 2260 2120}%
\special{pa 2224 2062}%
\special{pa 2220 2086}%
\special{pa 2198 2092}%
\special{pa 2260 2120}%
\special{fp}%
%
\special{pn 8}%
\special{pa 2490 1190}%
\special{pa 2700 1270}%
\special{da 0.070}%
\special{sh 1}%
\special{pa 2700 1270}%
\special{pa 2646 1228}%
\special{pa 2650 1252}%
\special{pa 2632 1266}%
\special{pa 2700 1270}%
\special{fp}%
%
\special{pn 8}%
\special{pa 4460 1390}%
\special{pa 3840 1920}%
\special{da 0.070}%
\special{sh 1}%
\special{pa 3840 1920}%
\special{pa 3904 1892}%
\special{pa 3882 1886}%
\special{pa 3878 1862}%
\special{pa 3840 1920}%
\special{fp}%
\put(18.9000,-17.9000){\makebox(0,0)[lb]{$V$}}%
\put(24.6000,-12.6000){\makebox(0,0)[rb]{$V^{\perp}$}}%
\put(44.0000,-13.4000){\makebox(0,0)[lb]{$\widetilde C_V$}}%
%
\special{pn 8}%
\special{pa 4630 840}%
\special{pa 3600 1400}%
\special{da 0.070}%
\special{sh 1}%
\special{pa 3600 1400}%
\special{pa 3668 1386}%
\special{pa 3648 1376}%
\special{pa 3650 1352}%
\special{pa 3600 1400}%
\special{fp}%
\put(46.7000,-9.1000){\makebox(0,0)[lb]{$\widetilde C$}}%
%
\special{pn 4}%
\special{pa 2760 2100}%
\special{pa 2794 2066}%
\special{fp}%
\special{pa 2820 2100}%
\special{pa 2896 2024}%
\special{fp}%
\special{pa 2880 2100}%
\special{pa 2998 1984}%
\special{fp}%
\special{pa 2940 2100}%
\special{pa 3098 1942}%
\special{fp}%
\special{pa 3000 2100}%
\special{pa 3200 1900}%
\special{fp}%
\special{pa 3060 2100}%
\special{pa 3302 1860}%
\special{fp}%
\special{pa 3120 2100}%
\special{pa 3402 1818}%
\special{fp}%
\special{pa 3180 2100}%
\special{pa 3490 1790}%
\special{fp}%
\special{pa 3240 2100}%
\special{pa 3550 1792}%
\special{fp}%
\special{pa 3300 2100}%
\special{pa 3608 1792}%
\special{fp}%
\special{pa 3360 2100}%
\special{pa 3668 1794}%
\special{fp}%
\special{pa 3420 2100}%
\special{pa 3726 1794}%
\special{fp}%
\special{pa 3480 2100}%
\special{pa 3786 1794}%
\special{fp}%
\special{pa 3540 2100}%
\special{pa 3846 1796}%
\special{fp}%
\special{pa 3600 2100}%
\special{pa 3904 1796}%
\special{fp}%
\special{pa 3660 2100}%
\special{pa 3964 1798}%
\special{fp}%
\special{pa 3720 2100}%
\special{pa 4022 1798}%
\special{fp}%
\special{pa 3850 2032}%
\special{pa 4082 1800}%
\special{fp}%
\special{pa 4048 1892}%
\special{pa 4142 1800}%
\special{fp}%
%
\special{pn 4}%
\special{pa 3470 740}%
\special{pa 2700 1030}%
\special{pa 2710 2110}%
\special{pa 3720 2100}%
\special{pa 4190 1790}%
\special{pa 4190 740}%
\special{pa 4190 740}%
\special{pa 3470 740}%
\special{fp}%
%
\special{pn 4}%
\special{pa 2720 2110}%
\special{pa 2710 2090}%
\special{fp}%
\special{pa 2766 2110}%
\special{pa 2710 1998}%
\special{fp}%
\special{pa 2810 2110}%
\special{pa 2708 1906}%
\special{fp}%
\special{pa 2854 2110}%
\special{pa 2708 1816}%
\special{fp}%
\special{pa 2900 2108}%
\special{pa 2706 1724}%
\special{fp}%
\special{pa 2944 2108}%
\special{pa 2706 1632}%
\special{fp}%
\special{pa 2990 2108}%
\special{pa 2706 1540}%
\special{fp}%
\special{pa 3034 2108}%
\special{pa 2704 1448}%
\special{fp}%
\special{pa 3078 2106}%
\special{pa 2704 1356}%
\special{fp}%
\special{pa 3124 2106}%
\special{pa 2702 1264}%
\special{fp}%
\special{pa 3168 2106}%
\special{pa 2702 1174}%
\special{fp}%
\special{pa 3214 2106}%
\special{pa 2700 1082}%
\special{fp}%
\special{pa 3258 2106}%
\special{pa 2718 1024}%
\special{fp}%
\special{pa 3302 2104}%
\special{pa 2756 1010}%
\special{fp}%
\special{pa 3348 2104}%
\special{pa 2794 996}%
\special{fp}%
\special{pa 3392 2104}%
\special{pa 2830 982}%
\special{fp}%
\special{pa 3436 2104}%
\special{pa 2868 968}%
\special{fp}%
\special{pa 3482 2102}%
\special{pa 2906 952}%
\special{fp}%
\special{pa 3526 2102}%
\special{pa 2944 938}%
\special{fp}%
\special{pa 3572 2102}%
\special{pa 2982 924}%
\special{fp}%
\special{pa 3616 2102}%
\special{pa 3020 910}%
\special{fp}%
\special{pa 3660 2102}%
\special{pa 3058 896}%
\special{fp}%
\special{pa 3706 2100}%
\special{pa 3096 882}%
\special{fp}%
\special{pa 3744 2086}%
\special{pa 3134 868}%
\special{fp}%
\special{pa 3776 2064}%
\special{pa 3172 854}%
\special{fp}%
\special{pa 3810 2040}%
\special{pa 3210 838}%
\special{fp}%
\special{pa 3844 2018}%
\special{pa 3248 824}%
\special{fp}%
\special{pa 3878 1996}%
\special{pa 3286 810}%
\special{fp}%
\special{pa 3912 1974}%
\special{pa 3324 796}%
\special{fp}%
\special{pa 3946 1952}%
\special{pa 3362 782}%
\special{fp}%
%
\special{pn 4}%
\special{pa 3980 1930}%
\special{pa 3398 768}%
\special{fp}%
\special{pa 4014 1908}%
\special{pa 3436 754}%
\special{fp}%
\special{pa 4048 1884}%
\special{pa 3476 740}%
\special{fp}%
\special{pa 4082 1862}%
\special{pa 3520 740}%
\special{fp}%
\special{pa 4116 1840}%
\special{pa 3566 740}%
\special{fp}%
\special{pa 4150 1818}%
\special{pa 3610 740}%
\special{fp}%
\special{pa 4182 1796}%
\special{pa 3656 740}%
\special{fp}%
\special{pa 4190 1720}%
\special{pa 3700 740}%
\special{fp}%
\special{pa 4190 1630}%
\special{pa 3746 740}%
\special{fp}%
\special{pa 4190 1540}%
\special{pa 3790 740}%
\special{fp}%
\special{pa 4190 1450}%
\special{pa 3836 740}%
\special{fp}%
\special{pa 4190 1360}%
\special{pa 3880 740}%
\special{fp}%
\special{pa 4190 1270}%
\special{pa 3926 740}%
\special{fp}%
\special{pa 4190 1180}%
\special{pa 3970 740}%
\special{fp}%
\special{pa 4190 1090}%
\special{pa 4016 740}%
\special{fp}%
\special{pa 4190 1000}%
\special{pa 4060 740}%
\special{fp}%
\special{pa 4190 910}%
\special{pa 4106 740}%
\special{fp}%
\special{pa 4190 820}%
\special{pa 4150 740}%
\special{fp}%
\end{picture}%
\hspace{1.5truecm}}

\vspace{0.5truecm}

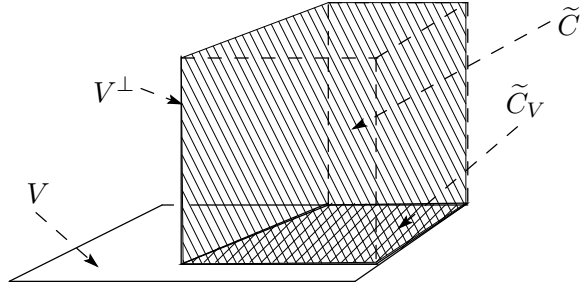
\centerline{{\bf Figure 4.}}


\vspace{0.5truecm}

\noindent
This space $V$ is $W$-invariant and $W$ acts trivially on the orthogonal 
complement $V^{\perp}$ of $V$.  Let $\{\phi_1,\cdots,\phi_{r'}\}$ be a base 
of the space of all $W$-invariant polynomial functions over $V$, where 
we note that $r'={\rm dim}\,V$.  
Set $\Phi:=(\phi_1,\cdots,\phi_{r'})$, which 
is a polynomial map from $V$ to ${\Bbb R}^{r'}$.  It is shown that $\Phi$ is 
a homeomorphism of the closure $\overline{\widetilde C_V}$ of 
$\widetilde C_V$ 
onto $\Phi(\overline{\widetilde C_V})$.  
Set $\xi_w(t):=\psi_t(w)$ and $\bar{\xi}_w(t):=\Phi(\psi_t(w))$, where 
$w\in\widetilde C_V$.  Let $(x_1,\cdots,x_{r'})$ be a Euclidean coordinate of 
$V$ and $(y_1,\cdots,y_{r'})$ the natural coordinate of ${\Bbb R}^{r'}$.  
Set $\xi^i_w(t):=x_i(\xi_w(t))$ and 
$\bar{\xi}^i_w(t):=y_i(\bar{\xi}_w(t))$ ($i=1,\cdots,r'$).  
Then we have 
$$\begin{array}{l}
\hspace{0.6truecm}
\displaystyle{(\bar{\xi}^i_w)'(t)=\langle{\rm grad}(y_i\circ\Phi)_{\xi_w(t)},
X_{\xi_w(t)}\rangle}\\
\displaystyle{=\sum_{j\in I^R_+}m_j^+\coth({\rm arctanh}\,c_j^+-\beta_j(\xi_w(t)))
\beta_j({\rm grad}(y_i\circ\Phi)_{\xi_w(t)})}\\
\hspace{0.6truecm}
\displaystyle{+\sum_{j\in I^R_-}m_j^-\tanh({\rm arctanh}\,c_j^--\beta_j
(\xi_w(t)))\beta_j({\rm grad}(y_i\circ\Phi)_{\xi_w(t)}).}
\end{array}$$
Let $f_i$ be the $W$-invariant $C^{\infty}$-function over $V$ such that 
$$\begin{array}{l}
\displaystyle{f_i(v):=
\sum_{j\in I^R_+}m_j^+\coth({\rm arctanh}\,c_j^+-\beta_j(v))
\beta_j({\rm grad}(y_i\circ\Phi)_v)}\\
\hspace{1.9truecm}
\displaystyle{+\sum_{j\in I^R_-}m_j^-\tanh({\rm arctanh}\,c_j^--\beta_j(v))
\beta_j({\rm grad}(y_i\circ\Phi)_v)}
\end{array}$$
for all $v\in\widetilde C_V$.  
It is easy to show that such a $W$-invariant $C^{\infty}$-function exists 
uniquely.  
According to the Schwarz's theorem in [S], we can describe $f_i$ as 
$f_i=Y_i\circ \Phi$ in terms of some $C^{\infty}$-function $Y_i$ over 
${\Bbb R}^{r'}$.  
Set $Y:=(Y_1,\cdots,Y_r)$, which is regarded as 
a $C^{\infty}$-vector field on ${\Bbb R}^{r'}$.  
Then we have $Y_{\Phi(w)}=\Phi_{\ast}(X_w)$ ($w\in\widetilde C_V$), that is, 
$Y\vert_{\Phi(\widetilde C_V)}=\Phi_{\ast}(X)$.  Also we can show that 
$Y\vert_{\partial\Phi(\widetilde C_V)}$ has no zero point.  
From these facts, we see that, for any $w\in\partial\widetilde C_V$, 
the set $\{w'\in\widetilde C_V\,\vert\,{\rm the}\,\,{\rm flow}\,\,\psi_t(w')
\,\,{\rm converges}\,\,{\rm to}\,\,w\}$ is equal to the image of a flow of 
$X$ (see Figure 5).  
In more general, we obtain the statement $(\ast_2)$ from this fact.  

\vspace{0.5truecm}

\centerline{
\unitlength 0.1in
\begin{picture}( 48.2900, 24.3400)(  8.5000,-34.0400)
%
\special{pn 13}%
\special{pa 1054 2344}%
\special{pa 1080 2326}%
\special{pa 1108 2310}%
\special{pa 1136 2294}%
\special{pa 1164 2278}%
\special{pa 1190 2262}%
\special{pa 1218 2246}%
\special{pa 1248 2234}%
\special{pa 1276 2218}%
\special{pa 1306 2204}%
\special{pa 1334 2192}%
\special{pa 1364 2180}%
\special{pa 1394 2168}%
\special{pa 1424 2156}%
\special{pa 1442 2154}%
\special{sp}%
%
\special{pn 13}%
\special{pa 1762 1874}%
\special{pa 1742 1898}%
\special{pa 1718 1922}%
\special{pa 1696 1944}%
\special{pa 1672 1966}%
\special{pa 1648 1986}%
\special{pa 1624 2008}%
\special{pa 1600 2028}%
\special{pa 1574 2048}%
\special{pa 1550 2068}%
\special{pa 1524 2086}%
\special{pa 1498 2104}%
\special{pa 1472 2124}%
\special{pa 1446 2142}%
\special{pa 1420 2160}%
\special{pa 1418 2160}%
\special{sp}%
%
\special{pn 13}%
\special{pa 1756 1878}%
\special{pa 1784 1862}%
\special{pa 1812 1846}%
\special{pa 1838 1830}%
\special{pa 1866 1814}%
\special{pa 1894 1798}%
\special{pa 1922 1782}%
\special{pa 1950 1768}%
\special{pa 1980 1754}%
\special{pa 2008 1740}%
\special{pa 2038 1728}%
\special{pa 2066 1714}%
\special{pa 2096 1704}%
\special{pa 2126 1694}%
\special{pa 2144 1688}%
\special{sp}%
%
\special{pn 13}%
\special{pa 2466 1408}%
\special{pa 2444 1432}%
\special{pa 2422 1456}%
\special{pa 2400 1478}%
\special{pa 2376 1500}%
\special{pa 2352 1520}%
\special{pa 2328 1542}%
\special{pa 2304 1562}%
\special{pa 2278 1582}%
\special{pa 2252 1600}%
\special{pa 2228 1620}%
\special{pa 2202 1640}%
\special{pa 2176 1658}%
\special{pa 2150 1676}%
\special{pa 2124 1694}%
\special{pa 2120 1696}%
\special{sp}%
%
\special{pn 13}%
\special{pa 1740 1714}%
\special{pa 1764 1692}%
\special{pa 1788 1672}%
\special{pa 1812 1652}%
\special{pa 1836 1630}%
\special{pa 1858 1608}%
\special{pa 1882 1586}%
\special{pa 1906 1564}%
\special{pa 1928 1542}%
\special{pa 1950 1518}%
\special{pa 1972 1496}%
\special{pa 1992 1472}%
\special{pa 2014 1448}%
\special{pa 2034 1422}%
\special{pa 2052 1396}%
\special{pa 2060 1386}%
\special{sp}%
%
\special{pn 13}%
\special{pa 2412 1174}%
\special{pa 2384 1186}%
\special{pa 2354 1200}%
\special{pa 2326 1214}%
\special{pa 2300 1232}%
\special{pa 2272 1248}%
\special{pa 2244 1264}%
\special{pa 2218 1280}%
\special{pa 2192 1300}%
\special{pa 2166 1318}%
\special{pa 2138 1336}%
\special{pa 2114 1354}%
\special{pa 2088 1374}%
\special{pa 2062 1394}%
\special{pa 2044 1408}%
\special{sp}%
%
\special{pn 13}%
\special{pa 1944 2104}%
\special{pa 1970 2086}%
\special{pa 1998 2070}%
\special{pa 2026 2054}%
\special{pa 2052 2036}%
\special{pa 2080 2022}%
\special{pa 2108 2006}%
\special{pa 2138 1992}%
\special{pa 2166 1978}%
\special{pa 2194 1964}%
\special{pa 2224 1950}%
\special{pa 2254 1940}%
\special{pa 2284 1928}%
\special{pa 2314 1918}%
\special{pa 2332 1912}%
\special{sp}%
%
\special{pn 13}%
\special{pa 2652 1632}%
\special{pa 2630 1656}%
\special{pa 2610 1680}%
\special{pa 2586 1702}%
\special{pa 2562 1724}%
\special{pa 2538 1744}%
\special{pa 2514 1766}%
\special{pa 2490 1786}%
\special{pa 2464 1806}%
\special{pa 2440 1826}%
\special{pa 2414 1844}%
\special{pa 2388 1864}%
\special{pa 2362 1882}%
\special{pa 2336 1900}%
\special{pa 2310 1918}%
\special{pa 2308 1920}%
\special{sp}%
%
\special{pn 13}%
\special{ar 2350 2396 652 388  3.1583234 4.0651315}%
%
\special{pn 13}%
\special{ar 1520 1438 330 388  0.7712003 1.8046062}%
%
\special{pn 13}%
\special{pa 2206 1632}%
\special{pa 2078 1754}%
\special{fp}%
\special{sh 1}%
\special{pa 2078 1754}%
\special{pa 2140 1722}%
\special{pa 2116 1718}%
\special{pa 2112 1694}%
\special{pa 2078 1754}%
\special{fp}%
%
\special{pn 13}%
\special{pa 1824 1826}%
\special{pa 1664 1928}%
\special{fp}%
\special{sh 1}%
\special{pa 1664 1928}%
\special{pa 1732 1908}%
\special{pa 1710 1900}%
\special{pa 1710 1876}%
\special{pa 1664 1928}%
\special{fp}%
%
\special{pn 13}%
\special{pa 1478 2112}%
\special{pa 1368 2202}%
\special{fp}%
\special{sh 1}%
\special{pa 1368 2202}%
\special{pa 1432 2176}%
\special{pa 1408 2168}%
\special{pa 1406 2144}%
\special{pa 1368 2202}%
\special{fp}%
%
\special{pn 13}%
\special{pa 1044 2346}%
\special{pa 884 2426}%
\special{fp}%
\special{sh 1}%
\special{pa 884 2426}%
\special{pa 954 2414}%
\special{pa 932 2402}%
\special{pa 934 2378}%
\special{pa 884 2426}%
\special{fp}%
%
\special{pn 13}%
\special{pa 2486 1796}%
\special{pa 2366 1918}%
\special{fp}%
\special{sh 1}%
\special{pa 2366 1918}%
\special{pa 2428 1884}%
\special{pa 2404 1880}%
\special{pa 2398 1856}%
\special{pa 2366 1918}%
\special{fp}%
%
\special{pn 13}%
\special{pa 2180 1958}%
\special{pa 2002 2050}%
\special{fp}%
\special{sh 1}%
\special{pa 2002 2050}%
\special{pa 2070 2036}%
\special{pa 2050 2026}%
\special{pa 2052 2002}%
\special{pa 2002 2050}%
\special{fp}%
%
\special{pn 13}%
\special{pa 1842 2152}%
\special{pa 1664 2316}%
\special{fp}%
\special{sh 1}%
\special{pa 1664 2316}%
\special{pa 1728 2286}%
\special{pa 1704 2280}%
\special{pa 1700 2256}%
\special{pa 1664 2316}%
\special{fp}%
%
\special{pn 13}%
\special{pa 1698 2366}%
\special{pa 1672 2610}%
\special{fp}%
\special{sh 1}%
\special{pa 1672 2610}%
\special{pa 1698 2546}%
\special{pa 1678 2556}%
\special{pa 1658 2542}%
\special{pa 1672 2610}%
\special{fp}%
%
\special{pn 13}%
\special{pa 1528 2620}%
\special{pa 1486 2854}%
\special{fp}%
\special{sh 1}%
\special{pa 1486 2854}%
\special{pa 1518 2792}%
\special{pa 1496 2800}%
\special{pa 1478 2784}%
\special{pa 1486 2854}%
\special{fp}%
%
\special{pn 13}%
\special{pa 1190 2508}%
\special{pa 1028 2690}%
\special{fp}%
\special{sh 1}%
\special{pa 1028 2690}%
\special{pa 1088 2654}%
\special{pa 1064 2650}%
\special{pa 1058 2628}%
\special{pa 1028 2690}%
\special{fp}%
%
\special{pn 13}%
\special{pa 2230 1266}%
\special{pa 2062 1348}%
\special{fp}%
\special{sh 1}%
\special{pa 2062 1348}%
\special{pa 2130 1336}%
\special{pa 2110 1324}%
\special{pa 2112 1300}%
\special{pa 2062 1348}%
\special{fp}%
%
\special{pn 13}%
\special{pa 1986 1480}%
\special{pa 1858 1622}%
\special{fp}%
\special{sh 1}%
\special{pa 1858 1622}%
\special{pa 1918 1586}%
\special{pa 1894 1582}%
\special{pa 1888 1560}%
\special{pa 1858 1622}%
\special{fp}%
%
\special{pn 13}%
\special{pa 1714 1744}%
\special{pa 1604 1856}%
\special{fp}%
\special{sh 1}%
\special{pa 1604 1856}%
\special{pa 1666 1822}%
\special{pa 1642 1818}%
\special{pa 1636 1794}%
\special{pa 1604 1856}%
\special{fp}%
%
\special{pn 13}%
\special{pa 1426 1816}%
\special{pa 1198 1786}%
\special{fp}%
\special{sh 1}%
\special{pa 1198 1786}%
\special{pa 1262 1814}%
\special{pa 1250 1792}%
\special{pa 1266 1774}%
\special{pa 1198 1786}%
\special{fp}%
%
\special{pn 13}%
\special{pa 1122 2040}%
\special{pa 850 2070}%
\special{fp}%
\special{sh 1}%
\special{pa 850 2070}%
\special{pa 920 2082}%
\special{pa 904 2064}%
\special{pa 914 2044}%
\special{pa 850 2070}%
\special{fp}%
%
\special{pn 13}%
\special{pa 2342 2112}%
\special{pa 2164 2242}%
\special{fp}%
\special{sh 1}%
\special{pa 2164 2242}%
\special{pa 2230 2220}%
\special{pa 2206 2210}%
\special{pa 2206 2186}%
\special{pa 2164 2242}%
\special{fp}%
%
\special{pn 13}%
\special{pa 2070 2294}%
\special{pa 1978 2580}%
\special{fp}%
\special{sh 1}%
\special{pa 1978 2580}%
\special{pa 2018 2522}%
\special{pa 1994 2528}%
\special{pa 1980 2510}%
\special{pa 1978 2580}%
\special{fp}%
%
\special{pn 13}%
\special{pa 2664 1990}%
\special{pa 2510 2192}%
\special{fp}%
\special{sh 1}%
\special{pa 2510 2192}%
\special{pa 2566 2152}%
\special{pa 2542 2150}%
\special{pa 2534 2128}%
\special{pa 2510 2192}%
\special{fp}%
%
\special{pn 13}%
\special{pa 1874 1408}%
\special{pa 1656 1530}%
\special{fp}%
\special{sh 1}%
\special{pa 1656 1530}%
\special{pa 1724 1516}%
\special{pa 1702 1504}%
\special{pa 1704 1480}%
\special{pa 1656 1530}%
\special{fp}%
%
\special{pn 13}%
\special{pa 1588 1612}%
\special{pa 1256 1572}%
\special{fp}%
\special{sh 1}%
\special{pa 1256 1572}%
\special{pa 1320 1600}%
\special{pa 1310 1578}%
\special{pa 1326 1560}%
\special{pa 1256 1572}%
\special{fp}%
%
\special{pn 13}%
\special{pa 2086 1204}%
\special{pa 1892 1246}%
\special{fp}%
\special{sh 1}%
\special{pa 1892 1246}%
\special{pa 1962 1252}%
\special{pa 1944 1234}%
\special{pa 1954 1212}%
\special{pa 1892 1246}%
\special{fp}%
%
\special{pn 13}%
\special{pa 1808 1276}%
\special{pa 1570 1224}%
\special{fp}%
\special{sh 1}%
\special{pa 1570 1224}%
\special{pa 1632 1258}%
\special{pa 1622 1236}%
\special{pa 1640 1220}%
\special{pa 1570 1224}%
\special{fp}%
%
\special{pn 13}%
\special{pa 2392 2284}%
\special{pa 2298 2610}%
\special{fp}%
\special{sh 1}%
\special{pa 2298 2610}%
\special{pa 2336 2550}%
\special{pa 2314 2558}%
\special{pa 2298 2540}%
\special{pa 2298 2610}%
\special{fp}%
%
\special{pn 13}%
\special{pa 2154 2498}%
\special{pa 2096 2804}%
\special{fp}%
\special{sh 1}%
\special{pa 2096 2804}%
\special{pa 2128 2742}%
\special{pa 2106 2752}%
\special{pa 2088 2734}%
\special{pa 2096 2804}%
\special{fp}%
%
\special{pn 13}%
\special{pa 1578 1398}%
\special{pa 1290 1358}%
\special{fp}%
\special{sh 1}%
\special{pa 1290 1358}%
\special{pa 1354 1386}%
\special{pa 1344 1366}%
\special{pa 1360 1348}%
\special{pa 1290 1358}%
\special{fp}%
%
\special{pn 13}%
\special{pa 1800 2640}%
\special{pa 1740 2926}%
\special{fp}%
\special{sh 1}%
\special{pa 1740 2926}%
\special{pa 1772 2864}%
\special{pa 1750 2874}%
\special{pa 1734 2856}%
\special{pa 1740 2926}%
\special{fp}%
%
\special{pn 13}%
\special{pa 1384 1480}%
\special{pa 1096 1428}%
\special{fp}%
\special{sh 1}%
\special{pa 1096 1428}%
\special{pa 1158 1460}%
\special{pa 1150 1438}%
\special{pa 1166 1420}%
\special{pa 1096 1428}%
\special{fp}%
\put(16.2100,-31.5900){\makebox(0,0)[lt]{$Y$}}%
\put(31.7100,-16.6300){\makebox(0,0)[lb]{$(\Phi\vert_{\bar{\widetilde C_V}})^{-1}$}}%
%
\special{pn 13}%
\special{pa 4146 2422}%
\special{pa 4912 1092}%
\special{fp}%
%
\special{pn 13}%
\special{pa 4154 2422}%
\special{pa 5680 2422}%
\special{fp}%
%
\special{pn 13}%
\special{pa 4980 1722}%
\special{pa 5006 1702}%
\special{pa 5030 1682}%
\special{pa 5054 1662}%
\special{pa 5080 1642}%
\special{pa 5106 1624}%
\special{pa 5134 1608}%
\special{pa 5160 1592}%
\special{pa 5190 1576}%
\special{pa 5218 1564}%
\special{pa 5248 1550}%
\special{pa 5278 1540}%
\special{pa 5284 1540}%
\special{sp}%
%
\special{pn 13}%
\special{pa 4762 2000}%
\special{pa 4740 2024}%
\special{pa 4718 2046}%
\special{pa 4696 2068}%
\special{pa 4672 2092}%
\special{pa 4650 2114}%
\special{pa 4626 2134}%
\special{pa 4602 2156}%
\special{pa 4576 2176}%
\special{pa 4552 2196}%
\special{pa 4526 2214}%
\special{pa 4500 2232}%
\special{pa 4474 2250}%
\special{pa 4446 2268}%
\special{pa 4418 2284}%
\special{pa 4392 2302}%
\special{pa 4364 2316}%
\special{pa 4334 2330}%
\special{pa 4306 2346}%
\special{pa 4278 2358}%
\special{pa 4248 2370}%
\special{pa 4218 2384}%
\special{pa 4190 2396}%
\special{pa 4162 2404}%
\special{sp}%
%
\special{pn 13}%
\special{pa 4758 1996}%
\special{pa 4782 1974}%
\special{pa 4806 1952}%
\special{pa 4832 1934}%
\special{pa 4856 1912}%
\special{pa 4882 1894}%
\special{pa 4908 1876}%
\special{pa 4934 1858}%
\special{pa 4960 1838}%
\special{pa 4988 1822}%
\special{pa 5014 1804}%
\special{pa 5042 1790}%
\special{pa 5070 1774}%
\special{pa 5098 1760}%
\special{pa 5128 1744}%
\special{pa 5156 1730}%
\special{pa 5186 1718}%
\special{pa 5216 1708}%
\special{pa 5246 1696}%
\special{pa 5276 1686}%
\special{pa 5306 1678}%
\special{pa 5336 1668}%
\special{pa 5368 1662}%
\special{pa 5382 1658}%
\special{sp}%
%
\special{pn 13}%
\special{ar 5602 2422 798 500  3.1415927 4.3783906}%
%
\special{pn 13}%
\special{pa 4980 1720}%
\special{pa 4954 1738}%
\special{pa 4928 1756}%
\special{pa 4900 1774}%
\special{pa 4874 1790}%
\special{pa 4846 1806}%
\special{pa 4818 1820}%
\special{pa 4786 1830}%
\special{pa 4756 1842}%
\special{pa 4726 1850}%
\special{pa 4696 1858}%
\special{pa 4664 1864}%
\special{pa 4632 1866}%
\special{pa 4600 1866}%
\special{pa 4568 1862}%
\special{pa 4538 1856}%
\special{pa 4506 1846}%
\special{pa 4482 1834}%
\special{sp}%
%
\special{pn 13}%
\special{pa 1680 2088}%
\special{pa 1708 2070}%
\special{pa 1734 2054}%
\special{pa 1762 2038}%
\special{pa 1790 2022}%
\special{pa 1818 2006}%
\special{pa 1846 1990}%
\special{pa 1874 1976}%
\special{pa 1904 1962}%
\special{pa 1932 1948}%
\special{pa 1960 1934}%
\special{pa 1990 1924}%
\special{pa 2020 1912}%
\special{pa 2050 1902}%
\special{pa 2068 1898}%
\special{sp}%
%
\special{pn 13}%
\special{ar 2086 2380 654 388  3.1583399 4.0659205}%
%
\special{pn 13}%
\special{pa 1918 1944}%
\special{pa 1740 2034}%
\special{fp}%
\special{sh 1}%
\special{pa 1740 2034}%
\special{pa 1808 2022}%
\special{pa 1786 2010}%
\special{pa 1790 1986}%
\special{pa 1740 2034}%
\special{fp}%
%
\special{pn 13}%
\special{pa 1578 2136}%
\special{pa 1400 2300}%
\special{fp}%
\special{sh 1}%
\special{pa 1400 2300}%
\special{pa 1464 2270}%
\special{pa 1440 2264}%
\special{pa 1436 2240}%
\special{pa 1400 2300}%
\special{fp}%
%
\special{pn 13}%
\special{pa 2406 1612}%
\special{pa 2384 1636}%
\special{pa 2362 1660}%
\special{pa 2340 1682}%
\special{pa 2316 1704}%
\special{pa 2292 1724}%
\special{pa 2268 1746}%
\special{pa 2244 1766}%
\special{pa 2218 1786}%
\special{pa 2194 1806}%
\special{pa 2168 1826}%
\special{pa 2142 1844}%
\special{pa 2116 1862}%
\special{pa 2090 1880}%
\special{pa 2064 1898}%
\special{pa 2062 1900}%
\special{sp}%
%
\special{pn 13}%
\special{pa 2240 1776}%
\special{pa 2120 1898}%
\special{fp}%
\special{sh 1}%
\special{pa 2120 1898}%
\special{pa 2182 1864}%
\special{pa 2158 1860}%
\special{pa 2152 1836}%
\special{pa 2120 1898}%
\special{fp}%
%
\special{pn 13}%
\special{ar 5446 2412 816 576  3.1415927 4.3760971}%
%
\special{pn 13}%
\special{pa 1434 2376}%
\special{pa 1368 2660}%
\special{fp}%
\special{sh 1}%
\special{pa 1368 2660}%
\special{pa 1402 2600}%
\special{pa 1380 2608}%
\special{pa 1364 2592}%
\special{pa 1368 2660}%
\special{fp}%
\put(10.6200,-34.0400){\makebox(0,0)[lt]{(The extension of $\Phi_{\ast}(X)$)}}%
\put(53.7300,-26.6000){\makebox(0,0)[rt]{Flows of $X\vert_{\widetilde C_V}$}}%
%
\special{pn 13}%
\special{pa 5170 1870}%
\special{pa 5374 1840}%
\special{fp}%
%
\special{pn 20}%
\special{pa 3070 1826}%
\special{pa 3850 1826}%
\special{fp}%
\special{sh 1}%
\special{pa 3850 1826}%
\special{pa 3782 1806}%
\special{pa 3796 1826}%
\special{pa 3782 1846}%
\special{pa 3850 1826}%
\special{fp}%
\put(25.6900,-11.5300){\makebox(0,0)[lb]{$\Phi(V)$}}%
\put(56.6900,-12.9500){\makebox(0,0)[rt]{$\bar{\widetilde C_V}$}}%
%
\special{pn 13}%
\special{ar 1040 2490 620 160  4.7256245 5.8715140}%
%
\special{pn 13}%
\special{ar 2170 2370 620 160  1.5840318 2.7299213}%
%
\special{pn 13}%
\special{ar 2160 2370 620 160  0.4116713 1.5575608}%
%
\special{pn 13}%
\special{pa 1070 2332}%
\special{pa 1088 2304}%
\special{pa 1104 2278}%
\special{pa 1122 2252}%
\special{pa 1142 2224}%
\special{pa 1160 2198}%
\special{pa 1178 2172}%
\special{pa 1198 2148}%
\special{pa 1216 2122}%
\special{pa 1236 2096}%
\special{pa 1256 2072}%
\special{pa 1278 2048}%
\special{pa 1300 2024}%
\special{pa 1322 2002}%
\special{pa 1344 1978}%
\special{pa 1368 1958}%
\special{pa 1394 1938}%
\special{pa 1420 1918}%
\special{pa 1444 1904}%
\special{sp}%
%
\special{pn 13}%
\special{pa 1792 1484}%
\special{pa 1774 1512}%
\special{pa 1758 1540}%
\special{pa 1740 1566}%
\special{pa 1722 1592}%
\special{pa 1704 1618}%
\special{pa 1684 1644}%
\special{pa 1666 1670}%
\special{pa 1646 1694}%
\special{pa 1626 1720}%
\special{pa 1606 1744}%
\special{pa 1584 1768}%
\special{pa 1564 1792}%
\special{pa 1540 1814}%
\special{pa 1518 1838}%
\special{pa 1494 1858}%
\special{pa 1470 1878}%
\special{pa 1444 1898}%
\special{pa 1418 1914}%
\special{sp}%
%
\special{pn 13}%
\special{pa 1994 970}%
\special{pa 1990 1002}%
\special{pa 1986 1034}%
\special{pa 1980 1066}%
\special{pa 1972 1096}%
\special{pa 1962 1126}%
\special{pa 1952 1156}%
\special{pa 1940 1186}%
\special{pa 1930 1216}%
\special{pa 1916 1246}%
\special{pa 1902 1274}%
\special{pa 1890 1304}%
\special{pa 1874 1332}%
\special{pa 1860 1360}%
\special{pa 1846 1390}%
\special{pa 1830 1418}%
\special{pa 1816 1446}%
\special{pa 1800 1474}%
\special{pa 1784 1498}%
\special{sp}%
\end{picture}%
\hspace{0.5truecm}}

\vspace{0.5truecm}

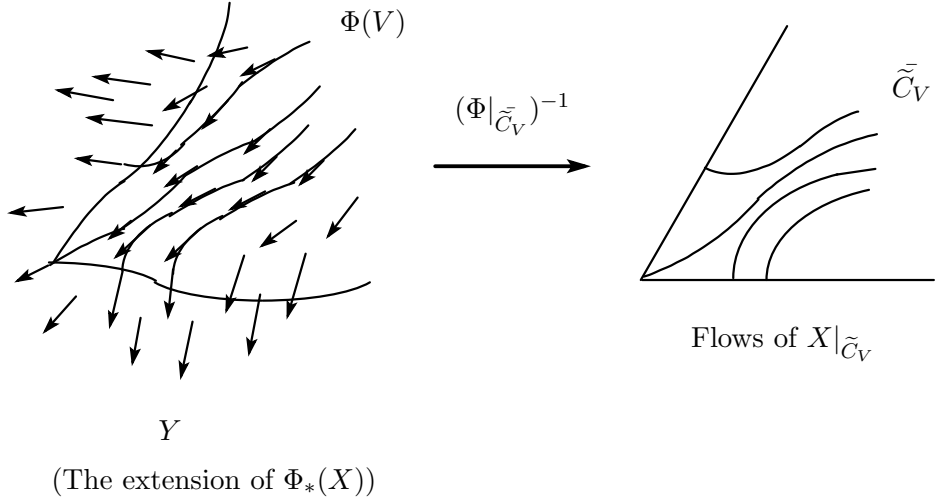
\centerline{{\bf Figure 5.}}

\vspace{0.5truecm}

\noindent
Take an arbitrary focal submanifold $F$ of $M$.  Let $\exp^{\perp}(w_1)$ be 
the only intersection point of $F$ and $\exp^{\perp}(\partial\widetilde C)$.  
According to the above fact $(\ast_2)$, 
the set of all parallel submanifolds of $M$ 
collapsing to $F$ along the mean curvature flow is a one-parameter 
$C^{\infty}$-family.  Thus the statement (iii) of Theorem A is shown.  
\hspace{2.5truecm}q.e.d.

\vspace{0.5truecm}

Next we prove the statement (ii) of Theorem A.  

\vspace{0.5truecm}

\noindent
{\it Proof of {\rm(ii)} of Theorem A.} 
Let $M$ and $F$ be as in (ii) of Theorem A.  
Since the natural fibration of $M$ onto $F$ is spherical, so is also 
the natural fibration of $\widetilde M$ onto $\widetilde F$.  
Hence $\widetilde F$ meets one of 
$(\partial\widetilde C\cap\beta_i^{-1}({\rm arctanh}\,c_i^+))^{\circ}$'s 
($i\in I^R_+$) (at one point).  Assume that $\widetilde F$ meets 
$(\partial\widetilde C\cap\beta_{i_0}^{-1}({\rm arctanh}
\,c_{i_0}^+))^{\circ}$.  
Let $u_0$ be the intersection point.  
Let $T$ be the explosion time of the flow $M_t$.  
Denote by $A^t$ (rep. $\widetilde A^t$) 
the shape tensor of $M_t$ (resp. $\widetilde M_t$), $H^t$ the mean curvature 
vector of $M_t$ and $\widetilde H^t$ the regularized mean curvature vector of 
$\widetilde M_t$.  We have 
$$\begin{array}{l}
\displaystyle{{\rm Spec}(\widetilde A^t_v)^{\bf c}\setminus\{0\}=
\left\{\left.\frac{\beta_i(v)}{{\rm arctanh}\,c_i^++j\pi\sqrt{-1}
-\beta_i(\psi_t(\hat0))}\,\right\vert\,i\in I^R_+,\,\,j\in{\Bbb Z}\right\}}\\
\hspace{2truecm}
\displaystyle{\cup\left\{
\left.\frac{\beta_i(v)}{{\rm arctanh}\,c_i^-+(j+\frac12)\pi\sqrt{-1}
-\beta_i(\psi_t(\hat0))}\,\right\vert\,i\in I^R_-,\,\,j\in{\Bbb Z}
\right\}}
\end{array}\leqno{(4.6)}$$
for each $v\in T^{\perp}_{\psi_t(\hat0)}\widetilde M_t
(=T_{\hat 0}\widetilde M)$.  
Since $\lim\limits_{t\to T-0}\psi_t(\hat0)=u_0\in
(\partial\widetilde C\cap\beta_{i_0}^{-1}({\rm arctanh}\,c_{i_0}^+))^{\circ}$, 
we have 
$\lim\limits_{t\to T-0}\beta_{i_0}(\psi_t(\hat0))
={\rm arctanh}\,c_{i_0}^+$ 
and $\lim\limits_{t\to T-0}\beta_i(\psi_t(\hat0))<{\rm arctanh}
\,c_i^+$ ($i\in I^R_+\setminus\{i_0\}$).  Hence we have 
$$\begin{array}{l}
\displaystyle{\lim_{t\to T-0}
\vert\vert(\widetilde A^t_v)^{\bf c}\vert\vert
^2_{\infty}(T-t)=\lim_{t\to T-0}\frac{\beta_{i_0}(v)^2}
{({\rm arctanh}\,c_{i_0}^+-\beta_{i_0}(\psi_t(\hat0)))^2}(T-t)}\\
\displaystyle{=\frac12\beta_{i_0}(v)^2\lim_{t\to T-0}\frac{1}
{({\rm arctanh}\,c_{i_0}^+-\beta_{i_0}(\psi_t(\hat0)))\beta_{i_0}
(\frac{d}{dt}\psi_t(\hat0))}.}
\end{array}\leqno{(4.7)}$$
Since $\frac{d}{dt}\psi_t(\hat0)=(\widetilde H^t)_{\psi_t(\hat0)}$, it follows 
from $(4.4)$ that 
$$\begin{array}{l}
\hspace{0.6truecm}\displaystyle{\lim_{t\to T-0}
({\rm arctanh}\,c_{i_0}^+-
\beta_{i_0}(\psi_t(\hat0))))\beta_{i_0}(\frac{d}{dt}\psi_t(\hat0))}\\
\displaystyle{=\lim_{t\to T-0}\left(\sum_{i\in I^R_+}m_i^+{\rm coth}
({\rm arctanh}\,c_i^+-\beta_i(\psi_t(\hat0)))\langle\beta_i^{\sharp},
\beta_{i_0}^{\sharp}\rangle({\rm arctanh}\,c_{i_0}^+
-\beta_{i_0}(\psi_t(\hat0)))\right.}\\
\hspace{0.6truecm}
\displaystyle{\left.
+\sum_{i\in I^R_-}m_i^-{\rm tanh}({\rm arctanh}\,c_i^--\beta_i
(\psi_t(\hat0)))\langle\beta_i^{\sharp},\beta_{i_0}^{\sharp}\rangle
({\rm arctanh}\,c_{i_0}^+-\beta_{i_0}(\psi_t(\hat0)))\right)}\\
\displaystyle{=m_{i_0}^+\langle\beta_{i_0}^{\sharp},
\beta_{i_0}^{\sharp}\rangle\lim_{t\to T-0}
{\rm coth}({\rm arctanh}\,c_{i_0}^+-\beta_{i_0}(\psi_t(\hat0)))
({\rm arctanh}\,c_{i_0}^+-\beta_{i_0}(\psi_t(\hat0)))}\\
\displaystyle{=m_{i_0}^+\langle\beta_{i_0}^{\sharp},
\beta_{i_0}^{\sharp}\rangle\lim_{t\to T-0}
\cosh^2({\rm arctanh}\,c_{i_0}^+-\beta_{i_0}(\psi_t(\hat0)))}\\
\displaystyle{=m_{i_0}^+\langle\beta_{i_0}^{\sharp},
\beta_{i_0}^{\sharp}\rangle,}
\end{array}$$
which together with $(4.7)$ deduces 
$$\lim_{t\to T-0}\vert\vert(\widetilde A^t_v)^{\bf c}
\vert\vert_{\infty}^2(T-t)
=\frac{\beta_{i_0}(v)^2}{2m_{i_0}^+\langle\beta_{i_0}^{\sharp},
\beta_{i_0}^{\sharp}\rangle}$$
and hence 
$$\lim_{t\to T-0}\mathop{\max}_{v\in S^{\perp}_{\psi_t(\hat0)}
\widetilde M_t}
\vert\vert(\widetilde A^t_v)^{\bf c}\vert\vert_{\infty}^2(T-t)
=\frac{1}{2m_{i_0}^+}.\leqno{(4.8)}$$
Thus $\widetilde M_t$ has type I singularity.  
Denote by $\exp_G$ the exponential map of $G$ and ${\rm Exp}$ 
the exponential map of $G/K$ at $eK$.  Also, denote by $S(1)$ the unit 
hypersphere in $\mathfrak b$ centered at $0$.  
Set $g_t:=\exp_G(\psi_t(\hat0))$ and $\bar v_t:=g_{t\ast}(v)$ for each 
$v\in S(1)$.  The relation 
$\bar v_t=(\pi\circ\phi)_{\ast\psi_t(\hat0)}(v)$ holds.  Since $M$ is proper 
complex equifocal and curvature-adapted by the assumption and since $M_t$ is a parallel 
submanifold of $M$, $M_t$ is also proper complex equifocal and 
curvature-adapted (see Lemma 3.4 of [Koi9]).  It is easy to show that 
$T_{{\rm Exp}(\psi_{t}(\hat0))}M_t=g_{t\ast}(\mathfrak m)$ and that 
$T_{{\rm Exp}(\psi_{t}(\hat0))}M_t=g_{t\ast}(\mathfrak m^R_0)
+\sum\limits_{i\in I^R}g_{t\ast}(\mathfrak m^R_i)$ is the common-eigenspace 
decomposition of $R(\cdot,\bar v_t)\bar v_t$'s ($v\in\mathfrak b$).  
In similar to $\beta_i$ ($i\in I^R$), $\lambda_i^+$ ($i\in I^R_+$) and 
$\lambda_i^-$ ($i\in I^R_-$), we define linear functions $\beta_i^t$ 
($i\in I^R$), $(\lambda_i^t)^+$ ($i\in I^R_+$) and $(\lambda_i^t)^-$ 
($i\in I^R_-$) on $T^{\perp}_{\psi_{t}(\hat0)}M_t=g_{t\ast}\mathfrak b$ by 
$$\begin{array}{l}
\displaystyle{R(\cdot,\bar v_t)\bar v_t\vert_{g_{t\ast}(\mathfrak m^R_i)}
=\beta_i^t(\bar v_t)^2{\rm id}\,\,\,\,(v\in\mathfrak b),}\\
\displaystyle{\{\lambda\in{\rm Spec}(A^t_{\bar v_t}
\vert_{g_{t\ast}(\mathfrak m^R_i)})\,\vert\,\vert\lambda\vert\,>\,
\vert\beta_i^t(\bar v_t)\vert\}=\{(\lambda_i^t)^+(\bar v_t)\}\,\,\,\,
(v\in\mathfrak b)}\\
\displaystyle{\{\lambda\in{\rm Spec}(A^t_{\bar v_t}
\vert_{g_{t\ast}(\mathfrak m^R_i)})\,\vert\,\vert\lambda\vert\,<\,
\vert\beta_i^t(\bar v_t)\vert\}=\{(\lambda_i^t)^-(\bar v_t)\}\,\,\,\,
(v\in\mathfrak b).}
\end{array}$$
It is clear that $\beta_i^t=\beta_i\circ g_{t\ast}^{-1}$ ($i\in I^R$).  
The values 
$\beta_i^t(\bar v_t)/(\lambda_i^t)^+(\bar v_t)$ ($i\in I^R_+$) and 
$(\lambda_i^t)^-(\bar v_t)/\beta_i^t(\bar v_t)$ ($i\in I^R_-$) are 
independent of the choice of $v\in\mathfrak b$.  
Denote by $(c_i^t)^+$ and $(c_i^t)^-$ these constants, respectively.  
If $i\in I^R_+\cap I^R_-$, then we have $(c_i^t)^+=(c_i^t)^-$.  
Hence we shall denote $(c_i^t)^+$ ($i\in I^R_+$) and $(c_i^t)^-$ 
($i\in I^R_-$) by $c_i^t$ for simplicity.  
In the sequel, we use this notation.  
The spectrum of $(\widetilde A^t_v)^{\bf c}$ other than zero is given by 
$$\begin{array}{l}
\displaystyle{{\rm Spec}\,(\widetilde A^t_v)^{\bf c}\setminus\{0\}
=\left\{\left.\frac{\beta_i(v)}{{\rm arctanh}\,c_i^t+j\pi\sqrt{-1}}\,\right\vert\,i\in I^R_+,
\,\,j\in{\Bbb Z}\right\}}\\
\hspace{3truecm}\displaystyle{\cup\left\{\left.\frac{\beta_i(v)}{{\rm arctanh}\,c_i^t
+(j+\frac12)\pi\sqrt{-1}}\,\right\vert\,i\in I^R_-,\,\,j\in{\Bbb Z}\right\}.}
\end{array}$$
On the other hand, we have 
$\displaystyle{\lim_{t\to T-0}\,\mathop{{\rm max}}_{v\in S(1)}
\vert(\lambda_{i_0}^t)^+(\bar v_t)\vert=\infty}$ and hence 
$\lim\limits_{t\to T-0}c_{i_0}^t=0$.  Also we have 
$\displaystyle{\lim\limits_{t\to T-0}\,\mathop{{\rm max}}_{v\in S(1)}
\vert(\lambda_i^t)^+(\bar v_t)\vert<\infty}$ and hence 
$\lim\limits_{t\to T-0}
\vert c_i^t\vert>0$ ($i\in I^R_+\setminus\{i_0\}$).  Therefore we obtain 
$$\begin{array}{l}
\hspace{0.6truecm}\displaystyle{\lim_{t\to T-0}(T-t)
\mathop{\max}_{v\in S(1)}\vert\vert(\widetilde A_v^t)^{\bf c}
\vert\vert_{\infty}^2
=\lim_{t\to T-0}(T-t)\mathop{\max}_{v\in S(1)}
\left(\frac{\beta_{i_0}(v)}{{\rm arctanh}\,c_{i_0}^t}\right)^2}\\
\displaystyle{=\mathop{\max}_{v\in S(1)}\beta_{i_0}(v)^2
\lim_{t\to T-0}\frac{T-t}
{{\rm arctanh}^2c_{i_0}^t}}\\
\displaystyle{=\mathop{\max}_{v\in S(1)}
\lim_{t\to T-0}\left(\frac{T-t}
{{\rm arctanh}^2(\beta_{i_0}(v)/(\lambda_{i_0}^t)^+(\bar v_t))}
\left(\frac{\beta_{i_0}(v)}{(\lambda_{i_0}^t)^+(\bar v_t)}\right)^2
(\lambda_{i_0}^t)^+(\bar v_t)^2\right)}\\
\displaystyle{=\lim_{t\to T-0}(T-t)\mathop{\max}_{v\in S(1)}
(\lambda_{i_0}^t)^+(\bar v_t)^2}\\
\displaystyle{=\lim_{t\to T-0}(T-t)\mathop{\max}_{v\in S(1)}
\vert\vert A^t_{\bar v_t}\vert\vert_{\infty}^2
=\lim_{t\to T-0}(T-t)\mathop{\max}_{v\in S(1)}
\vert\vert A^t_{\bar v_t}\vert\vert_{\infty}^2.}
\end{array}\leqno{(4.9)}$$
From this relation and $(4.8)$, we obtain 
$$\lim_{t\to T-0}(T-t)\mathop{\max}_{v\in S(1)}
\vert\vert A^t_{\bar v_t}\vert\vert_{\infty}^2
=\frac{1}{2m^+_{i_0}}\,<\,\infty.$$
Thus the mean curvature flow $M_t$ ($0\leq t<T$) has type I singularity.  
\hspace{1.5truecm}q.e.d.

\vspace{0.5truecm}

For each $S\subset I^R_+$, we set 
$$\begin{array}{l}
\displaystyle{\widetilde{\sigma}_S:=\{w\in\partial\widetilde C\,\vert\,
(\widetilde{\lambda}_i^+)_{\hat 0}(w)<1\,\,
(i\in I^R_+\setminus S)\,\,\&\,\,
(\widetilde{\lambda}_i^+)_{\hat 0}(w)=1\,\,(i\in S)\}}\\
\hspace{0.8truecm}
\displaystyle{=\{w\in\widetilde C\,\vert\,
\beta_i(w)<{\rm arctanh}\,c_i^+\,\,(i\in I^R_+\setminus S)\,\,\&\,\,
\beta_i(w)={\rm arctanh}\,c_i^+\,\,(i\in S)\},}
\end{array}$$
which is a simplex of $\widetilde C$.  
Take $w\in\widetilde{\sigma}_S$.  
Let $\widetilde w$ be the parallel 
normal vector field of $\widetilde M$ with $\widetilde w_{\hat0}=w$.  
Denote by $\eta_{\widetilde w}$ the end-point map for $\widetilde w$ and 
$\widetilde F_w:=\eta_{\widetilde w}(\widetilde M)$, which is a focal 
submanifold of $\widetilde M$.  
We have 
$$T_w\widetilde F_w
=\left(\mathop{\oplus}_{i\in I^R_+\setminus S}\mathop{\oplus}_{j\in{\Bbb Z}}
\eta_{\widetilde w\ast}(E_{ij}^+)\right)\oplus
\left(\mathop{\oplus}_{i\in I^R_-}\mathop{\oplus}_{j\in{\Bbb Z}}
\eta_{\widetilde w\ast}(E_{ij}^-)\right).$$
Denote by $\widetilde A^w$ the shape tensor of $\widetilde F_w$.  
In similar to $(4.3)$, we have 
$$\begin{array}{l}
\displaystyle{{\rm Tr}\,(\widetilde A^w_v)^{\bf c}
=\sum_{i\in I^R\setminus S}m_i^+\coth({\rm arctanh}\,c_i^+
-\beta_i(w))\beta_i(v)}\\
\hspace{1.8truecm}\displaystyle{+\sum_{i\in I^R_-}m_i^
-\tanh({\rm arctanh}\,c_i^+-\beta_i(w))\beta_i(v)\,\,(\in{\Bbb R})}
\end{array}
\leqno{(4.10)}$$
for any $v\in\mathfrak b$, where $\mathfrak b$ is regarded as a subspace of 
$T_w^{\perp}\widetilde F_w$.  
Set $L:=\widetilde M\cap T^{\perp}_w\widetilde F_w$, which is a focal leaf 
of $\widetilde M$.  For any $u\in L$, let $\mathfrak b_u$ be the section of 
$\widetilde M$ through $u$.  We can show 
$\displaystyle{(\widetilde H^w)_w\in
\mathop{\cap}_{u\in L}\mathfrak b_u}$.  
Hence, from $(4.10)$, the regularized mean curvature vector 
$\widetilde H^w$ of $\widetilde F_w$ exists and 
$(\widetilde H^w)_w$ is given by 
$$\begin{array}{l}
\displaystyle{(\widetilde H^w)_w=\sum_{i\in I^R\setminus S}
m_i^+\coth({\rm arctanh}\,c_i^+-\beta_i(w))\beta_i^{\sharp}}\\
\hspace{2.2truecm}\displaystyle{+\sum_{i\in I^R_-}m_i^-
\tanh({\rm arctanh}\,c_i^+-\beta_i(w))\beta_i^{\sharp}.}
\end{array}\leqno{(4.11)}$$
Define a vector field 
$X^{\widetilde{\sigma}_S}$ on $\widetilde{\sigma}_S$ by 
$X^{\widetilde{\sigma}_S}_w:=(\widetilde H^w)_w$ 
($w\in\widetilde{\sigma}_S$).  
This vector field $X^{\widetilde{\sigma}_S}$ is tangent to 
$\widetilde{\sigma}_S$.  
Let $\{\psi^{\widetilde{\sigma}_S}_t\}$ be the local one-parameter 
transformation group of $X^{\widetilde{\sigma}_S}$.  

\vspace{0.5truecm}

\noindent
{\it Proof of Theorem B.} 
First we shall show the statement (i) of Theorem B.  
Let $F$ be as in the satement (i) of Theorem B.  Set 
$\widetilde F:=(\pi\circ\phi)^{-1}(F)$.  
Since the lowest dimensional focal submanifold $F_{\it l}$ of $M$ is 
a one-point set by the assumption, we have $I^R_-=\emptyset$.  
Let $w_0$ be the intersection point of $\widetilde F$ and 
$\widetilde{\sigma}$.  Set 
$S_0:=\{i\in I^R_+(=I^R)\,\vert\,\beta_i(w_0)={\rm arctanh}\,c_i^+\}$.  
Since ${\rm dim}\,\widetilde{\sigma}\geq1$, we have 
$I^R\setminus S_0\not=\emptyset$.  
According to $(4.11)$, we have 
$$\begin{array}{l}
\displaystyle{(X^{\widetilde{\sigma}})_w=(\widetilde H^w)_w
=\sum_{i\in I^R\setminus S_0}
m_i^+\coth({\rm arctanh}\,c_i^+-\beta_i(w))\beta_i^{\sharp}}\\
\hspace{8.5truecm}\displaystyle{(w\in\widetilde{\sigma}).}
\end{array}\leqno{(4.12)}$$
We can show that $X^{\widetilde{\sigma}}$ is as in Figure 6 on a sufficiently 
small collar neighborhood of each maximal dimensional stratum of 
$\partial\widetilde{\sigma}$.  
Define a function $\rho_{\widetilde{\sigma}}$ over $\widetilde{\sigma}$ by 
$$\rho_{\widetilde{\sigma}}(w)
:=-\sum_{i\in I^R\setminus S_0}m_i^+\log\sinh
({\rm arctanh}\,c_i^+-\beta_i(w))\qquad(w\in\widetilde{\sigma}).$$
Easily we can show 
${\rm grad}\,\rho_{\widetilde{\sigma}}=X^{\widetilde{\sigma}}$.  
Let $(x_1,\cdots,x_{r''})$ be the Euclidean coordinate of 
$\displaystyle{\mathop{\cap}_{i\in S_0}\beta_i^{-1}({\rm arctanh}\,c_i^+)}$.  
For simplicity, set 
$\partial_i:=\frac{\partial}{\partial x_i}$ ($i=1,\cdots,r''$).  
Then we have 
$$(\partial_j\partial_k\rho_{\widetilde{\sigma}})(w)=
\sum_{i\in I^R\setminus S_0}\frac{m_i^+}
{\sinh^2({\rm arctanh}\,c_i^+-\beta_i(w))}
\beta_i(\partial_j)\beta_i(\partial_k).$$
Hence we see that $\rho_{\widetilde{\sigma}}$ is downward convex on 
$\widetilde{\sigma}$.  
Also, it is clear that $\rho_{\widetilde{\sigma}}(w)\to\infty$ as 
$w\to\partial\widetilde{\sigma}$ 
and that $\rho_{\widetilde{\sigma}}(tw)\to-\infty$ as $t\to\infty$ for each 
$w\in\widetilde{\sigma}$.  
From these facts, it follows that 
$\psi^{\widetilde{\sigma}}_t(w_0)$ converges to a point $w_1$ of 
$\partial\widetilde{\sigma}$ in a finite time.  
The mean curvature flow $F_t$ collapses to the focal submanifold of $M$ 
through $\exp^{\perp}(w_1) (\in\exp^{\perp}(\partial\widetilde{\sigma}))$.  
This completes the proof of the first-half part of the statement (i).  
The second-half part of the statement (i) is proved by imitating the proof of 
the statement (ii) of Theorem A.  

Next we shall show the statement (ii) of Theorem B.  
Set $V:={\rm Span}\{\beta_i^{\sharp}\,\vert\,i\in I^R_+\}$ and 
$\widetilde{\sigma}_V:=\widetilde{\sigma}\cap V$.  
Denote by $V_{\widetilde{\sigma}}$ be the affine subspace of $V$ containing 
$\widetilde{\sigma}_V$ as an open subset.  
Let $W_{\widetilde{\sigma}}$ be a finite Coxeter group generated by 
the reflections with respect to the (real) hyperplanes 
${\it l}^{\widetilde{\sigma}}_i$'s ($i\in I^R\setminus S_0$) in 
$V_{\widetilde{\sigma}}$ containing 
$\widetilde{\sigma}_i\cap V_{\widetilde{\sigma}}$.  
Let $\{\phi^{\widetilde{\sigma}}_1,\cdots,\phi^{\widetilde{\sigma}}_{r'}\}$ 
be a base of the space of all $W_{\widetilde{\sigma}}$-invariant polynomial 
functions over $V_{\widetilde{\sigma}}$, where 
we note that $r'={\rm dim}\,V_{\widetilde{\sigma}}$.  
Set $\Phi_{\widetilde{\sigma}}:=(\phi^{\widetilde{\sigma}}_1,\cdots,
\phi^{\widetilde{\sigma}}_{r'})$, which 
is a polynomial map from $V_{\widetilde{\sigma}}$ to ${\Bbb R}^{r'}$.  
It is shown that $\Phi_{\widetilde{\sigma}}$ is 
a homeomorphism of the closure $\overline{\widetilde{\sigma}_V}$ of 
$\widetilde{\sigma}_V$ onto 
$\Phi_{\widetilde{\sigma}}(\widetilde{\sigma}_V)$.  
Set $\xi_w(t):=\psi_t(w)$ and $\bar{\xi}_w(t):=\Phi_{\widetilde{\sigma}}
(\psi^{\widetilde{\sigma}}_t(w))$, where 
$w\in\widetilde{\sigma}_V$.  Let $(x_1,\cdots,x_{r'})$ be a Euclidean 
coordinate of $V_{\widetilde{\sigma}}$ and $(y_1,\cdots,y_{r'})$ 
the natural coordinate of ${\Bbb R}^{r'}$.  
Set $\xi^i_w(t):=x_i(\xi_w(t))$ and 
$\bar{\xi}^i_w(t):=y_i(\bar{\xi}_w(t))$ ($i=1,\cdots,r'$).  
Then we have 
$$\begin{array}{l}
\hspace{0.6truecm}
\displaystyle{(\bar{\xi}^i_w)'(t)=\langle{\rm grad}(y_i\circ
\Phi_{\widetilde{\sigma}})_{\xi_w(t)},X^{\widetilde{\sigma}}_{\xi_w(t)}
\rangle}\\
\displaystyle{=\sum_{j\in I^R\setminus S_0}m_j^+\coth({\rm arctanh}\,c_j^+-
\beta_j(\xi_w(t)))\beta_j({\rm grad}(y_i\circ
\Phi_{\widetilde{\sigma}})_{\xi_w(t)}).}
\end{array}$$
Let $f^{\widetilde{\sigma}}_i$ be the 
$W_{\widetilde{\sigma}}$-invariant $C^{\infty}$-function over 
$V_{\widetilde{\sigma}}$ such that 
$$\begin{array}{l}
\displaystyle{f^{\widetilde{\sigma}}_i(v):=
\sum_{j\in I^R\setminus S_0}m_j^+\coth({\rm arctanh}\,c_j^+-\beta_j(v))
\beta_j({\rm grad}(y_i\circ\Phi_{\widetilde{\sigma}})_v)}
\end{array}$$
for all $v\in\widetilde{\sigma}_V$.  
It is easy to show that such a $W_{\widetilde{\sigma}}$-invariant 
$C^{\infty}$-function exists uniquely.  
According to the Schwarz's theorem in [S], we can describe 
$f^{\widetilde{\sigma}}_i$ as 
$f^{\widetilde{\sigma}}_i=Y^{\widetilde{\sigma}}_i\circ
\Phi_{\widetilde{\sigma}}$ in terms of some 
$C^{\infty}$-function $Y^{\widetilde{\sigma}}_i$ over 
${\Bbb R}^{r'}$.  
Set $Y^{\widetilde{\sigma}}:=(Y^{\widetilde{\sigma}}_1,\cdots,
Y^{\widetilde{\sigma}}_r)$, which is regarded as 
a $C^{\infty}$-vector field on ${\Bbb R}^{r'}$.  
Then we have $Y^{\widetilde{\sigma}}_{\Phi_{\widetilde{\sigma}}(w)}
=(\Phi_{\widetilde{\sigma}})_{\ast}(X^{\widetilde{\sigma}}_w)$ 
($w\in\widetilde{\sigma}_V$), that is, 
$Y^{\widetilde{\sigma}}\vert_{\Phi_{\widetilde{\sigma}}(\widetilde{\sigma}_V)}
=(\Phi_{\widetilde{\sigma}})_{\ast}(X^{\widetilde{\sigma}})$.  
Also we can show that 
$Y^{\widetilde{\sigma}}\vert_{\partial\Phi_{\widetilde{\sigma}}
(\widetilde{\sigma}_V)}$ has no zero point.  
From these facts and the fact that $X^{\widetilde{\sigma}}$ is as in 
Figure 6 on a sufficiently small collar neighborhood of each maximal dimensional 
stratum of $\partial\widetilde{\sigma}_V$, we see that, for any 
$w\in\partial\widetilde{\sigma}_V$, 
the set $\{w'\in\widetilde{\sigma}_V\,\vert\,{\rm the}\,\,{\rm flow}\,\,
\psi^{\widetilde{\sigma}}_t(w')\,\,{\rm converges}\,\,{\rm to}\,\,w\}$ 
is equal to the image of a flow of $X^{\widetilde{\sigma}}$.  
In more general, the same fact holds 
for any $w\in\partial\widetilde{\sigma}$.  
From this fact, the statement (ii) of Theorem B follows.  
\hspace{4.3truecm}q.e.d.

\vspace{0.5truecm}

\centerline{
\unitlength 0.1in
\begin{picture}( 25.5000, 24.1900)( 22.4000,-38.3900)
%
\special{pn 8}%
\special{pa 2800 3000}%
\special{pa 4790 3000}%
\special{fp}%
%
\special{pn 8}%
\special{pa 2800 3000}%
\special{pa 3800 1420}%
\special{fp}%
\put(40.4000,-23.1000){\makebox(0,0)[rt]{$\widetilde{\sigma}$}}%
%
\special{pn 8}%
\special{ar 4790 3170 2230 370  3.6087401 3.6548939}%
\special{ar 4790 3170 2230 370  3.6825862 3.7287401}%
\special{ar 4790 3170 2230 370  3.7564324 3.8025862}%
\special{ar 4790 3170 2230 370  3.8302786 3.8764324}%
\special{ar 4790 3170 2230 370  3.9041247 3.9502786}%
\special{ar 4790 3170 2230 370  3.9779709 4.0241247}%
\special{ar 4790 3170 2230 370  4.0518170 4.0979709}%
\special{ar 4790 3170 2230 370  4.1256632 4.1718170}%
\special{ar 4790 3170 2230 370  4.1995093 4.2456632}%
\special{ar 4790 3170 2230 370  4.2733555 4.3195093}%
\special{ar 4790 3170 2230 370  4.3472016 4.3933555}%
\special{ar 4790 3170 2230 370  4.4210478 4.4672016}%
\special{ar 4790 3170 2230 370  4.4948939 4.5410478}%
\special{ar 4790 3170 2230 370  4.5687401 4.6148939}%
\special{ar 4790 3170 2230 370  4.6425862 4.6887401}%
%
\special{pn 8}%
\special{pa 4610 2810}%
\special{pa 4610 3700}%
\special{fp}%
\special{sh 1}%
\special{pa 4610 3700}%
\special{pa 4630 3632}%
\special{pa 4610 3646}%
\special{pa 4590 3632}%
\special{pa 4610 3700}%
\special{fp}%
%
\special{pn 8}%
\special{pa 4200 2820}%
\special{pa 4200 3710}%
\special{fp}%
\special{sh 1}%
\special{pa 4200 3710}%
\special{pa 4220 3642}%
\special{pa 4200 3656}%
\special{pa 4180 3642}%
\special{pa 4200 3710}%
\special{fp}%
%
\special{pn 8}%
\special{pa 3790 2850}%
\special{pa 3790 3740}%
\special{fp}%
\special{sh 1}%
\special{pa 3790 3740}%
\special{pa 3810 3672}%
\special{pa 3790 3686}%
\special{pa 3770 3672}%
\special{pa 3790 3740}%
\special{fp}%
%
\special{pn 8}%
\special{pa 3420 2880}%
\special{pa 3420 3770}%
\special{fp}%
\special{sh 1}%
\special{pa 3420 3770}%
\special{pa 3440 3702}%
\special{pa 3420 3716}%
\special{pa 3400 3702}%
\special{pa 3420 3770}%
\special{fp}%
%
\special{pn 8}%
\special{pa 3040 2950}%
\special{pa 3040 3840}%
\special{fp}%
\special{sh 1}%
\special{pa 3040 3840}%
\special{pa 3060 3772}%
\special{pa 3040 3786}%
\special{pa 3020 3772}%
\special{pa 3040 3840}%
\special{fp}%
\end{picture}%
\hspace{3truecm}}

\vspace{0.2truecm}

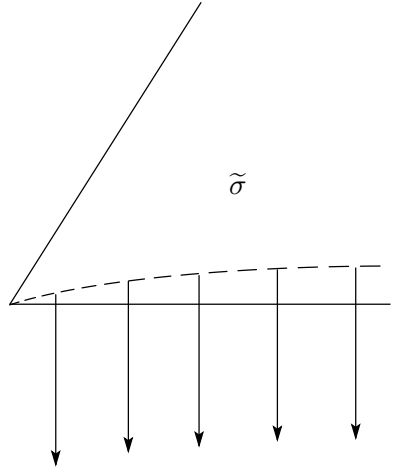
\centerline{{\bf Figure 6.}}

\vspace{0.4truecm}

We shall explain that, in the statement of Theorem B, we cannot weaken 
the condition that $F_{\it l}$ is a one-point set to the condition 
($\triangle'=\triangle$ and ${\rm dim}(\mathfrak p_{\alpha}\cap\mathfrak p')
\geq\frac12{\rm dim}\,\mathfrak p_{\alpha}$ ($\alpha\in\triangle$)) in 
the statement of Theorem A.  
Assume that $M$ satisfies the condition in the statement of Theorem A.  
Let $S_0$ be as above and $\widetilde{\sigma}
:=\widetilde{\sigma}_{S_0}$.  Define a function $\rho_{\widetilde{\sigma}}$ 
over $\widetilde{\sigma}$ by 
$$\begin{array}{l}
\displaystyle{\rho_{\widetilde{\sigma}}(w)
:=-\sum_{i\in I^R\setminus S_0}m_i^+\log\sinh
({\rm arctanh}\,c_i^+-\beta_i(w))}\\
\hspace{1.6truecm}\displaystyle{-\sum_{i\in I^R_-}m_i^-\log\cosh
({\rm arctanh}\,c_i^+-\beta_i(w))\qquad(w\in\widetilde{\sigma}).}
\end{array}$$
We have ${\rm grad}\,\rho_{\widetilde{\sigma}}=X^{\widetilde{\sigma}}$.  
Also, it follows from $m_i^+\geq m_i^-$ and $c_i^+=c_i^-$ ($i\in I^R_-$) that 
$$\begin{array}{l}
\displaystyle{(\partial_j\partial_k\rho_{\widetilde{\sigma}})(w)\geq
\sum_{i\in I^R\setminus(S_0\cup I^R_-)}\frac{m_i^+}
{\sinh^2({\rm arctanh}\,c_i^+-\beta_i(w))}
\beta_i(\partial_j)\beta_i(\partial_k)}\\
\hspace{2.5truecm}\displaystyle{
+\sum_{i\in(I^R\setminus S_0)\cap I^R_-}\frac{4m_i^+}
{\sinh^22({\rm arctanh}\,c_i^+-\beta_i(w))}
\beta_i(\partial_j)\beta_i(\partial_k)}\\
\hspace{2.5truecm}\displaystyle{
-\sum_{i\in S_0\cap I^R_-}\frac{m_i^+}
{\cosh^2({\rm arctanh}\,c_i^+-\beta_i(w))}
\beta_i(\partial_j)\beta_i(\partial_k)}.
\end{array}$$
Thus we cannot conclude whether $\rho$ is downward convex or not bacause of 
the existence of the third term in the right-hand side of this relation.  
Hence, in Theorem B, we cannot weaken the condition that $F_{\it l}$ is a one-point set 
to the condition in the statement of Theorem A.  

\section{Examples} 
Principal orbits of Hermann actions on a symmetric space $G/K$ of 
non-compact type are curvature-adapted isoparametric submanifolds and they 
have no focal point of non-Euclidean type on the ideal boundary of $G/K$.  
In particular, principal orbits of the isotropy action $K\curvearrowright G/K$ and 
those of Hermann actions $H\curvearrowright G/K$ 
as in Table 1 satisfy the additional conditions 
``${\rm codim}\,M={\rm rank}\,G/K$ and 
${\rm dim}(\mathfrak p_{\alpha}\cap\mathfrak p')\geq\frac12{\rm dim}\,
\mathfrak p_{\alpha}\,(\alpha\in\triangle)$" in the statement of Theorem A.  
In Table 1, $L$ is the fixed point group of 
$\theta\circ\tau$, where $\theta$ is a Cartan involution of $G$ with 
$({\rm Fix}\,\theta)_0\subset K\subset{\rm Fix}\,\theta$ and 
$\tau$ is an involution of $G$ with $({\rm Fix}\,\tau)_0\subset H\subset
{\rm Fix}\,\tau$.  
Then, for a Hermann action $H\curvearrowright G/K$, $F_{\it l}:=H(eK)$ 
is one of the lowest dimensional focal submanifolds of principal orbits of 
$H\curvearrowright G/K$.  The submanifolds $F_{\it l}$ and $F_{\it l}^{\perp}
:=\exp^{\perp}(T^{\perp}_{eK}H(eK))$ are reflective and hence 
they are symmetric spaces.  Explicitly they are described as 
$F_{\it l}=H/H\cap K$ and $F^{\perp}_{\it l}=L/H\cap K$, respectively (see Figure 7).  
In particular, in case of the 
isotropy action $K\curvearrowright G/K$, $F_{\it l}$ is a one-point set.  
Hence the principal orbits of the isotropy action satisfy the conditions 
in the statement of Theorem B.  

\vspace{0.1truecm}

$$\begin{tabular}{|c|c|c|c|}
\hline
{\scriptsize$H$} & {\scriptsize$G/K$} & {\scriptsize$F_{\it l}=H/H\cap K$} & 
{\scriptsize$F_{\it l}^{\perp}=L/H\cap K$}\\
\hline
{\scriptsize$SO^{\ast}(2n)$} & {\scriptsize$SU^{\ast}(2n)/Sp(n)$} & 
{\scriptsize $SO^{\ast}(2n)/U(n)$} & {\scriptsize $SL(n,{\Bbb C})/SU(n)$}\\
\hline
{\scriptsize$SO^{\ast}(2p)$} & {\scriptsize$SU(p,p)/S(U(p)\times U(p))$} & 
{\scriptsize $SO^{\ast}(2p)/U(p)$} & {\scriptsize $Sp(p,{\Bbb R})/U(p)$}\\
\hline
{\scriptsize$SO(n,{\Bbb C})$} & {\scriptsize$SL(n,{\Bbb C})/SU(n)$} & 
{\scriptsize$SO(n,{\Bbb C})/SO(n)$} & {\scriptsize$SL(n,{\Bbb R})/SO(n)$}\\
\hline
{\scriptsize$SU^{\ast}(2p)\cdot U(1)$} 
& {\scriptsize$Sp(p,p)/Sp(p)\times Sp(p)$} & 
{\scriptsize$SU^{\ast}(2p)/Sp(p)$} & 
{\scriptsize$Sp(p,{\Bbb C})/Sp(p)$}\\
\hline
{\scriptsize$SL(n,{\Bbb C})\cdot SO(2,{\Bbb C})$} 
& {\scriptsize$Sp(n,{\Bbb C})/Sp(n)$} & 
{\scriptsize$SL(n,{\Bbb C})/SU(n)$} & {\scriptsize$Sp(n,{\Bbb R})/U(n)$}\\
&&{\scriptsize$\times SO(2,{\Bbb C})/SO(2)$} & \\
\hline
{\scriptsize$Sp(1,3)$} 
& {\scriptsize$E_6^2/SU(6)\cdot SU(2)$} 
& {\scriptsize$Sp(1,3)/Sp(1)\times Sp(3)$} 
& {\scriptsize$F_4^4/Sp(3)\cdot Sp(1)$}\\
\hline
{\scriptsize$SU(1,5)\cdot SL(2,{\Bbb R})$} 
& {\scriptsize$E_6^{-14}/Spin(10)\cdot U(1)$} & 
{\scriptsize$SU(1,5)/S(U(1)\times U(5))$} & {\scriptsize$SO^{\ast}(10)/U(5)$}\\
&&{\scriptsize$\times SL(2,{\Bbb R})/SO(2)$}&\\
\hline
{\scriptsize$Sp(4,{\Bbb C})$} & {\scriptsize$E_6^{\bf c}/E_6$} & 
{\scriptsize$Sp(4,{\Bbb C})/Sp(4)$} & {\scriptsize$E_6^6/Sp(4)$}\\
\hline
{\scriptsize$SU(2,6)$} & {\scriptsize$E_7^{-5}/SO'(12)\cdot SU(2)$} & 
{\scriptsize$SU(2,6)/S(U(2)\times U(6))$} & 
{\scriptsize$E_6^2/SU(6)\cdot SU(2)$}\\
\hline
{\scriptsize$SL(8,{\Bbb C})$} & {\scriptsize$E_7^{\bf c}/E_7$} & 
{\scriptsize$SL(8,{\Bbb C})/SU(8)$} & {\scriptsize$E_7^7/SU(8)$}\\
\hline
{\scriptsize$SO(16,{\Bbb C})$} & 
{\scriptsize$E_8^{\bf c}/E_8$} & 
{\scriptsize$SO(16,{\Bbb C})/SO(16)$} & {\scriptsize$E_8^8/SO(16)$}\\
\hline
{\scriptsize$Sp(3,{\Bbb C})\cdot SL(2,{\Bbb C})$} 
& {\scriptsize$F_4^{\bf C}/F_4$} & 
{\scriptsize$Sp(3,{\Bbb C})/Sp(3)$} & 
{\scriptsize$F_4^4/Sp(3)\cdot Sp(1)$}\\
&&{\scriptsize$\times SL(2,{\Bbb C})/SU(2)$}&\\
\hline
{\scriptsize$SL(2,{\Bbb C})\times SL(2,{\Bbb C})$} 
& {\scriptsize$G_2^{\bf c}/G_2$} & 
{\scriptsize$SL(2,{\Bbb C})/SU(2)$} & {\scriptsize$G_2^2/SO(4)$}\\
&&{\scriptsize$\times SL(2,{\Bbb C})/SU(2)$}&\\
\hline
\end{tabular}$$

\vspace{0.2truecm}

\centerline{{\bf Table 1.}}


\vspace{0.5truecm}

\centerline{
\unitlength 0.1in
\begin{picture}( 49.0200, 27.9900)(  1.0200,-32.1900)
%
\special{pn 8}%
\special{ar 1928 2368 2718 212  4.7576159 5.9631853}%
%
\special{pn 13}%
\special{sh 1}%
\special{ar 3390 1888 10 10 0  6.28318530717959E+0000}%
\special{sh 1}%
\special{ar 3390 1888 10 10 0  6.28318530717959E+0000}%
%
\special{pn 8}%
\special{pa 3958 1906}%
\special{pa 3938 1932}%
\special{pa 3908 1944}%
\special{pa 3878 1952}%
\special{pa 3848 1960}%
\special{pa 3816 1966}%
\special{pa 3784 1972}%
\special{pa 3752 1976}%
\special{pa 3720 1978}%
\special{pa 3688 1982}%
\special{pa 3656 1984}%
\special{pa 3624 1986}%
\special{pa 3592 1988}%
\special{pa 3560 1990}%
\special{pa 3528 1990}%
\special{pa 3496 1992}%
\special{pa 3464 1992}%
\special{pa 3432 1992}%
\special{pa 3400 1992}%
\special{pa 3368 1992}%
\special{pa 3336 1990}%
\special{pa 3304 1988}%
\special{pa 3272 1988}%
\special{pa 3240 1986}%
\special{pa 3210 1984}%
\special{pa 3178 1982}%
\special{pa 3146 1978}%
\special{pa 3114 1976}%
\special{pa 3082 1972}%
\special{pa 3050 1966}%
\special{pa 3018 1962}%
\special{pa 2986 1956}%
\special{pa 2956 1950}%
\special{pa 2924 1942}%
\special{pa 2894 1934}%
\special{pa 2864 1922}%
\special{pa 2836 1906}%
\special{pa 2830 1878}%
\special{pa 2856 1860}%
\special{pa 2886 1848}%
\special{pa 2916 1840}%
\special{pa 2948 1834}%
\special{pa 2980 1826}%
\special{pa 3010 1822}%
\special{pa 3042 1818}%
\special{pa 3074 1816}%
\special{pa 3106 1812}%
\special{pa 3138 1810}%
\special{pa 3170 1808}%
\special{pa 3202 1806}%
\special{pa 3234 1806}%
\special{pa 3266 1804}%
\special{pa 3298 1804}%
\special{pa 3330 1804}%
\special{pa 3362 1804}%
\special{pa 3394 1804}%
\special{pa 3426 1806}%
\special{pa 3458 1806}%
\special{pa 3490 1808}%
\special{pa 3522 1810}%
\special{pa 3554 1812}%
\special{pa 3586 1814}%
\special{pa 3618 1816}%
\special{pa 3650 1818}%
\special{pa 3682 1822}%
\special{pa 3714 1826}%
\special{pa 3744 1832}%
\special{pa 3776 1836}%
\special{pa 3808 1842}%
\special{pa 3840 1848}%
\special{pa 3870 1856}%
\special{pa 3902 1864}%
\special{pa 3930 1878}%
\special{pa 3954 1898}%
\special{pa 3958 1906}%
\special{sp}%
%
\special{pn 8}%
\special{pa 3838 1030}%
\special{pa 3820 1054}%
\special{pa 3792 1070}%
\special{pa 3762 1080}%
\special{pa 3730 1088}%
\special{pa 3700 1096}%
\special{pa 3668 1102}%
\special{pa 3636 1108}%
\special{pa 3604 1112}%
\special{pa 3572 1116}%
\special{pa 3542 1122}%
\special{pa 3510 1126}%
\special{pa 3478 1128}%
\special{pa 3446 1130}%
\special{pa 3414 1132}%
\special{pa 3382 1134}%
\special{pa 3350 1136}%
\special{pa 3318 1138}%
\special{pa 3286 1140}%
\special{pa 3254 1142}%
\special{pa 3222 1142}%
\special{pa 3190 1142}%
\special{pa 3158 1142}%
\special{pa 3126 1142}%
\special{pa 3094 1142}%
\special{pa 3062 1140}%
\special{pa 3030 1140}%
\special{pa 2998 1138}%
\special{pa 2966 1134}%
\special{pa 2934 1132}%
\special{pa 2902 1128}%
\special{pa 2870 1124}%
\special{pa 2840 1120}%
\special{pa 2808 1114}%
\special{pa 2776 1106}%
\special{pa 2746 1096}%
\special{pa 2718 1082}%
\special{pa 2710 1054}%
\special{pa 2736 1034}%
\special{pa 2764 1022}%
\special{pa 2794 1012}%
\special{pa 2826 1004}%
\special{pa 2856 996}%
\special{pa 2888 990}%
\special{pa 2920 984}%
\special{pa 2952 980}%
\special{pa 2984 974}%
\special{pa 3014 972}%
\special{pa 3046 968}%
\special{pa 3078 964}%
\special{pa 3110 962}%
\special{pa 3142 960}%
\special{pa 3174 958}%
\special{pa 3206 956}%
\special{pa 3238 956}%
\special{pa 3270 954}%
\special{pa 3302 954}%
\special{pa 3334 952}%
\special{pa 3366 952}%
\special{pa 3398 952}%
\special{pa 3430 954}%
\special{pa 3462 954}%
\special{pa 3494 954}%
\special{pa 3526 956}%
\special{pa 3558 958}%
\special{pa 3590 962}%
\special{pa 3622 964}%
\special{pa 3654 968}%
\special{pa 3686 972}%
\special{pa 3716 978}%
\special{pa 3748 984}%
\special{pa 3780 992}%
\special{pa 3810 1002}%
\special{pa 3834 1022}%
\special{pa 3838 1030}%
\special{sp}%
%
\special{pn 8}%
\special{ar 2352 1888 472 1238  5.2665761 6.2831853}%
\special{ar 2352 1888 472 1238  0.0000000 0.2158999}%
%
\special{pn 8}%
\special{ar 3076 1888 880 1736  5.4964901 6.2831853}%
\special{ar 3076 1888 880 1736  0.0000000 0.1704352}%
%
\special{pn 8}%
\special{ar 2462 1882 354 1134  0.2839483 0.4717501}%
%
\special{pn 8}%
\special{ar 3114 1888 834 1686  0.2191680 0.3730137}%
%
\special{pn 8}%
\special{ar 2824 1888 566 1426  5.3744802 6.2831853}%
\special{ar 2824 1888 566 1426  0.0000000 0.0516644}%
%
\special{pn 8}%
\special{ar 2824 1888 566 1426  0.0875770 0.1938940}%
%
\special{pn 8}%
\special{ar 2824 1888 566 1438  0.2261168 0.4184180}%
%
\special{pn 13}%
\special{sh 1}%
\special{ar 3280 1048 10 10 0  6.28318530717959E+0000}%
\special{sh 1}%
\special{ar 3280 1048 10 10 0  6.28318530717959E+0000}%
%
\special{pn 8}%
\special{ar 4018 2760 2208 1332  3.6118369 4.0945117}%
%
\special{pn 8}%
\special{ar 2744 2678 4276 1002  5.0107308 5.2693944}%
%
\special{pn 8}%
\special{ar 2744 2678 4204 1008  4.8820866 4.9892178}%
%
\special{pn 8}%
\special{ar 2744 2678 4126 1008  4.7425313 4.8566809}%
%
\special{pn 8}%
\special{ar 5354 2666 920 914  3.5486241 4.3229023}%
%
\special{pn 8}%
\special{ar 2966 2118 2176 244  4.4267539 5.7014649}%
%
\special{pn 13}%
\special{ar 2336 2068 2744 198  5.1082287 5.8023466}%
%
\special{pn 13}%
\special{sh 1}%
\special{ar 3956 1906 10 10 0  6.28318530717959E+0000}%
\special{sh 1}%
\special{ar 3956 1906 10 10 0  6.28318530717959E+0000}%
%
\special{pn 8}%
\special{pa 4538 1640}%
\special{pa 4308 1932}%
\special{dt 0.045}%
\special{sh 1}%
\special{pa 4308 1932}%
\special{pa 4366 1892}%
\special{pa 4342 1890}%
\special{pa 4334 1866}%
\special{pa 4308 1932}%
\special{fp}%
%
\special{pn 8}%
\special{pa 2392 1590}%
\special{pa 2580 1888}%
\special{dt 0.045}%
\special{sh 1}%
\special{pa 2580 1888}%
\special{pa 2560 1820}%
\special{pa 2552 1842}%
\special{pa 2528 1842}%
\special{pa 2580 1888}%
\special{fp}%
%
\special{pn 8}%
\special{pa 3428 570}%
\special{pa 3216 838}%
\special{dt 0.045}%
\special{sh 1}%
\special{pa 3216 838}%
\special{pa 3274 798}%
\special{pa 3250 796}%
\special{pa 3242 772}%
\special{pa 3216 838}%
\special{fp}%
%
\special{pn 8}%
\special{pa 4174 906}%
\special{pa 3728 1348}%
\special{dt 0.045}%
\special{sh 1}%
\special{pa 3728 1348}%
\special{pa 3788 1314}%
\special{pa 3766 1310}%
\special{pa 3760 1286}%
\special{pa 3728 1348}%
\special{fp}%
%
\special{pn 8}%
\special{pa 4874 1496}%
\special{pa 4678 1864}%
\special{dt 0.045}%
\special{sh 1}%
\special{pa 4678 1864}%
\special{pa 4728 1814}%
\special{pa 4704 1816}%
\special{pa 4692 1796}%
\special{pa 4678 1864}%
\special{fp}%
%
\special{pn 13}%
\special{sh 1}%
\special{ar 3120 2898 10 10 0  6.28318530717959E+0000}%
\special{sh 1}%
\special{ar 3120 2898 10 10 0  6.28318530717959E+0000}%
%
\special{pn 8}%
\special{sh 1}%
\special{ar 3756 2860 10 10 0  6.28318530717959E+0000}%
\special{sh 1}%
\special{ar 3756 2860 10 10 0  6.28318530717959E+0000}%
\put(25.3200,-14.5900){\makebox(0,0)[rt]{$\exp^{\perp}(\mathfrak b)$}}%
\put(34.9900,-32.1900){\makebox(0,0)[rt]{$\exp^{\perp}(\mathfrak b)$}}%
\put(47.3300,-4.2000){\makebox(0,0)[rt]{$F_{\it l}=H(eK)=H/H\cap K$}}%
\put(48.2000,-7.5000){\makebox(0,0)[rt]{$M=H(gK)$}}%
\put(46.1500,-15.2100){\makebox(0,0)[rt]{$C$}}%
\put(56.6000,-13.5300){\makebox(0,0)[rt]{$F_{\it l}^{\perp}=L/H\cap K$}}%
\put(30.7500,-28.5800){\makebox(0,0)[rt]{$eK$}}%
\put(39.2300,-28.7000){\makebox(0,0)[rt]{$gK$}}%
%
\special{pn 8}%
\special{pa 4120 1676}%
\special{pa 3956 1900}%
\special{dt 0.045}%
\special{sh 1}%
\special{pa 3956 1900}%
\special{pa 4012 1858}%
\special{pa 3988 1858}%
\special{pa 3978 1834}%
\special{pa 3956 1900}%
\special{fp}%
%
\special{pn 8}%
\special{pa 3540 1646}%
\special{pa 3390 1882}%
\special{dt 0.045}%
\special{sh 1}%
\special{pa 3390 1882}%
\special{pa 3442 1836}%
\special{pa 3418 1836}%
\special{pa 3408 1814}%
\special{pa 3390 1882}%
\special{fp}%
\put(36.7200,-15.2100){\makebox(0,0)[rt]{$eK$}}%
\put(41.9800,-15.7000){\makebox(0,0)[rt]{$gK$}}%
\put(25.2500,-24.2300){\makebox(0,0)[rt]{in fact}}%
%
\special{pn 8}%
\special{ar 3712 1888 1250 1234  2.3089287 2.3185944}%
\special{ar 3712 1888 1250 1234  2.3475916 2.3572573}%
\special{ar 3712 1888 1250 1234  2.3862545 2.3959202}%
\special{ar 3712 1888 1250 1234  2.4249174 2.4345831}%
\special{ar 3712 1888 1250 1234  2.4635803 2.4732460}%
\special{ar 3712 1888 1250 1234  2.5022432 2.5119089}%
\special{ar 3712 1888 1250 1234  2.5409061 2.5505718}%
\special{ar 3712 1888 1250 1234  2.5795690 2.5892347}%
\special{ar 3712 1888 1250 1234  2.6182319 2.6278977}%
\special{ar 3712 1888 1250 1234  2.6568948 2.6665606}%
\special{ar 3712 1888 1250 1234  2.6955577 2.7052235}%
\special{ar 3712 1888 1250 1234  2.7342206 2.7438864}%
\special{ar 3712 1888 1250 1234  2.7728836 2.7825493}%
\special{ar 3712 1888 1250 1234  2.8115465 2.8212122}%
\special{ar 3712 1888 1250 1234  2.8502094 2.8598751}%
\special{ar 3712 1888 1250 1234  2.8888723 2.8985380}%
\special{ar 3712 1888 1250 1234  2.9275352 2.9372009}%
\special{ar 3712 1888 1250 1234  2.9661981 2.9758638}%
\special{ar 3712 1888 1250 1234  3.0048610 3.0145267}%
\special{ar 3712 1888 1250 1234  3.0435239 3.0531896}%
\special{ar 3712 1888 1250 1234  3.0821868 3.0918525}%
\special{ar 3712 1888 1250 1234  3.1208497 3.1305155}%
%
\special{pn 8}%
\special{pa 2878 2796}%
\special{pa 2894 2810}%
\special{dt 0.045}%
\special{sh 1}%
\special{pa 2894 2810}%
\special{pa 2856 2752}%
\special{pa 2854 2776}%
\special{pa 2830 2782}%
\special{pa 2894 2810}%
\special{fp}%
%
\special{pn 8}%
\special{pa 2886 2610}%
\special{pa 2328 3046}%
\special{pa 4106 3082}%
\special{pa 4576 2654}%
\special{pa 4576 2654}%
\special{pa 2886 2610}%
\special{fp}%
%
\special{pn 8}%
\special{pa 4622 2504}%
\special{pa 4254 2722}%
\special{dt 0.045}%
\special{sh 1}%
\special{pa 4254 2722}%
\special{pa 4322 2704}%
\special{pa 4300 2694}%
\special{pa 4302 2670}%
\special{pa 4254 2722}%
\special{fp}%
\put(46.7800,-23.9800){\makebox(0,0)[lt]{$C$}}%
%
\special{pn 8}%
\special{pa 3790 2640}%
\special{pa 3130 2890}%
\special{pa 4184 3002}%
\special{pa 4568 2654}%
\special{pa 4568 2654}%
\special{pa 3790 2640}%
\special{fp}%
%
\special{pn 4}%
\special{pa 3168 2894}%
\special{pa 3158 2878}%
\special{fp}%
\special{pa 3206 2898}%
\special{pa 3188 2868}%
\special{fp}%
\special{pa 3244 2902}%
\special{pa 3216 2858}%
\special{fp}%
\special{pa 3282 2906}%
\special{pa 3244 2846}%
\special{fp}%
\special{pa 3320 2910}%
\special{pa 3274 2836}%
\special{fp}%
\special{pa 3358 2914}%
\special{pa 3302 2826}%
\special{fp}%
\special{pa 3396 2918}%
\special{pa 3330 2814}%
\special{fp}%
\special{pa 3434 2922}%
\special{pa 3358 2804}%
\special{fp}%
\special{pa 3472 2926}%
\special{pa 3388 2792}%
\special{fp}%
\special{pa 3508 2930}%
\special{pa 3416 2782}%
\special{fp}%
\special{pa 3548 2934}%
\special{pa 3444 2772}%
\special{fp}%
\special{pa 3586 2938}%
\special{pa 3474 2760}%
\special{fp}%
\special{pa 3624 2942}%
\special{pa 3502 2750}%
\special{fp}%
\special{pa 3660 2946}%
\special{pa 3530 2738}%
\special{fp}%
\special{pa 3700 2950}%
\special{pa 3558 2728}%
\special{fp}%
\special{pa 3738 2954}%
\special{pa 3588 2718}%
\special{fp}%
\special{pa 3774 2958}%
\special{pa 3616 2706}%
\special{fp}%
\special{pa 3812 2962}%
\special{pa 3644 2696}%
\special{fp}%
\special{pa 3850 2966}%
\special{pa 3674 2686}%
\special{fp}%
\special{pa 3888 2970}%
\special{pa 3702 2674}%
\special{fp}%
\special{pa 3926 2974}%
\special{pa 3730 2664}%
\special{fp}%
\special{pa 3964 2978}%
\special{pa 3758 2652}%
\special{fp}%
\special{pa 4002 2982}%
\special{pa 3788 2642}%
\special{fp}%
\special{pa 4040 2986}%
\special{pa 3822 2642}%
\special{fp}%
\special{pa 4078 2990}%
\special{pa 3858 2642}%
\special{fp}%
\special{pa 4116 2994}%
\special{pa 3892 2642}%
\special{fp}%
\special{pa 4154 2998}%
\special{pa 3930 2644}%
\special{fp}%
\special{pa 4188 2998}%
\special{pa 3964 2644}%
\special{fp}%
\special{pa 4210 2976}%
\special{pa 4002 2644}%
\special{fp}%
\special{pa 4234 2956}%
\special{pa 4036 2644}%
\special{fp}%
%
\special{pn 4}%
\special{pa 4256 2936}%
\special{pa 4072 2646}%
\special{fp}%
\special{pa 4278 2916}%
\special{pa 4108 2646}%
\special{fp}%
\special{pa 4302 2896}%
\special{pa 4144 2646}%
\special{fp}%
\special{pa 4322 2876}%
\special{pa 4180 2648}%
\special{fp}%
\special{pa 4346 2854}%
\special{pa 4214 2648}%
\special{fp}%
\special{pa 4368 2834}%
\special{pa 4250 2648}%
\special{fp}%
\special{pa 4390 2814}%
\special{pa 4286 2648}%
\special{fp}%
\special{pa 4414 2794}%
\special{pa 4322 2650}%
\special{fp}%
\special{pa 4436 2774}%
\special{pa 4358 2650}%
\special{fp}%
\special{pa 4458 2754}%
\special{pa 4394 2650}%
\special{fp}%
\special{pa 4482 2732}%
\special{pa 4430 2650}%
\special{fp}%
\special{pa 4504 2712}%
\special{pa 4466 2652}%
\special{fp}%
\special{pa 4526 2692}%
\special{pa 4500 2652}%
\special{fp}%
\special{pa 4548 2672}%
\special{pa 4538 2652}%
\special{fp}%
\end{picture}%
\hspace{6truecm}}

\vspace{0.5truecm}

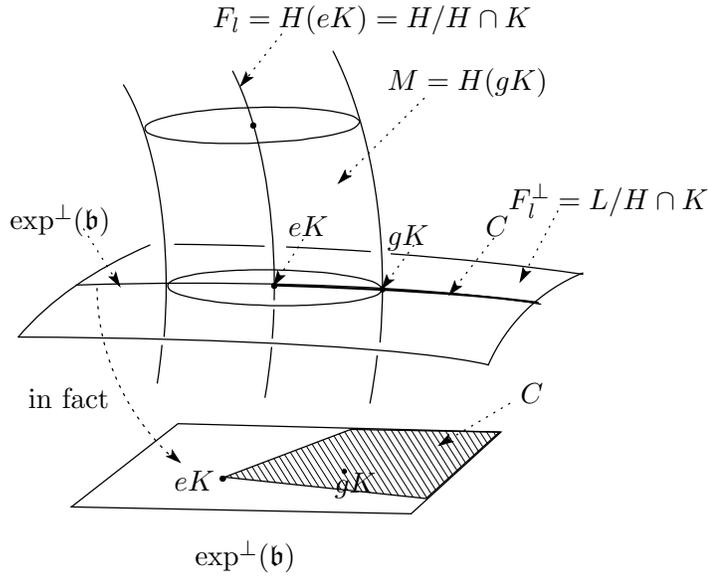
\centerline{{\bf Figure 7.}}

\vspace{1truecm}

\centerline{{\bf References}}

\vspace{0.5truecm}



\noindent
[BT] J. Berndt and H. Tamaru, Cohomogeneity one actions on noncompact 
sym-

metric spaces with a totally geodesic singular orbit, 
Tohoku Math. J. (2) {\bf 56} 

(2004), 163--177.

\noindent
{\small [BV] J. Berndt and L. Vanhecke, 
Curvature-adapted submanifolds, 
Nihonkai Math. J. {\bf 3} 

(1992), 177--185.



\noindent
[Ch] U. Christ, 
Homogeneity of equifocal submanifolds, J. Differential Geom. {\bf 62} (2002), 

1--15.

\noindent
[Co] H. S. M. Coxeter, 
Discrete groups generated by reflections, Ann. of Math. (2) 
{\bf 35} 

(1934),588--621.

\noindent
[E] H. Ewert, 
Equifocal submanifolds in Riemannian symmetric spaces, Doctoral thesis.

\noindent
[Ge1] L. Geatti, 
Invariant domains in the complexfication of a noncompact Riemannian 
sym-

metric space, J. Algebra {\bf 251} (2002), 619--685.

\noindent
[Ge2] L. Geatti, 
Complex extensions of semisimple symmetric spaces, manuscripta math. 

{\bf 120} (2006) 1-25.

\noindent
[Go] C. Gorodski, 
Polar actions on compact symmetric spaces which admit a totally geodesic 

principal orbit, Geom. Dedicata {\bf 103} (2004), 193--204.

\noindent
[HL] E. Heintze and X. Liu, 
Homogeneity of infinite dimensional isoparametric submanifolds, 

Ann. of Math. {\bf 149} (1999), 149-181.

\noindent
[HLO] E. Heintze, X. Liu and C. Olmos, 
Isoparametric submanifolds and a Chevalley type 

restricction theorem, Integrable systems, geometry, and topology, 151-190, 
AMS/IP 

Stud. Adv. Math. 36, Amer. Math. Soc., Providence, RI, 2006.



\noindent
[HPTT] E. Heintze, R. S. Palais, C. L. Terng and G. Thorbergsson, 
Hyperpolar actions on 

symmetric spaces, Geometry, topology and physics, 214--245 Conf. Proc. 
Lecture Notes 

Geom. Topology {\bf 4}, Internat. Press, Cambridge, Mass., 1995.

\noindent
[He] S. Helgason, 
Differential geometry, Lie groups and symmetric spaces, Pure Appl. Math. 

80, Academic Press, New York, 1978.

\noindent
[Hug] M. C. Hughes, 
Complex reflection groups, Comm. Algebra {\bf 18} (1990), 3999--4029.

\noindent
[Hui1] G. Huisken, 
Flow by mean curvature of convex surfaces into spheres, J. Differential 

Geom. {\bf 20} (1984) 237-266.

\noindent
[Hui2] G. Huisken, 
Contracting convex hypersurfaces in Riemannian manifolds by their mean 

curvature, Invent. math. {\bf 84} (1986) 463-480.

\noindent
[KN] S. Kobayashi and K. Nomizu, 
Foundations of differential geometry, Interscience Tracts 

in Pure and Applied Mathematics 15, Vol. II, New York, 1969.

\noindent
[Koi1] N. Koike, 
Submanifold geometries in a symmetric space of non-compact 
type and a 

pseudo-Hilbert space, Kyushu J. Math. {\bf 58} (2004), 167--202.

\noindent
[Koi2] N. Koike, 
Complex equifocal submanifolds and infinite dimensional 
anti-Kaehlerian 

isoparametric submanifolds, Tokyo J. Math. {\bf 28} (2005), 
201--247.

\noindent
[Koi3] N. Koike, 
Actions of Hermann type and proper complex equifocal submanifolds, 
Osaka 

J. Math. {\bf 42} (2005) 599-611.

\noindent
[Koi4] N. Koike, 
A splitting theorem for proper complex equifocal submanifolds, Tohoku 

Math. J. {\bf 58} (2006) 393-417.

\noindent
[Koi5] N. Koike, 
Complex hyperpolar actions with a totally geodesic orbit, Osaka J. Math. 

{\bf 44} (2007), 491-503.  

\noindent
[Koi6] N. Koike, On curvature-adapted and proper complex equifocal 
submanifolds, Kyung

pook Math. J. {\bf 50} (2010), 509-536.

\noindent
[Koi7] N. Koike, The homogeneous slice theorem for the complete 
complexification of a 

proper complex equifocal submanifold, Tokyo J. Math. {\bf 33} (2010), 1-30.

\noindent
[Koi8] N. Koike, Hermann type actions on a pseudo-Riemannian symmetric 
space, Tsukuba 

J. Math. {\bf 34} (2010), 137-172.

\noindent
[Koi9] N. Koike, Examples of a complex hyperpolar action without singular 
orbit, Cubo A 

Math. J. {\bf 12} (2010), 131-147.

\noindent
[Koi10] N. Koike, Collapse of the mean curvature flow for equifocal 
submanifolds, Asian J. 

Math. {\bf 15} (2011), 101-128.

\noindent
[Koi11] N. Koike, The complexifications of pseudo-Riemannian manifolds and 
anti-Kaehler 

geometry, arXiv:math.DG/0807.1601.v3.

\noindent
[Koi12] N. Koike, A Cartan type identity for isoparametric hypersurfaces 
in symmetric 

spaces, Tohoku Math. J. (to appear).

\noindent
[Koi13] N. Koike, The classifications of certain kind of isoparametric 
submanifolds in 

non-compact symmetric spaces, arXiv:math.DG/1209.1933v1.




\noindent
[Kol] A. Kollross, A classification of hyperpolar and cohomogeneity one 
actions, Trans. 

Amer. Math. Soc. {\bf 354} (2002), 571--612.

\noindent
[LT] X. Liu and C. L. Terng, The mean curvature flow for isoparametric 
submanifolds, 

Duke Math. J. {\bf 147} (2009), 157-179.

\noindent
[O'N] B. O'Neill, 
Semi-Riemannian Geometry, with applications to relativity, 
Pure Appl. 

Math. 103, Academic Press, New York, 1983.




\noindent
[P] R. S. Palais, 
Morse theory on Hilbert manifolds, Topology {\bf 2} (1963), 299--340.



\noindent
[PaTe] R. S. Palais and C. L. Terng, Critical point theory and submanifold 
geometry, Lecture 

Notes in Math. {\bf 1353}, Springer-Verlag, Berlin, 1988.

\noindent
[PoTh] F. Podest$\acute a$ and G. Thorbergsson, Polar actions on rank-one 
symmetric spaces, 

J. Differential Geometry {\bf 53} (1999) 131-175.

\noindent
[S] G. Schwarz, 
Smooth functions invariant under the action of a compact Lie group, 

Topology {\bf 14} (1975) 63-68.



\noindent
[Te1] C. L. Terng, 
Isoparametric submanifolds and their Coxeter groups, 
J. Differential 

Geom. {\bf 21} (1985), 79--107.

\noindent
[Te2] C. L. Terng, 
Proper Fredholm submanifolds of Hilbert space, 
J. Differential Geom. {\bf 29} 

(1989), 9--47.

\noindent
[Te3] C. L. Terng, 
Polar actions on Hilbert space, J. Geom. Anal. {\bf 5} (1995), 
129--150.

\noindent
[TT] C. L. Terng and G. Thorbergsson, 
Submanifold geometry in symmetric spaces, 

J. Differential Geom. {\bf 42} (1995), 665--718.

\noindent
[Z] X. P. Zhu, Lectures on mean curvature flows, Studies in Advanced Math., 
AMS/IP, 2002.






\begin{flushleft}
\textsc{Department of Mathematics, Faculty of Science, 
Tokyo University of Science, 1-3 Kagurazaka Shinjuku-ku, Tokyo 162-8601, 
Japan}\\
\textit{E-mail address}: koike@ma.kagu.tus.ac.jp
\end{flushleft}

\end{document}